\documentclass[10pt]{article}
\usepackage[utf8]{inputenc}
\usepackage[T1]{fontenc}

\usepackage{amsfonts, amsbsy, amssymb, amsmath, graphicx, float, subfigure, authblk}
\usepackage{setspace, ifdraft}

\usepackage[osf,sc]{mathpazo}

\usepackage[scaled=0.90]{helvet}

\usepackage[scaled=0.85]{beramono}

\ifdraft{\onehalfspace}{\singlespace} 

\usepackage{hyperref}
\hypersetup{
	colorlinks,
	linkcolor={red!50!black},
	citecolor={blue!50!black},
	urlcolor={blue!80!black}
}

\usepackage{xargs}                      
\usepackage[pdftex,dvipsnames]{xcolor}  
\usepackage[colorinlistoftodos,prependcaption,textsize=small]{todonotes}
\newcommandx{\unsure}[2][1=]{\todo[linecolor=red,backgroundcolor=red!25,bordercolor=red,#1]{#2}}
\newcommandx{\change}[2][1=]{\todo[linecolor=blue,backgroundcolor=blue!25,bordercolor=blue,#1]{#2}}
\newcommandx{\info}[2][1=]{\todo[linecolor=OliveGreen,backgroundcolor=OliveGreen!25,bordercolor=OliveGreen,#1]{#2}}
\newcommandx{\improvement}[2][1=]{\todo[linecolor=Plum,backgroundcolor=Plum!25,bordercolor=Plum,#1]{#2}}
\newcommandx{\thiswillnotshow}[2][1=]{\todo[disable,#1]{#2}}

\usepackage[square, sort, numbers]{natbib}
\usepackage{enumitem}
\usepackage{wasysym, tikz}

\ifdraft{
\usepackage[mathlines]{lineno} 
\linenumbers\relax
\newcommand*\patchAmsMathEnvironmentForLineno[1]{%
	\expandafter\let\csname old#1\expandafter\endcsname\csname #1\endcsname
	\expandafter\let\csname oldend#1\expandafter\endcsname\csname 
	end#1\endcsname
	\renewenvironment{#1}%
	{\linenomath\csname old#1\endcsname}%
	{\csname oldend#1\endcsname\endlinenomath}}%
\newcommand*\patchBothAmsMathEnvironmentsForLineno[1]{%
	\patchAmsMathEnvironmentForLineno{#1}%
	\patchAmsMathEnvironmentForLineno{#1*}}%
\AtBeginDocument{%
	\patchBothAmsMathEnvironmentsForLineno{equation}%
	\patchBothAmsMathEnvironmentsForLineno{align}%
	\patchBothAmsMathEnvironmentsForLineno{flalign}%
	\patchBothAmsMathEnvironmentsForLineno{alignat}%
	\patchBothAmsMathEnvironmentsForLineno{gather}%
	\patchBothAmsMathEnvironmentsForLineno{multline}}%
}

\newcommand{\newpoint}{\tikz\draw[black,fill=black] (0,0) circle (.4ex);}

\begin{document}

\title{Finding NHIM: Identifying High Dimensional Phase Space Structures in Reaction Dynamics using Lagrangian Descriptors}

\author[1]{Shibabrat Naik\thanks{s.naik@bristol.ac.uk}}
\author[2]{V{\'i}ctor J. Garc{\'i}a-Garrido}
\author[1]{Stephen Wiggins}
\affil[1]{\small{School of Mathematics, University of Bristol, Bristol BS8 1TW, United Kingdom}}
\affil[2]{\small{Departamento de F{\'i}sica y Matem{\'a}ticas, Universidad de Alcal{\'a}, 28871, Alcal{\'a} de Henares, Spain}}

\date{}

\maketitle

\begin{abstract}
Phase space structures such as dividing surfaces, normally hyperbolic invariant manifolds, their stable and unstable manifolds have been an integral part of computing quantitative results such as transition fraction, stability erosion in multi-stable mechanical systems, and reaction rates in chemical reaction dynamics. Thus, methods that can reveal their geometry in high dimensional phase space (4 or more dimensions) need to be benchmarked by comparing with known results. In this study we assess the capability of one such method called {\em Lagrangian descriptor} for revealing the types of high dimensional phase space structures associated with index-1 saddle in Hamiltonian systems. The Lagrangian descriptor based approach is applied to two and three degree-of-freedom quadratic Hamiltonian systems where the high dimensional phase space structures are known, that is as closed-form analytical expressions. This leads to a direct comparison of features in the Lagrangian descriptor plots and the phase space structures' intersection with an isoenergetic two-dimensional surface, and hence provides a validation of the approach.


\end{abstract}


\ifdraft{\tableofcontents}{}

\section{Introduction}
	
Identifying invariant manifolds in high dimensional phase space and their anchor, normally hyperbolic invariant manifold (NHIM), is the stepping stone in applying phase space transport methods to a broad array of problems in physics, chemistry, and engineering~\cite{wiggins90,Gottwald1995,Jaffe1999,DeOliveira2002,Dellnitz2005,Gabern2005,Gabern2006,waalkens08,ww10,wiggins2013normally,wig2016}. These problems are typically formulated in phase space of more than 4 dimensions (that is, $N = 2$ or more degrees of freedom) and the geometric 
approach of computing invariant manifolds that are codimension-1 separatrices requires Poincar\'e-Birkoff normal form theory. The Poincar\'e-Birkoff normal form approach for computing NHIMs and their stable and unstable manifolds has only been developed for Hamiltonian systems in the neighborhood of  index-k saddle points. This approach in itself is successful and rigorous, but its implementation requires experience in writing or using algebraic manipulation programs. 
Furthermore, these computations become expensive as the dimensions increase since Poincar\'e-Birkoff normal form theory is based on Taylor expansion, a high dimensional vector valued polynomial, of the Hamiltonian in a neighborhood of the saddle point to sufficiently high order. Besides, the resulting visualization can be non-obvious leading to cumbersome interpretation and can not work for high dimensions. The succcess and shortcoming of the vizualization approach based on topological methods for four-dimensional space becomes apparent when trying to extend the results in Refs.~\cite{Kocak1986,Banks1992,Hoffmann1991}. This raises the question of detecting the signatures of the high dimensional structures, if they exist, when probed with low dimensional surfaces. If one does manage to detect these structures, how do they manifest on the low dimensional surfaces.

One such method is the Lagrangian descriptors (LDs) developed in Refs.~\cite{madrid2009,mendoza2010,mancho2013,lopesino2017}. Lagrangian descriptors are a class of trajectory diagnostic methods that can reveal phase space structures by encoding geometric property of trajectories (such as, phase space arc length, configuration space distance or displacement, cumulative action or kinetic energy) initialised on a two dimensional surface. The method was originally developed in the context of Lagrangian transport studies in fluid dynamics that requires identifying transport barriers which are the invariant manifolds in time dependent 2D fluid flow. Besides, the method applies to both Hamiltonian and non-Hamiltonian systems (\cite{lopesino2017}) as well as to systems with arbitrary, stochastic and dissipative, time-dependence (\cite{balibrea2016lagrangian,craven2017lagrangian,junginger2016lagrangian}). Lagrangian descriptor based detection of transport barriers has also been applied directly to data sets, such as those obtained from satellite observations or numerical simulations~(\cite{amism11,mendoza2014,ggmwm15,ramos2018}).
Furthermore, the method is straightforward to implement computationally and it provides a ``high resolution'' method for exploring the influence of high dimensional phase space structure on trajectory behaviour.  The method of Lagrangian descriptors takes an {\it opposite} approach to that of Lyapunov exponent type calculations by emphasizing the initial conditions of trajectories, rather than their advected locations that is involved in calculating normalized rate of divergence. This is achieved by considering a two dimensional section of the full phase space and discretizing with a dense grid of initial conditions. Even though the trajectories wander off in the phase space, as the initial conditions evolve in time, there is no loss in resolution of the two dimensional section. Our objective is to clarify the use of Lagrangian descriptors as a diagnostic on two dimensional sections of high dimensional phase space structures. This diagnostic is also meant to be used as the preliminary step in computing the NHIM, their stable and unstable manifolds using other computational means~\cite{junginger2016transition,bardakcioglu2018,ezra_2018}. In this article, we assess the capability of Lagrangian descriptors to detect the type of high dimensional phase space structures in Hamiltonian systems such as the NHIM, their stable, and unstable manifolds that are used in computing rates of chemical reactions. 

Computing chemical reaction rates is based on the fundamental framework of transition state theory as formulated in phase space by Polanyi, Evans, Wigner, Eyring~\cite{evans1935applications,wigner1938transition,eyring1938theory,garrett_2000}). Further research has shown that constructing a \emph{locally} recrosssing free, minimal flux orientable surface called a {\it dividing surface} (DS) is the phase space structure that provides the correct estimate for reaction rates. This dividing surface of geometry $\mathbb{S}^{2N - 2}$ ( a 2N-2 dimensional sphere) is constructed from the NHIM of geometry $\mathbb{S}^{2N - 3}$ which forms the equator of the dividing surface, on a constant energy surface~\cite{waalkens2004direct,waalkens08,ww10,wig2016}. Furthermore, the global dynamics of reactive and non-reactive trajectories is governed by the invariant (stable and unstable) manifolds, of the NHIM, are $\mathbb{S}^{2N - 3} \times \mathbb{R}$, and thus act as codimension-1 separatrices on the energy surface. Thus, the NHIM acts as the anchor for the local dynamics via the dividing surface from reactants to products or vice-versa and for the global dynamics via the stable/unstable invariant manifolds. Thus, detecting and constructing the NHIM forms a fundamental step in obtaining quantitative results in reaction dynamics~\cite{ezra_2009,ezra_2018}. In the present study, we will focus on detecting and verifying these high dimensional phase space structures. 


Recently the applicability of the Lagrangian descriptor based approach to time dependent problems with random and dissipative forcing in chemical reaction dynamics has been shown whereby the transition state trajectory is calculated using the extremal values in Lagrangian descriptor values; see Refs.~\cite{craven2015lagrangian,craven2016deconstructing,craven2017lagrangian,junginger2016uncovering,junginger2016transition,junginger2017chemical,junginger2017variational,feldmaier2017obtaining,revuelta2017transition,patra2018detecting}. The initial conditions for the transition state trajectory is identified by computing the extrema of the Lagrangian descriptor on a two dimensional domain. Comparing extremal and singular features in the Lagrangian descriptor plot with invariant manifolds in three dimensional vector fields, Ref.~\cite{garcia-garrido_2018} has also provided numerical evidence for detecting NHIM using the LD on two dimensional surfaces. The authors presented the comparison using a stacked version of a two dimensional linear saddle, the duffing oscillator, and a time-perturbed 3D geophysical model. This has the advantage that the NHIM and the invariant manifolds are known exactly or can be computed with other established methods. For the problems considered therein, the NHIM (a curve embedded in 3D) and its associated invariant manifolds (a 2D surface in 3D) can be visualized in the 3D space since the dimensionality of phase space structures is less than 3. We take this approach to the next logical step by applying this to benchmark problems of two and three degrees-of-freedom Hamiltonian system where the phase space structures are known exactly, that is they can be written as closed-form analytical expressions

This article is outlined as follows. In Section~\ref{sec:theory}, we describe the method of Lagrangian descriptor used in this study and present an analytical result on identifying invariant manifolds using features in the Lagrangian descriptor. In Section~\ref{sec:normal_form_2dof} and~\ref{sec:normal_form_3dof}, we discuss the benchmark systems and present numerical evidence of our claim using isoenergetic two-dimensional surfaces. In Section~\ref{sec:summ}, we summarize our results on detecting invariant manifolds and discuss related future directions. 

\section{Theory}
\label{sec:theory}

\subsection{Method of Lagrangian Descriptor}
\label{sec:LD}

The Lagrangian descriptor (LD) as presented in Ref.\cite{madrid2009} is the arc length of a trajectory calculated on a chosen initial time $t_0$ and measured for fixed forward and backward integration time, $\tau$. For continuous time dynamical systems, Ref.\cite{lopesino2017} gives an alternative definition of the LD which is useful for proving rigorous results and can be computed along with the trajectory. It provides a characterization of the notion of singular features of the LD that facilitates a proof for detecting invariant manifolds in certain model situations.  In addition, the ``additive nature'' of this new definition of LD provides 
an approach for assessing the influence of each degree-of-freedom separately on the Lagrangian descriptor.  This property was used in Ref.\cite{demian2017} which showed that a Lagrangian descriptor can be used to detect Lyapunov periodic orbits in the two degrees-of-freedom H{\'e}non-Heiles Hamiltonian system. We will describe this procedure for two and three degrees-of-freedom linear autonomous Hamiltonian systems. We begin by establishing notation in the general setting of a time-dependent vector field where 

\begin{equation}
  \frac{d\mathbf{x}}{dt} = \mathbf{v}(\mathbf{x},t), \quad \mathbf{x} \in \mathbb{R}^n \;,\; t \in \mathbb{R}
\end{equation}

where $\mathbf{v}(\mathbf{x},t) \in C^r$ ($r \geq 1$) in $\mathbf{x}$ and continuous in time. The definition of LDs depends on the initial condition $\mathbf{x}_{0} = \mathbf{x}(t_0)$, on the initial time $t_0$ (trivial for autonomous systems) and the integration time $\tau$, and the type of norm of the trajectory's components, and takes the form,

\begin{equation}
M_p(\mathbf{x}_{0},t_0,\tau) = \displaystyle{\int^{t_0+\tau}_{t_0-\tau} \sum_{i=1}^{n} |\dot{x}_{i}(t;\mathbf{x}_{0})|^p \; dt}
\label{eqn:M_function}
\end{equation}

\noindent where $p \in (0,1]$ and $\tau \in \mathbb{R}^{+}$ are freely chosen parameters,  and the overdot symbol represents the derivative with respect to time. It is to be noted here that there are three formulations of the function $M_p$ in the literature: the arc length of a trajectory in phase space~\cite{madrid2009}, the arc length of a trajectory projected on the configuration space~~\cite{junginger2016transition,junginger2016uncovering,junginger2017chemical,junginger2017variational}, and the sum of the $p$-norm of the vector field components~\cite{lopesino_2015,lopesino2017}.
Although the latter formulation of the Lagrangian descriptor~\eqref{eqn:M_function} developed in Ref.~\cite{lopesino_2015,lopesino2017} does not resemble the arc length, the numerical results using either of these forms have been shown to be in agreement and promise of predictive capability in geophysical flows~\cite{amism11,mendoza2014,ggmwm15,ramos2018}. The formulation we adopt here is motivated by the fact that this allows for proving rigorous result, which we will discuss in the next section, connecting the singular features and minimum in the LD plots with NHIM and its stable and unstable manifolds. 
It follows from the result that 

\begin{align}
\mathcal{W}^s(\mathbf{x}_0, t_0) & = \text{\rm argmin} \; \mathcal{L}^{(f)}(\mathbf{x}_0, t_0, \tau) \\
\mathcal{W}^u(\mathbf{x}_0, t_0) & = \text{\rm argmin} \; \mathcal{L}^{(b)}(\mathbf{x}_0, t_0, \tau)
\end{align}
where the stable and unstable manifolds ($\mathcal{W}^s(\mathbf{x}_0, t_0)$ and $\mathcal{W}^u(\mathbf{x}_0, t_0)$) denote the invariant manifolds at intial time $t_0$ and $\text{\rm argmin} \; (\cdot)$ denotes the argument that minimizes the function $\mathcal{L}^{(\cdot)}(\mathbf{x}_0, t_0, \tau)$ in forward and backward time, respectively. In addition, the coordinates of the NHIM at time $t_0$ is given by the intersection $\mathcal{W}^s(\mathbf{x}_0, t_0)$ and $\mathcal{W}^u(\mathbf{x}_0, t_0)$ of the stable and unstable manifolds, and thus given by

\begin{align}
\mathcal{M}(\mathbf{x}_0, t_0) & = \text{\rm argmin} \; \left( \mathcal{L}^{(f)}(\mathbf{x}_0, t_0, \tau) + \mathcal{L}^{(b)}(\mathbf{x}_0, t_0, \tau) \right) = \text{\rm argmin} \; \mathcal{L}(\mathbf{x}_0, t_0, \tau) \qquad \text{NHIM}
\end{align}
%

\subsection{Lagrangian descriptor and invariant manifolds}

We dedicate this section to prove that the method of Lagrangian descriptor recovers the NHIM and its stable and unstable manifolds for a linear quadratic Hamiltonian system with 3 DoF. We will show that given a sufficiently large integration time, the Lagrangian descriptor value reaches a minimum value at the NHIM and its stable and unstable manifolds, and also that singularities (singular features) in the Lagrangian descriptor values, that is points where the scalar function is non-differentiable, identify the stable and unstable manifolds of the NHIM, and hence also the NHIM. Therefore, this method can be used to detect the NHIM and thus a way to construct the dividing surface that is essential for the computation of reaction rates in chemistry. The arguments we follow to prove this mathematical connection are based on those described in Refs.~\cite{lopesino2017,demian2017}.

Consider the quadratic Hamiltonian with $3$ DoF that has a index-1 saddle at the origin
\begin{equation}
\mathcal{H}(\mathbf{q},\mathbf{p}) = \dfrac{\lambda}{2}\left(p_1^2 - q_1^2\right) +  \dfrac{\omega_2}{2} \left(p_2^2 + q_2^2 \right) + \dfrac{\omega_3}{2} \left(p_3^2 + q_3^2 \right)
\label{ham_ndof_lin}
\end{equation}
Observe that all the DoF that appear in the Hamiltonian $\mathcal{H}$ are uncoupled so we can write:
\[
\mathcal{H} = \sum_{i=1}^{3} \mathcal{H}_i
\]  
where the energy in each DoF is given by
\[
\mathcal{H}_i = \dfrac{\omega_i}{2} \left(p_i^2 + q_i^2 \right) \;,\quad i \in \lbrace 2,3\rbrace \quad,\quad \mathcal{H}_1 = \dfrac{\lambda}{2}\left(p_1^2 - q_1^2\right)
\]

Given the initial condition $\mathbf{x}_0 = \mathbf{x}(t_0) = (\mathbf{q}_0,\mathbf{p}_0) \in \mathbb{R}^6$ at time $t_0 = 0$, where $\mathbf{q}_0 = (q_1^{0},q_2^{0},q_3^{0})$ and $\mathbf{p}_0 = (p_1^{0},p_2^{0},p_3^{0})$, the general solution to Hamilton's equations is
\begin{equation}
\begin{split}
q_1(t) & = \frac{1}{2} \left[A e^{\lambda t} + B e^{-\lambda t}\right] \\[.1cm]
p_1(t) & = \frac{1}{2} \left[A e^{\lambda t} - B e^{-\lambda t}\right] \\[.1cm]
q_i(t) & = q^0_i \cos(\omega_i t) + p^0_i \sin(\omega_i t) \; , \; i \in \lbrace 2,3 \rbrace \\[.2cm]
p_i(t) & = p^0_i \cos(\omega_i t) - q^0_i \sin(\omega_i t) \; , \; i \in \lbrace 2,3\rbrace 
\end{split}
\label{gen_sol_linham}
\end{equation}
where the constants $A$ and $B$ are given by
\[
A = q^0_1 + p_1^0 \quad,\quad B = q^0_1 - p_1^0
\]
This initial condition gives rise to a phase space trajectory with total energy
\[
\mathcal{H}\left(\mathbf{x}_0\right) = \mathcal{H}_0 = \sum_{i=1}^{3} \mathcal{H}^0_i
\]
which is distributed among all the DoF of the system, and the energy in each DoF is
\[
\mathcal{H}_i^{0} = \dfrac{1}{2} \left(p_i^2 +  \omega_i^2 q_i^2 \right) \;,\quad i \in \lbrace 1,\ldots,n-1 \rbrace \quad,\quad \mathcal{H}_n^{0} = \dfrac{1}{2}\left(p_n^2 - \lambda^2 q_n^2\right)
\]
Observe that each bath mode $(q_i,p_i)$ will give rise to a periodic orbit in the $q_i$-$p_i$ plane, where $i \in \lbrace 2,3 \rbrace$. Moreover, this periodic orbit is a cricle
\[
q_i^2 + p_i^2 = \dfrac{2\mathcal{H}^{0}_i}{\omega_i}
\]
with radius $R = \sqrt{\dfrac{2\mathcal{H}^{0}_i}{\omega_i}}$. The isoenergetic NHIM for this Hamiltonian system is:
\[
\text{NHIM} = \lbrace \left(\mathbf{q},\mathbf{p}\right) \in \mathbb{R}^{6} \;|\; \mathcal{H}_0 =  \sum_{i=2}^{3} \dfrac{\omega_i}{2} \left( p_i^2 + q_i^2 \right) \; , \; q_1 = p_1 = 0 \rbrace 
\]
which is a $3$-dimensional sphere. The stable and unstable manifolds of the NHIM are
\begin{equation}
\begin{split}
\mathcal{W}^{u}\left(\text{NHIM}\right) & = \lbrace \left(\mathbf{q},\mathbf{p}\right) \in \mathbb{R}^{6} \;|\; \mathcal{H}_0 =  \sum_{i=2}^{3} \dfrac{\omega_i}{2} \left( p_i^2 + q_i^2 \right) \; , \; q_1 = p_1 \rbrace \\
\mathcal{W}^{s}\left(\text{NHIM}\right) & = \lbrace \left(\mathbf{q},\mathbf{p}\right) \in \mathbb{R}^{6} \;|\; \mathcal{H}_0 =  \sum_{i=2}^{3} \dfrac{\omega_i}{2} \left( p_i^2 + q_i^2 \right) \; , \; q_1 = -p_1 \rbrace 
\end{split}
\label{su_manifolds}
\end{equation}  

Take a fixed integration time $\tau > 0$ and $\gamma \in (0,1]$, The Lagrangian descriptor is
\begin{equation}
M_{\gamma}(\mathbf{x}_{0},t_0,\tau) = \int^{t_0+\tau}_{t_0-\tau} \sum_{i=1}^{3}  |\dot{q}_{i}(t;\mathbf{x}_0)|^\gamma + |\dot{p}_{i}(t;\mathbf{x}_0)|^\gamma \; dt
\label{Mp_func}
\end{equation}
Observe that we can decompose this integral into the hyperbolic and elliptic components
\begin{equation}
M_{\gamma}(\mathbf{x}_{0},t_0,\tau) = M^{h}_{\gamma}(\mathbf{x}_{0},t_0,\tau) + M^{e}_{\gamma}(\mathbf{x}_{0},t_0,\tau)
\label{hyp_comp_M}
\end{equation}
where the hyperbolic part is
\begin{equation}
M^{h}_{\gamma}(\mathbf{x}_{0},t_0,\tau) = \int^{t_0+\tau}_{t_0-\tau} |\dot{q}_{1}(t;\mathbf{x}_0)|^\gamma + |\dot{p}_{1}(t;\mathbf{x}_0)|^\gamma \; dt
\label{ell_comp_M}
\end{equation}
and the elliptic component is given by
\[
M^{e}_{\gamma}(\mathbf{x}_{0},t_0,\tau) = \int^{t_0+\tau}_{t_0-\tau} \; \sum_{i=2}^{3} |\dot{q}_{i}(t;\mathbf{x}_0)|^\gamma + |\dot{p}_{i}(t;\mathbf{x}_0)|^\gamma \; dt
\]

We focus our attention first on the hyperbolic contribution of LDs. Since the dynamical system is autonomous we can take without loss of generality $t_0 = 0$, and we eliminitae the dependence of LDs on the initial time to simplify notation. We can write
\[
M_{\gamma}^{h}\left(\mathbf{x}_0,\tau\right) = M_{\gamma}^{h,q_1}\left(\mathbf{x}_0,\tau\right) + M_{\gamma}^{h,p_1}\left(\mathbf{x}_0,\tau\right) = \int_{-\tau}^{\tau} |\dot{q}_1|^{\gamma} \, dt + \int_{-\tau}^{\tau} |\dot{p}_1|^{\gamma} \, dt
\]
for which we need the time derivative of the position and momentum coordinates
\[
\dot{q}_1(t) = \frac{\lambda}{2} \left[A e^{\lambda t} - B e^{-\lambda t}\right] \quad,\quad \dot{p}_1(t) = \frac{\lambda}{2} \left[A e^{\lambda t} + B e^{-\lambda t}\right]
\]
Observe that it is not possible to calculate the integrals $M_{\gamma}^{h,q_n}$ and $M_{\gamma}^{h,p_n}$ analytically, so we will accurately approximate their values by means of an asymptotic analysis. Take a small value of $\tau_0$, then we can write
\begin{equation}
\begin{split}
M^{h,q_1}_{\gamma}(\mathbf{x}_{0},\tau) & = \int^{-\tau_0}_{-\tau} |\dot{q}_1|^\gamma \; dt + M^{h,q_1}_{\gamma}(\mathbf{x}_{0},\tau_0) + \int^{\tau}_{\tau_0} |\dot{q}_1|^\gamma \; dt \\
M^{h,p_1}_{\gamma}(\mathbf{x}_{0},\tau) & = \int^{-\tau_0}_{-\tau} |\dot{p}_1|^\gamma \; dt +  M^{h,p_1}_{\gamma}(\mathbf{x}_{0},\tau_0) + \int^{\tau}_{\tau_0} |\dot{p}_1|^\gamma \; dt 
\end{split}
\label{tau_decomp}
\end{equation}
Expanding $M^{h,q_1}_{\gamma}(\mathbf{x}_{0},\tau_0)$ and  $M^{h,p_1}_{\gamma}(\mathbf{x}_{0},\tau_0)$ in a Taylor series about $\tau = 0$ gives
\[
\begin{split}
M^{h,q_1}_{\gamma}(\mathbf{x}_{0},\tau_0) & = M^{h,q_1}_{\gamma}(\mathbf{x}_{0},0) + \tau_0\frac{\partial M^{h,q_1}_{\gamma}}{\partial \tau_0}(\mathbf{x}_{0},0) + O(\tau^2_0) = \\[.1cm] & = 2\tau_0 \lambda^{\gamma} \, |p^0_1|^{\gamma} + O(\tau_0^2) \approx 2\tau_0 \lambda^{\gamma} \, |p^0_1|^{\gamma} \\[.2cm]
M^{h,p_1}_{\gamma}(\mathbf{x}_{0},\tau_0) & = M^{h,p_1}_{\gamma}(\mathbf{x}_{0},0) + \tau_0\frac{\partial M^{h,p_1}_{\gamma}}{\partial \tau_0}(\mathbf{x}_{0},0) + O(\tau^2_0) = \\[.1cm] & = 2\tau_0 \lambda^{\gamma} \, |q^0_1|^{\gamma} + O(\tau_0^2) \approx 2\tau_0 \lambda^{\gamma} \, |q^0_1|^{\gamma}
\end{split}
\]
This clearly shows that the hyperobolic component of LDs behaves for small $\tau_0$ as
\[
M_{\gamma}^{h}\left(\mathbf{x}_0,\tau_0\right) \approx 2\tau_0 \lambda^{\gamma} \left(|p^0_1|^\gamma + |q^0_1|^\gamma \right)  
\]
having singularities on the phase space sets $q_n^0 = 0$ or $p_n^0 = 0$ which do not coincide with the unstable and stable manifolds of the NHIM for the linear Hamiltonian. Nevertheless, as we describe next, when the integration time $\tau$ is sufficiently large, the leading order singularities that appear in LDs align with the true stable and unstalble manifolds of the NHIM. 

Take $\tau \gg 1$, we deal first with the integral that involves the time derivative of $q_1$
\begin{equation}
\begin{split}
\int^{\tau}_{\tau_0} |\dot{q}_1|^\gamma \; dt &= \dfrac{\lambda^{\gamma}}{2^{\gamma}} \int^{\tau}_{\tau_0} | A e^{\lambda t} - B e^{-\lambda t}|^{\gamma} \, dt = \dfrac{\lambda^{\gamma}}{2^{\gamma}} \int^{\tau}_{\tau_0} \left| A e^{\lambda t} \left(1 - \frac{B}{A}e^{-2\lambda t}\right)\right|^{\gamma} \, dt = \\[.2cm] &= \dfrac{\lambda^{\gamma}}{2^{\gamma}} \int^{\tau}_{\tau_0} |Ae^{\lambda t}|^{\gamma} \left|1 - \frac{B}{A}e^{-2\lambda t}\right|^{\gamma} dt = \dfrac{\lambda^{\gamma}}{2^{\gamma}} \int^{\tau}_{\tau_0} |Ae^{\lambda t}|^{\gamma} + O\left(\frac{|B|}{|A|^{1-p}}e^{-\gamma \lambda t}\right) \, dt
\end{split}
\end{equation}
where we have used the expansion $(1+x)^{p} = 1 + px^{p-1} + O(x^2)$ for small $x$. Then the asymptotic behavior of the integral is
\[
\int^{\tau}_{\tau_0} |\dot{q}_1|^\gamma \; dt \sim \dfrac{\lambda^{\gamma}}{2^{\gamma}} \int^{\tau}_{\tau_0} |Ae^{\lambda t}|^{\gamma} \, dt = \frac{|A|^{\gamma}\lambda^{\gamma}}{2^{\gamma}} \int^{\tau}_{\tau_0} e^{\gamma\lambda t} \, dt \sim \dfrac{|A|^{\gamma}\lambda^{\gamma-1}}{\gamma 2^{\gamma}} \, e^{\gamma \lambda \tau}
\]
Analogously, for negative times
\begin{equation}
\begin{split}
\int^{-\tau_0}_{-\tau} |\dot{q}_1|^\gamma \; dt &= \frac{\lambda^{\gamma}}{2^{\gamma}} \int^{-\tau_0}_{-\tau} |A e^{\lambda t} - B e^{-\lambda t}|^{\gamma} \, dt = \frac{\lambda^{\gamma}}{2^{\gamma}} \int^{\tau}_{\tau_0} |Be^{\lambda s}|^{\gamma} \left|1 - \frac{A}{B}e^{-2\lambda s}\right|^{\gamma} \, ds \\[.2cm] &= \frac{\lambda^{\gamma}}{2^{\gamma}} \int^{\tau}_{\tau_0} |Be^{\lambda s}|^{\gamma} + O\left(\frac{|A|}{|B|^{1-p}}e^{-\gamma \lambda s}\right) \, ds
\end{split}
\end{equation}
Then the asymptotic behavior of the integral is
\[
\int^{-\tau_0}_{-\tau} |\dot{q}_1|^\gamma \; dt \sim \frac{\lambda^{\gamma}}{2^{\gamma}} \int^{\tau}_{\tau_0} |B e^{\lambda s}|^{\gamma} ds = \frac{|B|^{\gamma} \lambda^{\gamma}}{2^{\gamma}} \int^{\tau}_{\tau_0} e^{\gamma\lambda s} ds \sim \dfrac{|B|^{\gamma} \lambda^{\gamma-1}}{\gamma 2^{\gamma}} \, e^{\gamma \lambda \tau}
\]
If we develop a similar asymptotic argument for the integral that involves the time derivative of $p_1$ it can be easily shown that
\[
\begin{split}
\int^{\tau}_{\tau_0} |\dot{p}_1|^\gamma \, dt &\sim \frac{\lambda^{\gamma}}{2^{\gamma}} \int^{\tau}_{\tau_0} |Ae^{\lambda t}|^{\gamma} \, dt \sim \dfrac{|A|^{\gamma} \lambda^{\gamma-1}}{\gamma 2^{\gamma}} \, e^{\gamma \lambda \tau} \\[.2cm]
\int^{-\tau_0}_{-\tau} |\dot{p}_1|^\gamma \, dt &\sim \frac{\lambda^{\gamma}}{2^{\gamma}} \int^{\tau}_{\tau_0} |B e^{\lambda t}|^{\gamma} \, dt \sim \dfrac{|B|^{\gamma} \lambda^{\gamma-1}}{\gamma 2^{\gamma}} \, e^{\gamma \lambda \tau}
\end{split}
\]
In summary, the asymptotic behavior as $\tau \gg 1$ of the hyperbolic component of LDs is
\begin{equation}
M_{\gamma}^{h}\left(\mathbf{x}_0,\tau\right) \sim \dfrac{\lambda^{\gamma - 1}}{\gamma 2^{\gamma - 1}} \left(|A|^{\gamma} + |B|^{\gamma}\right) \, e^{\gamma \lambda \tau}
\label{M_hyp_asymp}
\end{equation}
Therefore, we have show that $M_{\gamma}^{h}$ grows exponentially with $\tau$ and also that the leading order singularities in $M_{\gamma}^{h}$ occur when $|A| = 0$, that is, when $p^0_1 = -  q^0_1$, which corresponds to initial conditions on the stable manifold of the NHIM, or in the case where $|B| = 0$, that is, $p^0_1 = q^0_1$, corresponding to initial conditions on the unstable manifold of the NHIM. Moreover, $M_{\gamma}^{h}$ is non-differentiable at the NHIM, since it is given by the intersection of the stable and unstable manifolds. 

To finish the proof we focus now on the elliptic component of LDs. We have that
\[
M^{e}_{\gamma}(\mathbf{x}_{0},\tau) = \sum_{i=2}^{3} M^{e,i}_{\gamma}(\mathbf{x}_{0},\tau) = \sum_{i=2}^{3} \int^{\tau}_{-\tau} |\dot{q}_{i}(t;\mathbf{x}_0)|^\gamma + |\dot{p}_{i}(t;\mathbf{x}_0)|^\gamma \; dt
\]   
Since all the bath modes follow the same type of dynamics, we focus on analyzing the contribution of the bath mode $(q_i,p_i)$ where $i \in \lbrace 2,3 \rbrace$ to the elliptic component of LDs, that is, $M^{e,i}_{\gamma}(\mathbf{x}_{0},\tau)$. We recall that
\[
\dot{q}_i(t) =  p^0_i \omega_i \cos(\omega_i t) - q^0_i \omega_i \sin(\omega_i t) \quad,\quad \dot{p}_i(t) = -\omega_i q^0_i \cos(\omega_i t) - \omega_i p^0_i \sin(\omega_i t)
\]
Then we have to evaluate the integral
\[
M^{e,i}_{\gamma}(\mathbf{x}_{0},\tau) = \int_{-\tau}^{\tau} |\dot{q}_{i}(t;\mathbf{x}_0)|^\gamma + |\dot{p}_{i}(t;\mathbf{x}_0)|^\gamma \; dt
\]
We deal first with the integral
\[
\int_{-\tau}^{\tau} |\dot{q}_{i}(t;\mathbf{x}_0)|^\gamma \; dt = \int_{-\tau}^{\tau} |p^0_i \omega_i \cos(\omega_i t) - q^0_i \omega_i \sin(\omega_i t)|^{\gamma} \; dt
\]
Using that $|\dot{q}_i|$ is periodic with period $T_i = \pi / \omega_i$, we can always write $\tau = N T_i + r$ for some integer $N$ and $r \in [0,T_i]$. With this decomposition we can write
\[
\begin{split}
\int_{-\tau}^{\tau} |\dot{q}_{i}|^\gamma \; dt &= \int_{-N T_i - r}^{-N T_i} |\dot{q}_{i}|^\gamma \; dt + \int_{-N T_i}^{N T_i} |\dot{q}_{i}|^\gamma \; dt + \int_{N T_i}^{N T_i + r} |\dot{q}_{i}|^\gamma \; dt = \\[.2cm] 
&= 2N \int_{0}^{T_i} |\dot{q}_{i}|^\gamma \; dt + \int_{-r}^{r} |\dot{q}_{i}|^\gamma \; dt
\end{split}
\]
Since the first integral covers one whole period, without loss of generality we can compute its value by starting on the initial condition $p_i^0 = 0$ and $q_i^0 = \sqrt{2\mathcal{H}_i^0 /  \omega_i}$. We can do this because the bath mode $(q_i,p_i)$ follows a circle in the $q_i$-$p_i$ plane with that radius. This yields
\[
\begin{split}
\int_{0}^{T_i} |\dot{q}_{i}|^\gamma \; dt &= (2\mathcal{H}_i^0 \omega_i)^{\gamma / 2} \int_{0}^{T_i} |\sin(\omega_i t)|^\gamma \; dt = \dfrac{(2\mathcal{H}_i^0 \omega_i)^{\gamma / 2}}{\omega_i} \int_{0}^{\pi} |\sin u|^\gamma \; du = \\[.2cm]
&= \dfrac{2(2\mathcal{H}_i^0)^{\gamma / 2}}{\omega_i^{1 - \gamma / 2}} \int_{0}^{\frac{\pi}{2}} \sin^{\gamma}u \; du = \dfrac{(2\mathcal{H}_i^0)^{\gamma / 2}}{\omega_i^{1 - \gamma / 2}} B\left(\frac{\gamma+1}{2},\frac{1}{2}\right) = \\[.2cm]
& = \dfrac{(2\mathcal{H}_i^0)^{\gamma / 2}}{\omega_i^{1 - \gamma / 2}} \dfrac{\Gamma\left(\dfrac{\gamma+1}{2}\right)\Gamma\left(\dfrac{1}{2}\right)}{\Gamma\left(\dfrac{\gamma}{2} + 1\right)} = \dfrac{\sqrt{\pi} \, (2\mathcal{H}_i^0)^{\gamma / 2}}{\omega_i^{1 - \gamma / 2}} \dfrac{\Gamma\left(\dfrac{\gamma+1}{2}\right)}{\dfrac{\gamma}{2}\Gamma\left(\dfrac{\gamma}{2}\right)} = \\[.2cm]
&= \dfrac{2\sqrt{\pi} \, (2\mathcal{H}_i^0)^{\gamma / 2}}{\gamma \, \omega_i^{1 - \gamma / 2}} \dfrac{\Gamma\left(\dfrac{\gamma+1}{2}\right)}{\Gamma\left(\dfrac{\gamma}{2}\right)} 
\end{split}
\]
where $B$ is the beta function and $\Gamma$ is the gamma function. Threfore we have shown that
\begin{equation}
\int_{-\tau}^{\tau} |\dot{q}_{i}|^\gamma \; dt = \dfrac{4 N \sqrt{\pi} \, (2\mathcal{H}_i^0)^{\gamma / 2}}{\gamma \, \omega_i^{1 - \gamma / 2}} \dfrac{\Gamma\left(\dfrac{\gamma+1}{2}\right)}{\Gamma\left(\dfrac{\gamma}{2}\right)} +  \int_{-r}^{r} |\dot{q}_{i}|^\gamma \; dt
\label{ell_qi_int}
\end{equation}
Followng this argument, we can also compute the integral
\[
\int_{-\tau}^{\tau} |\dot{p}_{i}(t;\mathbf{x}_0)|^\gamma \; dt = \int_{-\tau}^{\tau} |\omega_i q^0_i \cos(\omega_i t) + \omega_i p^0_i \sin(\omega_i t)|^{\gamma} \; dt
\]
by using the decomposition
\[
\begin{split}
\int_{-\tau}^{\tau} |\dot{p}_{i}|^\gamma \; dt &= \int_{-N T_i - r}^{-N T_i} |\dot{p}_{i}|^\gamma \; dt + \int_{-N T_i}^{N T_i} |\dot{p}_{i}|^\gamma \; dt + \int_{N T_i}^{N T_i + r} |\dot{p}_{i}|^\gamma \; dt = \\[.2cm] 
&= 2N \int_{0}^{T_i} |\dot{p}_{i}|^\gamma \; dt + \int_{-r}^{r} |\dot{p}_{i}|^\gamma \; dt
\end{split}
\]
The integral over one period in this case, choosing $q_i^0 = 0$ and $p_i^0 = \sqrt{2\mathcal{H}_i^0 /  \omega_i}$, is
\[
\begin{split}
\int_{0}^{T_i} |\dot{p}_{i}|^\gamma \; dt &= (2\mathcal{H}_i^0 \omega_i)^{\gamma / 2} \int_{0}^{T_i} |\sin(\omega_i t)|^\gamma \; dt = \dfrac{2\sqrt{\pi} \, (2\mathcal{H}_i^0)^{\gamma / 2}}{\gamma \, \omega_i^{1 - \gamma / 2}} \dfrac{\Gamma\left(\dfrac{\gamma+1}{2}\right)}{\Gamma\left(\dfrac{\gamma}{2}\right)} 
\end{split}
\]
This gives:
\begin{equation}
\int_{-\tau}^{\tau} |\dot{p}_{i}|^\gamma \; dt = \dfrac{4 N \sqrt{\pi} \, (2\mathcal{H}_i^0)^{\gamma / 2}}{\gamma \, \omega_i^{1 - \gamma / 2}} \dfrac{\Gamma\left(\dfrac{\gamma+1}{2}\right)}{\Gamma\left(\dfrac{\gamma}{2}\right)} + \int_{-r}^{r} |\dot{p}_{i}|^\gamma \; dt
\label{ell_pi_int}
\end{equation}
Consequently, the elliptic component of LDs corresponding to the bath mode $(q_i,p_i)$ is
\[
\begin{split}
M^{e,i}_{\gamma}(\mathbf{x}_{0},\tau) &= \dfrac{8 N \sqrt{\pi} \, (2\mathcal{H}_i^0)^{\gamma / 2}}{\gamma \, \omega_i^{1 - \gamma / 2}} \dfrac{\Gamma\left(\dfrac{\gamma+1}{2}\right)}{\Gamma\left(\dfrac{\gamma}{2}\right)} +  M^{e,i}_{\gamma}(\mathbf{x}_{0},r) = \\[.2cm]
& = \dfrac{4 (2\omega_i\mathcal{H}_i^0)^{\gamma / 2}}{\sqrt{\pi}} \dfrac{\Gamma\left(\dfrac{\gamma+1}{2}\right)}{\Gamma\left(\dfrac{\gamma}{2} + 1\right)} \left(\tau - r\right) +  M^{e,i}_{\gamma}(\mathbf{x}_{0},r)
\end{split}
\]
which clearly shows that it is a linearly increasing function of $\tau$, in contrast to the hyperbolic component of LDs that grows exponentially with the integration time. If we divide the elliptic part of LDs by $2\tau$ (the integration period of the trajectory) and we let $\tau \to \infty$ we get that
\begin{equation}
\begin{split}
\lim_{\tau \to \infty} \dfrac{1}{2\tau} M^{e,i}_{\gamma}(\mathbf{x}_{0},\tau) &= \dfrac{2 (2\omega_i\mathcal{H}_i^0)^{\gamma / 2}}{\sqrt{\pi}} \dfrac{\Gamma\left(\dfrac{\gamma+1}{2}\right)}{\Gamma\left(\dfrac{\gamma}{2} + 1\right)} =  (2\omega_i\mathcal{H}_i^0)^{\gamma / 2} \dfrac{2}{\pi} B\left(\dfrac{\gamma+1}{2},\dfrac{1}{2}\right) = \\ &= \dfrac{2}{\pi} \left(\omega_i R_i\right)^{\gamma} B\left(\dfrac{\gamma+1}{2},\dfrac{1}{2}\right)
\end{split}
\label{lim_m_gamma}
\end{equation}
where $R_i$ is the radius of the circle described by the bath coordinate $(q_i,p_i)$. So the time average of the elliptic part of LDs converges to a value that depends on the energy $\mathcal{H}_i^0$ of the initial condition $(q_i^0,p_i^0)$, and therefore in the limit its value is constant on the periodic orbit described by the bath mode on the $q_i$-$p_i$ plane. This result is of interest because it allows us to use LDs as a tool to recover phase space KAM tori, given that the dynamical system under study satisfies the conditions of the Ergodic Partition Theorem \cite{mezic3}. Therefore, if one analyzes thr time average of LDs on a specific initial condition in phase space $\mathbf{x}_0$, and its value converges, then this initial copndition would lie in an invariant phase space set cosisting of points that share the same time average value. Therefore, the contours of the time average of LDs, when it converges, identify invariant phase space sets. It is also important to highlight that the limit value to which $M^{e,i}$ converges is directly related to the limit value to which the classical definition of LDs, representing the arclength of the trajectory of the bath mode, converges. Indeed:
\begin{equation}
M^{e,i}(\mathbf{x}_{0},\tau) = \int_{-\tau}^{\tau} \sqrt{\left(\dot{q}_{i}(t;\mathbf{x}_0)\right)^2 + \left(\dot{p}_{i}(t;\mathbf{x}_0)\right)^2 } \; dt = 2\tau \omega_i R_i
\label{lim_m}
\end{equation}
From Eqs. \eqref{lim_m_gamma} and \eqref{lim_m} it follows that
\[
\displaystyle{\lim_{\tau \to \infty} \dfrac{1}{2\tau} M^{e,i}_{\gamma}} = \dfrac{2}{\pi} \left(\lim_{\tau \to \infty} \dfrac{1}{2\tau} M^{e,i}\right)^{\gamma} B\left(\dfrac{\gamma+1}{2},\dfrac{1}{2}\right)
\]

To conclude the proof we will show how LDs attains a local minimum at the phase space points corresponding to the stable and unstable manifolds of the NHIM and a global minimum at the NHIM. Given an energy of the system $\mathcal{H}_0$ above that of the origin, we know that a family of NHIM parametrized by the energy bifurcates from the rank-1 saddle. The phase space points that lie on the stable (or unstable) manifold, that is $q_1^0 = -p_1^0$ ($q_1^0 = p_1^0$), contribute to the hyperbolic component of LDs described in Eq. \eqref{M_hyp_asymp} with $|A| = 0$ ($|B| = 0$) so that the Lagrangian descriptor has a local minimum at the manifold. Moreover, if the initial condition is on the NHIM, that is $q_1^0 = p_1^0 = 0$, then $|A| = |B| = 0$ and consequently $M^{h}_{\gamma}(\mathbf{x}_{0},\tau) = 0$. Furthermore, all the energy of the system for these points concentrates on the bath modes that evolve periodically, which implies that $M^{e}_{\gamma}(\mathbf{x}_{0},\tau) > 0$. As a result, Lagrangian descriptors attain a global minimum at these points. We would like to conclude by noting that due to Moser's theorem (generalization of the Liapunov's subcenter theorem) on index-1 saddle in a Hamiltonian system, this result on Lagrangian descriptor for a linear system extends to the local neighborhood of an index-1 saddle equilibrium point in the full nonlinear system~\cite{Moser1958b}.

\section{Two Degrees of Freedom}
\label{sec:normal_form_2dof}

In this section we use the Lagrangian descriptor method to identify the NHIM and its stable and unstable manifolds for a 2 DoF separable quadratic Hamiltonian with an index-1 saddle at the origin. The advantage of using this model Hamiltonian is that we can compare the analytical expressions for the phase space structures with features in LD plots.

\subsection{Decoupled quadratic Hamiltonian}
We consider a quadratic Hamiltonian for a two degree-of-freedom system given by
\begin{equation}
H(q_1, p_1, q_2, p_2) = \underbrace{\frac{\lambda}{2} (p_1^2 - q_1^2)}_\text{$H_r$} + 
\underbrace{\frac{\omega_2}{2}(q_2^2 + p^2_2)}_\text{$H_b$},
\quad 
\lambda,\omega_2 > 0 
\label{eqn:sqh_2dof}
\end{equation}
with the corresponding vector field
\begin{equation}
\begin{aligned}
\dot{q}_1 = & \phantom{-} \frac{\partial H_2}{\partial p_1} = \phantom{-}\lambda p_1, \\
\dot{p}_1 = & -\frac{\partial H_2}{\partial q_1} = \phantom{-}\lambda q_1,\\
\dot{q}_2 = & \phantom{-}\frac{\partial H_2}{\partial p_2} = \phantom{-}\omega_2 p_2,\\
\dot{p}_2 = & -\frac{\partial H_2}{\partial q_2} = -\omega_2 q_2\\
\end{aligned}
\label{eqn:eom_nf_2dof}
\end{equation}
The equilibirum point is located at $(0,0,0,0)$ on the zero total energy surface. It is trivial to check that the eigenvalues of the linearized system around the equilibrium point are $\pm \lambda$ and $\pm i \omega_2$, and hence the equilibrium point is of saddle $\times$ center type. In this form, the Hamiltonian~\eqref{eqn:sqh_2dof} is decoupled into the ``reactive'' mode given by $H_r$ and the ``bath'' mode given by $H_b$, hence it will be referred to as separable quadratic Hamiltonian (SQH). This representation lends to discussing the distribution of total energy between the two modes for the phase space structures in uncoupled coordinates. 
In this form, a chemical reaction is said to have occurred when the $q_1$ coordinate of a trajectory changes sign and thus, an isoenergetic, $H(q_1, p_1, q_2, p_2) = h$, dividing surface (DS) can be defined by $q_1 = 0$ hypersurface. The constant energy defines a three dimensional surface in the four dimensional phase space and is given by
\begin{equation}
\frac{\lambda}{2} (p_1^2 - q_1^2) + \frac{\omega_2}{2}(p_2^2 + q^2_2) = H_r + H_b = h > 0, \qquad H_r > 0, \quad H_b \geq 0
\end{equation}
The dividing surface, $q_1 = 0$, for a constant energy is 
\begin{equation}
\frac{\lambda}{2} p_1^2 + \frac{\omega_2}{2}(p_2^2 + q^2_2) = H_r + H_b = h > 0, \qquad H_r > 0, \quad H_b \geq 0 \label{eqn:ds_sqh_2dof}
\end{equation}
which is a two dimensional surface, or precisely of geometry $\mathbb{S}^2$, that is a 2-sphere on the three dimensional energy surface. Thus it is codimension-1 and partitions the energy surface into reactant $p_1 - q_1 > 0$ and product $p_1 - q_1 < 0$ regions by the forward and backward ``reaction'' dividing surfaces as shown in Ref.~\cite{waalkens2004direct} and given by
\begin{equation}
\begin{aligned}
p_1 & = & \pm \sqrt{\frac{2}{\lambda}\left( H_r + H_b - \frac{\omega_2}{2}(p_2^2 + q_2^2) \right)}, \quad \text{forward/backward DS} 
\end{aligned}
\end{equation}
The forward and backward DS is joined at $p_1 = 0$, thus
\begin{equation}
\mathcal{M}(h) = \left\{ (q_1, p_1, q_2, p_2) \; \vert \; q_1 = 0, p_1 = 0, \frac{\omega_2}{2}\left( p_2^2 + q_2^2 \right) = H_b = h \geq 0 \right\} \qquad \text{NHIM}
\label{eqn:sep_quad_ham2dof_nhim}
\end{equation}
which is of geometry $\mathbb{S}^1$, that is a circle centered at the origin and of radius $\sqrt{2h/\omega_2}$ in the $(q_2,p_2)$ plane. This is a {\em normally hyperbolic invariant manifold} (NHIM) associated with the index-1 saddle and parametrized by total energy $H(q_1, p_1, q_2, p_2) = h$~\cite{wiggins2013normally}. {\em Invariance} follows from the vector field~\eqref{eqn:eom_nf_2dof} since when $q_1 = p_1 =0$, $\dot{q}_1 = \dot{p}_1 =0$. Thus $q_1$ and $p_1$ always remain zero, and trajectories with these initial conditions remain on the NHIM, that is, $q_1 = p_1 =0$ is invariant. It is {\em normally hyperbolic} since the directions normal to the NHIM, that is, the $(q_1, p_1)$ surface, have linear saddle like dynamics. For a two degree-of-freedom system, this NHIM is more commonly referred in the literature as an unstable periodic orbit.

In order to understand the relationship between the NHIM and the index 1 saddle point, we note that for $H_r =0$ and $H_b =0$ the NHIM reduces to the point $(q_1, p_1, q_2, p_2) = (0, 0, 0, 0)$, which is the index-1 saddle point on the energy surface $H_r + H_b = 0$. Therefore, as the total energy is increased from $0$, with $H_b$ increasing from zero, the NHIM ``grows'' from the index-1 saddle point on the zero energy surface into an invariant 1-sphere. This shows how the ``influence'' of the index-1 saddle point is carried to higher energy sufaces on which the saddle point does not exist.

The stable and unstable manifolds of the NHIM~\eqref{eqn:sep_quad_ham2dof_nhim} are given by 
\begin{align}
\mathcal{W}^{\rm u}(\mathcal{M}(h)) =& \big\{  (q_1, p_1, q_2, p_2) \; | \; q_1 = p_1, \frac{\omega_2}{2}\left( p_2^2 + q_2^2 
\right) = H_b > 0  \big\}, \label{eqn:quad_ham2dof_umani}\\
\mathcal{W}^{\rm s}(\mathcal{M}(h)) =& \big\{  (q_1, p_1, q_2, p_2) \; | \; q_1 = -p_1, \frac{\omega_2}{2}\left( p_2^2 + q_2^2 
\right) = H_b > 0  \big\} \label{eqn:quad_ham2dof_smani}
\end{align}
which is two dimensional surface and have geometry $\mathbb{R} \times \mathbb{S}^1$ for a fixed energy. Thus, the codimension-1 geometry of the manifolds partition the phase space into ``reactive'' and ``non-reactive'' trajectories as shown in Fig.~\ref{fig:ps_struct_quad_ham2dof}. The detailed visualization is only possible for the two degrees-of-freedom system which is also the starting point for testing the Lagrangian descriptor based approach for detecting the phase space structures.

\begin{figure}[!ht]
	\centering
	\includegraphics[width=0.98\textwidth]{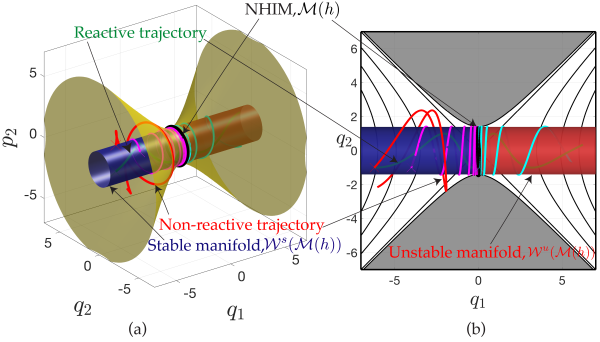}
	\caption{(a) Shows the phase space structures, namely the NHIM which is an unstable periodic orbit (as black circle), its stable and unstable manifolds (blue and red surfaces), and the energy surface (as yellow surface) for the two degrees-of-freedom decoupled quadratic Hamiltonian~\eqref{eqn:sqh_2dof}. Four example trajectories illustrate the dynamics mediated by these phase space structures: the red and green trajectories have the same initial configuration space coordinates, but are located outside and inside the codimension-1 invariant manifolds which separate the ``reactive'' (green) and ``non-reactive'' (red) trajectories. Also shown are the two trajectories on the manifolds (magenta on the stable manifold and cyan on the unstable manifold) that show the dynamics on the invariant surfaces. (b) Shows the projection of the stable (blue) and unstable (red) manifolds in the configuration space coordinates $(q_1, q_2)$ along with the projection of the four example trajectories. The equipotential contours are the black lines and the region of inaccessible motion is shadded as grey.}
	\label{fig:ps_struct_quad_ham2dof}
\end{figure}

\subsubsection{Detecting the unstable periodic orbit and its manifolds}

To identify NHIM and its invariant manifolds, we compute the Lagrangian descriptor in a square domain of size 2 units around the origin and discretize the coordinates of the two-dimensional surface. Next, we pick a constant value for one of the two remaining coordinates, and use the total energy equation to solve for the fourth coordinate. Due to the form of the Hamiltonian~\eqref{eqn:sqh_2dof}, obtaining the coordinate from the constant energy condition reduces simply to solving a quadratic equation. 


\newpoint~Isoenergetic two-dimensional surface parametrized by $(q_1, p_1)$~\textemdash~On the constant energy surface, $H(q_1, p_1, q_2, p_2) = h$, we compute Lagrangian descriptor on a two-dimensional surface parametrized by $(q_1, p_1)$ coordinates by defining 
\begin{equation}
U_{q_1p_1}^+ = \left\{ (q_1, p_1, q_2, p_2) \; | \; q_2 = 0, \dot{q}_2 > 0 : p_2(q_1, p_1, q_2; h) > 0 \right\}
\end{equation}
where
\begin{align}
p_2(q_1, p_1, q_2 = 0; h) = \sqrt{\frac{2}{\omega_2}\left( h - \frac{\lambda}{2}(p_1^2 - q_1^2) \right)} 
\end{align}
The intersection of the two-dimensional surface $U_{q_1p_1}^+$ with the NHIM~\eqref{eqn:sep_quad_ham2dof_nhim} becomes

\begin{align}
\mathcal{M}(h) \cap U_{q_1p_1}^+ = \left\{ (q_1, p_1, q_2, p_2) \; \vert \; \right. & \left. \kern - \nulldelimiterspace p_1 = 0, q_1 = 0, q_2 = 0, \dot{q}_2 > 0 : \right. \nonumber \\ 
& \left.\kern - \nulldelimiterspace  p_2(q_1, p_1, q_2; h) > 0 \right\}.
\end{align}
%
Thus, the NHIM is located at the origin $(0,0)$ and marked by a red cross in the LD plot (Fig.~\ref{fig:LD_SQH_2dof_q1p1}).  Furthermore, the intersection of the two-dimensional surface with the unstable~\eqref{eqn:quad_ham2dof_umani} and stable manifolds~\eqref{eqn:quad_ham2dof_smani} is given by

\begin{align}
\mathcal{W}^u(\mathcal{M}(h)) \cap U_{q_1p_1}^+ = \left\{ (q_1, p_1, q_2, p_2) \; \vert \; \right. & \left. \kern - \nulldelimiterspace p_1 = q_1, q_2 = 0, \dot{q}_2 > 0 : \right. \nonumber \\  
& \left.\kern - \nulldelimiterspace p_2(q_1, p_1, q_2;h) > 0 \right\}, \\
\mathcal{W}^s(\mathcal{M}(h)) \cap U_{q_1p_1}^+ = \left\{ (q_1, p_1, q_2, p_2) \; \vert \; \right. & \left. \kern - \nulldelimiterspace p_1 = -q_1, q_2 = 0, \dot{q}_2 > 0 : \right. \nonumber \\  
& \left.\kern - \nulldelimiterspace p_2(q_1, p_1, q_2;h) > 0 \right\},
\end{align}

which are one-dimensional for a fixed energy, and represent lines passing through the origin shown as dashed red (unstable) and white (stable) lines, respectively, in Fig.~\ref{fig:LD_SQH_2dof_q1p1}. The only points of local minima in the LD plot (Fig.~\ref{fig:LD_SQH_2dof_q1p1}) also lie along the lines passing through the origin and correspond to the manifolds of the NHIM. 

\begin{figure}[!th]
	\centering
	\subfigure[]{\includegraphics[width=0.32\textwidth]{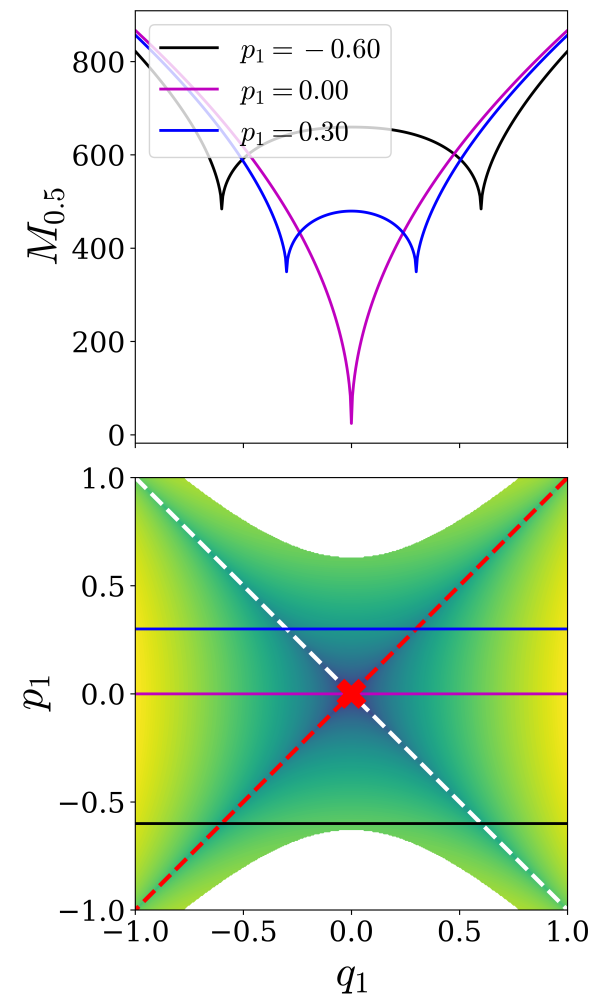}\label{fig:LD_SQH_2dof_q1p1}}
	\subfigure[]{\includegraphics[width=0.32\textwidth]{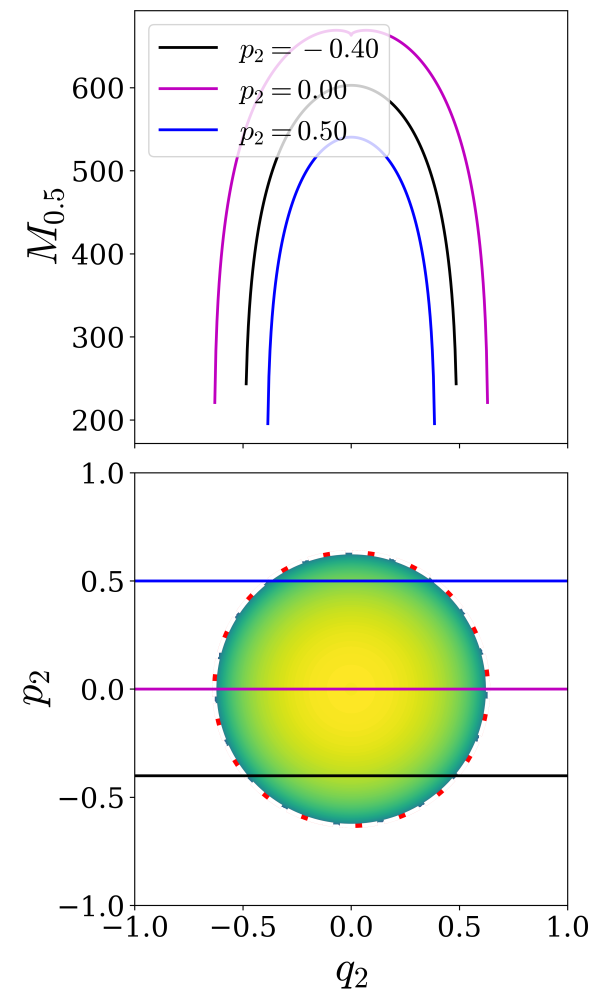}\label{fig:LD_SQH_2dof_q2p2}}
	\subfigure[]{\includegraphics[width=0.32\textwidth]{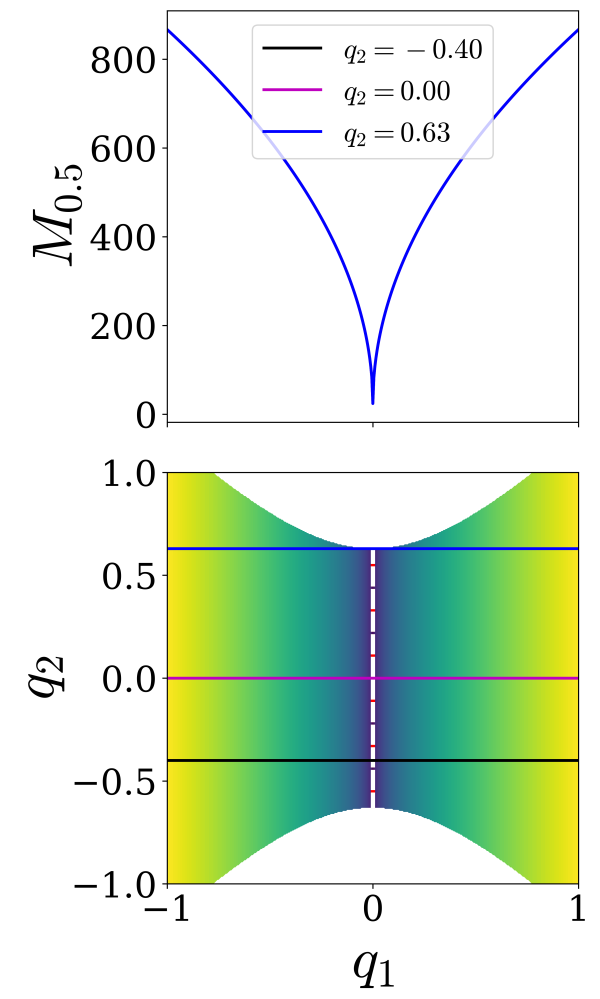}\label{fig:LD_SQH_2dof_q1q2}}
	\caption{Lagrangian descriptor slice for the two degrees of freedom separable quadratic Hamiltonian on the isoenergetic two-dimensional surfaces \protect\subref{fig:LD_SQH_2dof_q1p1} $U_{q_1p_1}^+$, \protect\subref{fig:LD_SQH_2dof_q2p2} $U_{q_2p_2}^+$, \protect\subref{fig:LD_SQH_2dof_q1q2} $U_{q_1q_2}^+$. Parameters used are $\lambda = 1.0, \omega_2 = 1.0$ for constant total energy $H_2 = 0.2$ and integration time $\tau = 10$ is fixed.}
	\label{fig:LD_SQH_2dof_1}
\end{figure}

\newpoint~Isoenergetic two-dimensional surface parametrized by $(q_2, p_2)$~\textemdash~Next, on the fixed energy surface, $H(q_1, p_1, q_2, p_2) = h$, we compute the Lagrangian descriptor on a two-dimensional surface parametrized by $(q_2, p_2)$ coordinates by defining 

\begin{equation}
U_{q_2p_2}^+ = \left\{ (q_1, p_1, q_2, p_2) \; | \; q_1 = 0, \dot{q}_1 \geqslant 0 : p_1(q_1, q_2, p_2; h) \geqslant 0 \right\}
\end{equation}
where
\begin{align}
p_1(q_1 = 0, q_2, p_2; h) =& \sqrt{\frac{2}{\lambda}\left( h  - \frac{\omega_2}{2}\left( p_2^2 + q_2^2 \right) \right)}
\end{align}

The intersection of the two-dimensional surface $U_{q_2p_2}^+$ with the NHIM~\eqref{eqn:sep_quad_ham2dof_nhim} is given by

\begin{align}
\mathcal{M}(h) \cap U_{q_2p_2}^+ = \left\{ (q_1, p_1, q_2, p_2) \; \vert \; \right. & \left. \kern - \nulldelimiterspace p_1 = 0, q_1 = 0, \dot{q}_1 > 0, \right. \nonumber \\ 
& \left.\kern - \nulldelimiterspace   \frac{\omega_2}{2}\left( p_2^2 + q_2^2 \right) = h \right\}.
\end{align}

which represents a circle of radius $\sqrt{2h/\omega_2}$. The radius is $\approx 0.632$ for $h = 0.2, \omega_2 = 1.0$ and is marked as a dashdot red circle in Fig.~\ref{fig:LD_SQH_2dof_q2p2}. The points on the circle are also locations of minima as shown by the one-dimensional slices at constant $p_2$.

Next, intersection of the two-dimensional $U_{q_2p_2}^+$ with the stable~\eqref{eqn:quad_ham2dof_smani} and unstable manifolds~\eqref{eqn:quad_ham2dof_umani} is given by

\begin{align}
\mathcal{W}^u(\mathcal{M}(h)) \cap U_{q_2p_2}^+ = \left\{ (q_1, p_1, q_2, p_2) \; \vert \; \right. & \left. \kern - \nulldelimiterspace p_1 = q_1, q_1 = 0, \dot{q}_1 > 0 : \right. \nonumber \\  
& \left.\kern - \nulldelimiterspace \frac{\omega_2}{2}\left( p_2^2 + q_2^2 \right) = h \right\}, \\
\mathcal{W}^s(\mathcal{M}(h)) \cap U_{q_2p_2}^+ = \left\{ (q_1, p_1, q_2, p_2) \; \vert \; \right. & \left. \kern - \nulldelimiterspace p_1 = -q_1, q_1 = 0, \dot{q}_1 > 0 : \right. \nonumber \\  
& \left.\kern - \nulldelimiterspace \frac{\omega_2}{2}\left( p_2^2 + q_2^2 \right) = h \right\},
\end{align}

which are circles of same radius $\sqrt{2h/\omega_2}$. These manifolds are identified along the boundary of the LD plot in Fig.~\ref{fig:LD_SQH_2dof_q2p2} along which the minima also occurs.

\newpoint~Isoenergetic two-dimensinal surface parametrized by $(q_1, q_2)$~\textemdash~Next, on the constant energy surface, $H(q_1, p_1, q_2, p_2) = h$, we compute the Lagrangian descriptor on a two-dimensional surface by defining

\begin{equation}
U_{q_1q_2}^+ = \left\{ (q_1, p_1, q_2, p_2) \; | \; p_1 = 0, p_2(q_1, p_1, q_2; h) > 0 \right\}
\end{equation}

where
\begin{align}
p_2(q_1, p_1 = 0, q_2; h) =& \sqrt{\frac{2}{\omega_2}\left( h  + \frac{\lambda}{2}q_1^2 - \frac{\omega_2}{2}q_2^2 \right)}
\end{align}

The intersection of the NHIM~\eqref{eqn:sep_quad_ham2dof_nhim} with the two-dimensional surface $U_{q_1q_2}^+$ is

\begin{align}
\mathcal{M}(h) \cap U_{q_1q_2}^+ = \left\{ (q_1, p_1, q_2, p_2) \; \vert \; \right. & \left. \kern - \nulldelimiterspace p_1 = 0, q_1 = 0, p_2 > 0 : \right. \nonumber \\ 
& \left.\kern - \nulldelimiterspace   \frac{\omega_2}{2}\left( p_2^2 + q_2^2 \right) = h \right\}.
\end{align}

which represent points on the line $q_1 = 0$ which is marked by the dashdot line in Fig.~\ref{fig:LD_SQH_2dof_q1q2}.

Next, intersection of the two-dimensional $U_{q_1q_2}^+$ with the stable~\eqref{eqn:quad_ham2dof_smani} and unstable manifolds~\eqref{eqn:quad_ham2dof_umani} is given by

\begin{align}
\mathcal{W}^u(\mathcal{M}(h)) \cap U_{q_1q_2}^+ = \left\{ (q_1, p_1, q_2, p_2) \; \vert \; \right. & \left. \kern - \nulldelimiterspace p_1 = q_1, p_1 = 0, p_2 > 0 : \right. \nonumber \\  
& \left.\kern - \nulldelimiterspace \frac{\omega_2}{2}\left( p_2^2 + q_2^2 \right) = h \right\}, \\
\mathcal{W}^s(\mathcal{M}(h)) \cap U_{q_1q_2}^+ = \left\{ (q_1, p_1, q_2, p_2) \; \vert \; \right. & \left. \kern - \nulldelimiterspace p_1 = -q_1, p_1 = 0, p_2 > 0 : \right. \nonumber \\  
& \left.\kern - \nulldelimiterspace \frac{\omega_2}{2}\left( p_2^2 + q_2^2 \right) = h \right\},
\end{align}

which denote lines parallel to the $q_2$ axis and marked by the dashed red (unstable) and white (stable) lines in Fig.~\ref{fig:LD_SQH_2dof_q1q2}. These manifolds are identified by minima in the Lagrangian descriptor plot in Fig.~\ref{fig:LD_SQH_2dof_q1q2} as also shown in the one-dimensional slice for constant $q_2$.

Now the question is how do these phase space structures manifest when the system is non-separable into ``reactant'' and ``product'' coordinates. We will answer this using a linear symplectic transformation of the separable quadratic Hamiltonian and then compute the LD on different 2D surfaces for a fixed energy.

%

\subsection{Coupled quadratic Hamiltonian}

To couple the coordinates in the separable quadratic Hamiltonian~\eqref{eqn:sqh_2dof}, we introduce a linear transformation, $C$, such that it satisfies the symplectic condition
\begin{equation}
C \mathcal{J} C^T = \mathcal{J} = \begin{pmatrix}
0_N & I_N \\
-I_N & 0_N
\end{pmatrix}
\label{eqn:symp_cond}
\end{equation}
where $\mathcal{J}$ is the 2N $\times$ 2N matrix and $I_N$ denotes the $N \times N$ identity matrix. The symplectic transformation $C$ acts on the non-separable (coupled) coordinates $(x, y, p_x, p_y)$ to give the decoupled coordinates $(q_1, q_2, p_1, p_2)$: 
\begin{equation}
\begin{bmatrix}
q_1 \\
q_2 \\
p_1 \\
p_2
\end{bmatrix}
= C
\begin{bmatrix}
x \\
y \\
p_x \\
p_y 
\end{bmatrix} 
\implies
\begin{bmatrix}
\dot{q_1} \\
\dot{q_2} \\
\dot{p_1} \\
\dot{p_2}
\end{bmatrix}
= C
\begin{bmatrix}
\dot{x} \\
\dot{y} \\
\dot{p_x} \\
\dot{p_y} 
\end{bmatrix},
\quad
\begin{bmatrix}
\dot{x} \\
\dot{y} \\
\dot{p_x} \\
\dot{p_y} 
\end{bmatrix} 
= -C 
\begin{bmatrix}
\dot{q_1} \\
\dot{q_2} \\
\dot{p_1} \\
\dot{p_2}
\end{bmatrix}
\end{equation}
since, $C$ satisfies $C^T = C^{-1} = -C$ and some examples are given in Appendix~\ref{sect:examples_C}. 

Let us consider the symplectic transformation
\begin{equation}
C =
\begin{pmatrix}
0 & 0 & 1 & 0 \\
0 & 0 & 0 & 1 \\
-1 & 0 & 1 & 1 \\
0 & -1 & 1 & 1
\end{pmatrix}
\label{eqn:two_dof_C}	
\end{equation}
The coupled (non-separable) coordinates become
\begin{equation}
\begin{aligned}
q_1 &= p_x \\
q_2 &= p_y \\
p_1 &= -x + p_x + p_y \\
p_2 &= -y + p_x + p_y
\end{aligned}
\end{equation}

This set of transformed coordinates applied to~\eqref{eqn:sqh_2dof} gives the non-separable quadratic Hamiltonian 
\begin{align}
\mathcal{H}(x, p_x, y, p_y) = \frac{\lambda}{2} x^2 + \frac{\omega_2}{2} y^2 + & \frac{\omega_2}{2} \left( 2p_y^2 + p_x^2 + 2 p_x p_y - 2 y p_x - 2 y p_y \right) + \nonumber \\ 
& \frac{\lambda}{2} \left( p_y^2 + 2p_xp_y - 2xp_x - 2xp_y \right) 
\label{eqn:ham_symp_2dof}
\end{align}
which gives the Hamiltonian vector field
\begin{equation}
\begin{aligned}
\dot{x} & = -\lambda x - \omega_2 y + \omega_2 p_x + (\lambda + \omega_2) p_y \\
\dot{p}_x & = \phantom{-}\lambda( -x + p_x + p_y )  \\
\dot{y} & = -\lambda x - \omega_2 y + (\lambda + \omega_2) p_x + (\lambda + 2\omega_2) p_y \\
\dot{p}_y & =  \phantom{-}\omega_2( -y + p_x + p_y )
\end{aligned}
\label{eqn:eom_symp_2dof}
\end{equation}
where the equilibrium point is at $(0,0,0,0)$, and its total energy is $0$. The Jacobian at this equilibrium point has eigenvalues of type $\lambda, -\lambda, i\omega_2, -i \omega_2$ (as shown in the Appendix~\ref{sect:linear_symp_transf_jac}) and is of saddle $\times$ center type, that is index-1. 


In the decoupled (separable) quadratic Hamiltonian, the dividing surface is given by $q_1 = 0$ Eqn.~\eqref{eqn:ds_sqh_2dof} which in the transformed coordinates becomes $p_x = 0$. The dividing surface on the fixed energy surface $\mathcal{H} = h$ is given by
\begin{equation}
{\rm DS} = \left\{ (x, p_x, y, p_y) \, | \, \frac{\lambda}{2} p_y^2 + \omega_2 p_y^2 - \lambda x p_y - \omega_2 y p_y + \frac{\lambda}{2} x^2 + \frac{\omega_2}{2} y^2 = h \right\}.
\end{equation}

In the decoupled coordinates, the NHIM is defined by $q_1 = 0, p_1 = 0$ on a fixed energy surface, which gives $-x + p_x + p_y = 0$ that is $p_y = x$ in the coupled coordinates, and is given by
\begin{equation}
\mathcal{M}(h) = \left\{ (x, p_x, y, p_y) \, | \, p_x = 0, p_y = x, \omega_2 x^2 - \omega_2 x y + \frac{\omega_2}{2} y^2 = h \right\}.
\label{eqn:nhim_symp_2dof}
\end{equation}

Next, in the coupled coordinates, the stable and unstable manifolds of the NHIM are given by
\begin{align}
\mathcal{W}^{\rm u}(\mathcal{M}(h)) &= \left\{ (x, p_x, y, p_y) \; | \; x = p_y, \frac{\omega_2}{2}\left( (-y + p_x + p_y)^2 + p_y^2 \right)  = h > 0 \right\} \label{eqn:umani_nsqh_2dof} \\
\mathcal{W}^{\rm s}(\mathcal{M}(h)) &= \left\{ (x, p_x, y, p_y) \; | \; x =  2p_x + p_y, \frac{\omega_2}{2}\left( ( -y + p_x + p_y)^2 + p_y^2 \right)  = h > 0 \right\} \label{eqn:smani_nsqh_2dof}
\end{align}
%


Invertible linear symplectic transformations are clearly $C^\infty$ diffeomorphisms. They not only preserve the Hamiltonian structure, but they preserve the phase space structures relevant to  reaction dynamics. In particular, Lyapunov exponents are preserved under symplectic diffeomorphisms which implies that normal hyperbolicity is preserved. The no-recrossing property of the dividing surface is a result of transversality of the Hamiltonian vector field to the dividing surface, and such transversality properties are preserved under diffeomorphisms. Hence the NHIM~\eqref{eqn:nhim_symp_2dof} also has a geometry $\mathbb{S}^1$ and the invariant manifolds~\eqref{eqn:umani_nsqh_2dof} and \eqref{eqn:smani_nsqh_2dof} are also of geometry $\mathbb{R} \times \mathbb{S}^1$.

\subsubsection{Detecting the unstable periodic orbit and its manifolds}

Now we illustrate the procedure for detecting the unstable periodic orbit, its stable and unstable manifolds by computing the Lagrangian descriptor on isoenergetic two-dimensional surfaces parametrized by a pair of coordinates of the separable quadratic Hamiltonian~\eqref{eqn:ham_symp_2dof}.

\newpoint~Isoenergetic two-dimensional surface parametrized by $(x,p_x)$~\textemdash~On a fixed energy surface, $\mathcal{H}(x, p_x, y, p_y) = h$, we consider a two-dimensional surface parametrized by $(x, p_x)$ coordinates by defining 

\begin{equation}
U_{xp_x}^+ = \left\{ (x,p_x,y,p_y) \; \vert \; y = 0, \; p_y = p_y(x,p_x,y;h) : \dot{y}(x,p_x,y,p_y) > 0 \right\} , 
\label{eqn:sos_symp_2dof_xpx}
\end{equation}
Combining the dividing surface condition $p_x = 0$, NHIM condition $p_y = x$, the directionality condition $\dot{y}(x,y,p_x,p_y;h) > 0$, gives $x > 0$. Thus, for a fixed energy, the intersection of the one-dimensional NHIM~\eqref{eqn:nhim_symp_2dof} with the isoenergetic two-dimensional surface~\eqref{eqn:sos_symp_2dof_xpx} is given by

\begin{align}
\mathcal{M}(h) \cap U^{+}_{xp_x} = \left\{ (x,p_x,y,p_y) \; \vert \; \right. & \left. \kern - \nulldelimiterspace y = 0, p_x = 0, p_y = x, \right. \nonumber \\ 
& \left.\kern - \nulldelimiterspace \omega_2 x^2 - \omega_2 x y + \frac{\omega_2}{2} y^2 = h : x > 0 \right\} \label{eqn:nhim_x_coord} 
\end{align}
which is a point on the line $p_x = 0$ and $x = \sqrt{2h/\omega_2}$, and is marked by the red cross in Fig.~\ref{fig:LD_NSQH_2dof_xpx}. This point is also identified by the minimum value of the Lagrangian descriptor at this coordinate (as shown along the one-dimensional slice at constant $p_x = 0$) and agrees with~\eqref{eqn:nhim_x_coord} to within the grid resolution used for the discretization of the two-dimensional surface. Since the NHIM is the intersection of stable and unstable manifolds, its coordinate is also a minimum and singular point in the LD values.

Next, the intersection of the unstable and stable manifolds~\eqref{eqn:smani_nsqh_2dof} with the isoenergetic two-dimensional surface $U_{xp_x}^+$ is given by

\begin{align}
\mathcal{W}^{\rm u}(\mathcal{M}(h)) \cap U_{xp_x}^+ =  \left\{ (x,p_x,y,p_y) \; \vert \; \right. & \left. \kern - \nulldelimiterspace y = 0, x = p_y, \right. \nonumber \\ 
& \left.\kern - \nulldelimiterspace \frac{\omega_2}{2}\left( ( x + p_x )^2 + x^2 \right)  = h \; : \; (\lambda + \omega_2) p_x + 2 \omega_2 x  > 0 \right\} \\
\mathcal{W}^{\rm s}(\mathcal{M}(h)) \cap U_{xp_x}^+ =  \left\{ (x,p_x,y,p_y) \; \vert \;  \right. & \left. \kern - \nulldelimiterspace y = 0, x = 2p_x + p_y, \right. \nonumber \\ 
& \left.\kern - \nulldelimiterspace \frac{\omega_2}{2}\left( ( x - p_x )^2 + ( x - 2p_x )^2 \right)  = h \; : \; (\lambda + 3 \omega_2) p_x - 2 \omega_2 x < 0  \right\}
\end{align}

where the inequalities are derived from the directionality condition $\dot{y} > 0$ for the surface. This is detected by the minima in the LD plot in Fig.~\ref{fig:LD_NSQH_2dof_xpx} and highlighted by the dashed red (unstable) and dashed white (stable) curves. Points on these manifolds are also picked up by the one-dimensional slices along the constant $p_x = -0.60, 0.30$.

\begin{figure}[!ht]
	\centering
	\subfigure[]{\includegraphics[width=0.32\textwidth]{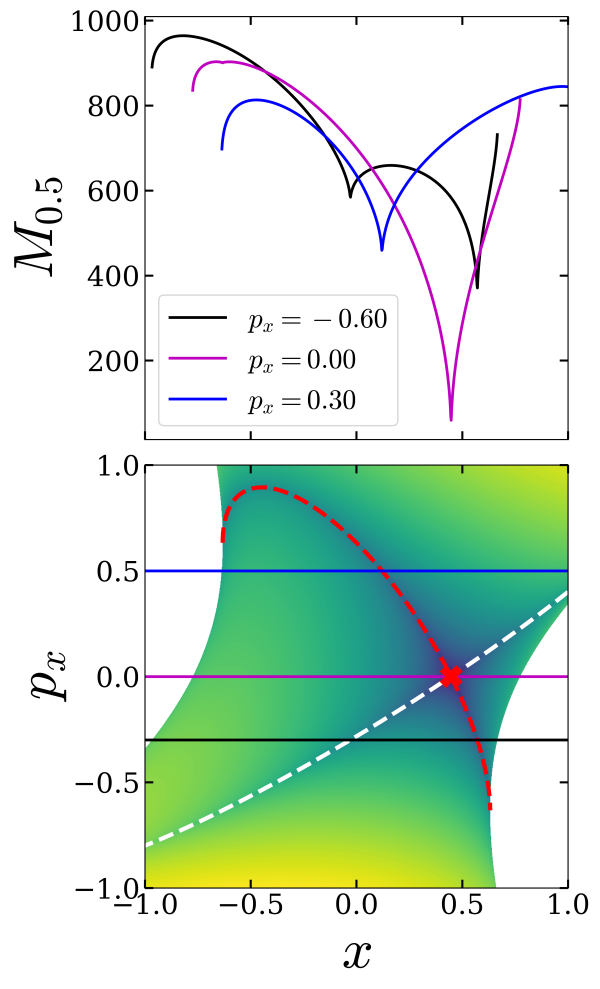}\label{fig:LD_NSQH_2dof_xpx}}
	\subfigure[]{\includegraphics[width=0.32\textwidth]{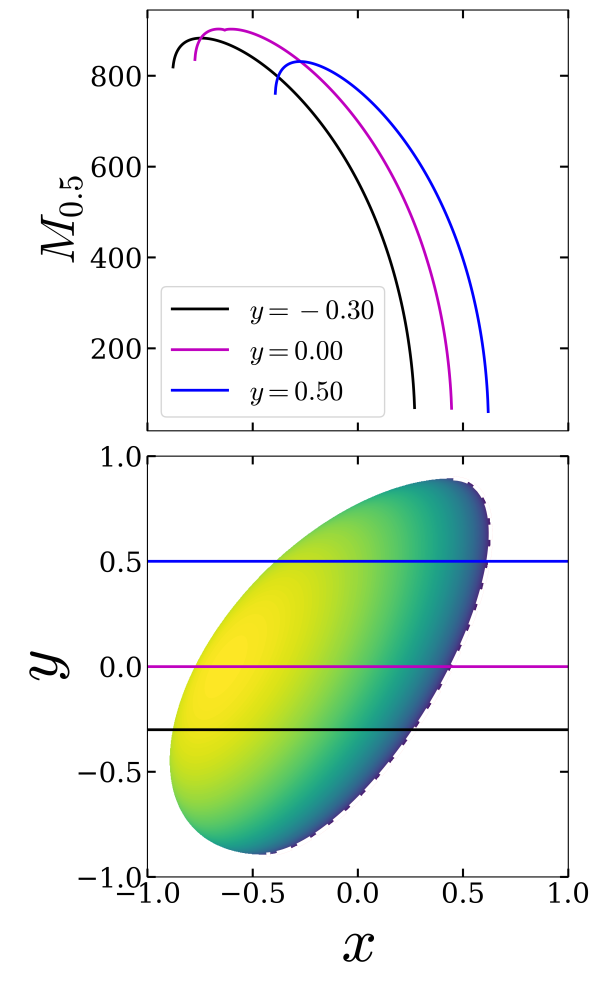}\label{fig:LD_NSQH_2dof_xy}}
	\subfigure[]{\includegraphics[width=0.32\textwidth]{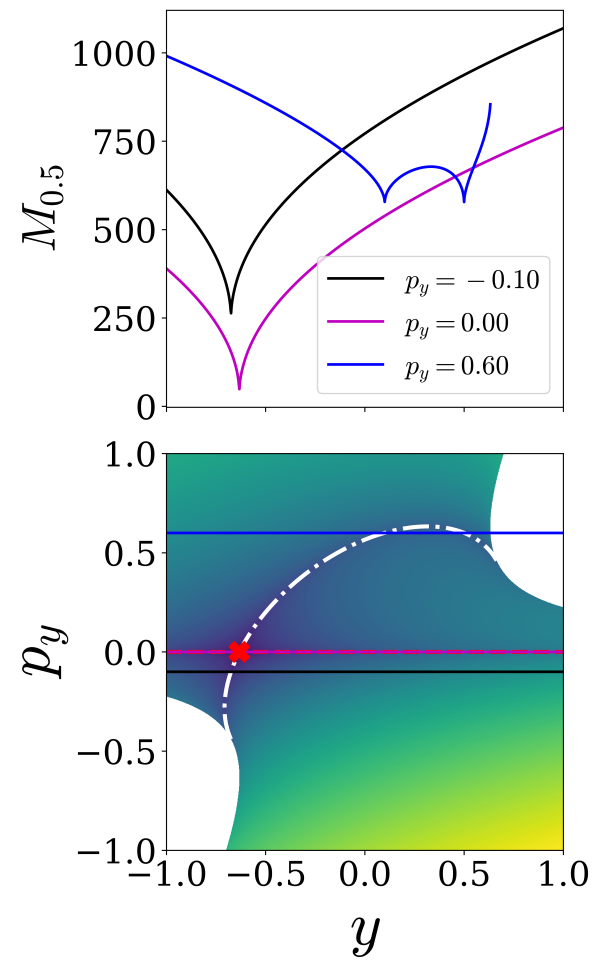}\label{fig:LD_NSQH_2dof_ypy}}
	\caption{Lagrangian descriptor plot of the non-separable quadratic Hamiltonian vector field~\eqref{eqn:eom_symp_2dof} on the isoenergetic two-dimensional surface \protect\subref{fig:LD_NSQH_2dof_xpx} $U_{xp_x}^+$, \protect\subref{fig:LD_NSQH_2dof_xy} $U_{xy}^+$, \protect\subref{fig:LD_NSQH_2dof_ypy} $U_{yp_y}^{+}$. The intersection of the NHIM and the isoenergetic two-dimensional surfaces is shown as a red cross (or dash-dot red line) and the corresponding for the manifolds is shown as dashed red (unstable) and dashed white (stable) curves.  The parameters used are $\lambda = \omega_2 = 1.0$, $h = 0.2$, and $\tau = 10$.}
	\label{fig:linear_symp_trans_2dof_M400x400_1}
\end{figure}


\newpoint~Isoenergetic two-dimensional surface parametrized by $(x,y)$~\textemdash~On a fixed energy surface $\mathcal{H}(x, p_x, y, p_y) = h$, we compute the Lagrangian descriptor on a two-dimensional surface parametrized by $(x, y)$ coordinates by defining

\begin{equation}
U_{xy}^+ = \left\{ (x,p_x,y,p_y) \; | \; p_x = 0, \; p_y(x,p_x,y;h) > 0 : \dot{p_x}(x, y, p_y) \geqslant 0 \right\} , 
\label{eqn:sos_symp_2dof_xy}
\end{equation}

Thus, the intersection of the NHIM~\eqref{eqn:nhim_symp_2dof} with this isoenergetic two-dimensional surface~\eqref{eqn:sos_symp_2dof_xy} is given by

\begin{align}
\mathcal{M}(h) \cap U_{xy}^+ =  \left\{ (x,p_x,y,p_y) \, | \, \right. & \left. \kern - \nulldelimiterspace p_y = x, p_x = 0, p_y > 0: p_y \geqslant x,  \right. \nonumber \\  
& \left.\kern - \nulldelimiterspace \omega_2 x^2 - \omega_2 x y + \frac{\omega_2}{2} y^2 = h \right\}
\end{align}

which represents a one-dimensional curve and is marked as a dashdot red line in Fig.~\ref{fig:LD_NSQH_2dof_xy}. The NHIM is also identified by the minima in Lagrangian descriptor values along one-dimensional sections at constant $y$ in Fig.~\ref{fig:LD_NSQH_2dof_xy}. 

Next, the intersection of the unstable and stable manifolds~\eqref{eqn:smani_nsqh_2dof} with the isoenergetic two-dimensional surface $U_{xy}^+$ is given by

\begin{align}
\mathcal{W}^{\rm u}(\mathcal{M}(h)) \cap U_{xy}^+ =  \left\{ (x,p_x,y,p_y) \; \vert \; \right. & \left. \kern - \nulldelimiterspace x = p_y, p_x = 0, p_y > 0: p_y \geqslant x, \right. \nonumber \\ 
& \left.\kern - \nulldelimiterspace \omega_2 x^2 - \omega_2 x y + \frac{\omega_2}{2} y^2 = h \right\} \\
\mathcal{W}^{\rm s}(\mathcal{M}(h)) \cap U_{xy}^+ =  \left\{ (x,p_x,y,p_y) \; \vert \;  \right. & \left. \kern - \nulldelimiterspace x = 2p_x + p_y, p_x = 0, p_y > 0: p_y \geqslant x, \right. \nonumber \\ 
& \left.\kern - \nulldelimiterspace \omega_2 x^2 - \omega_2 x y + \frac{\omega_2}{2} y^2 = h \right\}
\end{align}

where each manifold represents the same one-dimensional curve as the NHIM and are marked by dashed red (unstable) and dashed white (stable) lines in the Fig.~\ref{fig:LD_NSQH_2dof_xy}. These manifolds are also identified by points of minima in the Lagrangian descriptor values as shown by one-dimensional sections at constant $y$ values.


\newpoint~Isoenergetic two-dimensional surface parametrized by $(y, p_y)$~\textemdash~On the fixed energy surface $\mathcal{H}(x, p_x, y, p_y) = h$, we take the Lagrangian descriptor slice by defining a two-dimensional surface

\begin{equation}
U_{yp_y}^+ = \left\{ (x,p_x,y,p_y) \; | \; x = 0, \, p_x(x,y,p_y;h) > 0 : \dot{x}(x, p_x, y, p_y) > 0  \right\}
\label{eqn:2dsurf_nsqh_2dof_ypy}
\end{equation}

Thus, the intersection of the NHIM~\eqref{eqn:nhim_symp_2dof} with this isoenergetic two-dimensional surface~\eqref{eqn:2dsurf_nsqh_2dof_ypy} is given by

%
\begin{align}
\mathcal{M}(h) \cap U_{yp_y}^+ =  \left\{ (x,p_x,y,p_y) \; \vert \; \right. & \left. \kern - \nulldelimiterspace x = 0, p_y = x, p_x = 0, \dot{x}(x, p_x, y, p_y) > 0 : y < 0  \right. \nonumber \\ 
& \left.\kern - \nulldelimiterspace \frac{\omega_2}{2} y^2 = h \right\}  
\end{align}

which is a point $y = \sqrt{2h/\omega_2}$ on the $p_y = 0$ line and shown as a red cross in Fig.\ref{fig:LD_NSQH_2dof_ypy}. This point is also identified by the minima in the Lagrangian descriptor values as shown along one-dimensional slice in Fig.~\ref{fig:LD_NSQH_2dof_ypy}.

Next, the intersection of the unstable~\eqref{eqn:umani_nsqh_2dof} and stable manifolds~\eqref{eqn:smani_nsqh_2dof} with the isoenergetic surface $U_{yp_y}^+$ manifest as

\begin{align}
\mathcal{W}^{\rm u}(\mathcal{M}(h)) \cap U_{yp_y}^+ =  \left\{ (x,p_x,y,p_y) \; \vert \; \right. & \left. \kern - \nulldelimiterspace x = 0, x = p_y, \dot{x}(x, p_x, y, p_y) > 0 : y - p_x < 0, \right. \nonumber \\ 
& \left.\kern - \nulldelimiterspace \frac{\omega_2}{2}\left( - y + p_x \right)^2  = h \right\} \\
\mathcal{W}^{\rm s}(\mathcal{M}(h)) \cap U_{yp_y}^+ =  \left\{ (x,p_x,y,p_y) \; \vert \;  \right. & \left. \kern - \nulldelimiterspace x = 0, x = 2p_x + p_y, \right. \nonumber \\ 
& \left.\kern - \nulldelimiterspace \dot{x}(x, p_x, y, p_y) > 0 : -\omega_2y + ( \frac{\omega_2}{2} + \lambda )p_y > 0, \right. \nonumber \\ 
& \left.\kern - \nulldelimiterspace \frac{\omega_2}{2}\left( ( -y + \frac{p_y}{2} )^2 + p_y^2 \right)  = h \right\}
\end{align}

which are shown as dashed red (unstable) and dashed white (stable) curves in the Fig.~\ref{fig:LD_NSQH_2dof_ypy}. These curves are also identified by the points of minima in the Lagrangian descriptor values as evident by the one-dimensional slices in Fig.~\ref{fig:LD_NSQH_2dof_ypy}.

\section{Three Degrees of Freedom} 
\label{sec:normal_form_3dof}

In this section we use the Lagrangian Descriptor method to identify the NHIM and its stable and unstable manifolds for a 3 DoF separable quadratic Hamiltonian with an index-1 saddle point at the origin. The advantage of using this model Hamiltonian is that we can compare the analytical expressions for the phase space structures present with features in LD plots.

\subsection{Decoupled quadratic Hamiltonian}
We consider the three degrees of freedom benchmark example in Ref.\cite{wig2016} which is a linear system (quadratic Hamiltonian) with an equilibrium point of saddle-center-center equilibrium type at the origin.  The Hamiltonian of this system is of the form

\begin{equation}
H = \underbrace{\frac{\lambda}{2} \left(p_1^2 - q_1^2 \right)}_{H_r} + \underbrace{\frac{\omega_2}{2} \left(p_2^2 + q_2^2 \right)}_{H_{b_1}} + \underbrace{\frac{\omega_3}{2} \left(p_3^2 + q_3^2 \right)}_{H_{b_2}}, \quad \lambda, \, \omega_2, \, \omega_3 >0
\label{ham3}
\end{equation}

\noindent
with the corresponding Hamiltonian vector field given by:

\begin{equation}
\begin{aligned}
\dot{q}_1 & = \frac{\partial H}{\partial p_1}= \lambda p_1, \\
\dot{p}_1 & = -\frac{\partial H}{\partial q_1}= \lambda q_1, \\
\dot{q}_2 & = \frac{\partial H}{\partial p_2}= \omega_2 p_2, \\
\dot{p}_2 & = -\frac{\partial H}{\partial q_2}= -\omega_2 q_2, \\
\dot{q}_3 & = \frac{\partial H}{\partial p_3}= \omega_3 p_3, \\
\dot{p}_3 & = -\frac{\partial H}{\partial q_3}= -\omega_3 q_3, 
\end{aligned}
\label{eqn:hameq3}
\end{equation}

Since  the total Hamiltonian decouples into the independent 
Hamiltonians $H_r$, $H_{b_1}$ and $H_{b_2}$  we can analyze the phase portraits for each separately.  In 
the language of chemical reaction dynamics $H_r$ corresponds to the ``reactive mode'', $H_{b_1}$ and 
$H_{b_2}$   are ``bath modes''. The equilibrium point $(q_1, p_1, q_2, p_2, q_3, p_3)=(0, 0, 0, 0, 0, 
0)$ is an index-1 saddle point for the full three DoF system on the zero (total) energy surface. 

In this system a trajectory ``reacts'' when its $q_1$ coordinate changes sign.  Therefore  the 
surface $q_1=0$ is a {\em dividing surface} (DS) for trajectories, separating ``reactants'' from 
``products'' (i.e. before and after reaction). This DS is a  five dimensional surface in the six 
dimensional phase space. We discuss its geometrical structure, in both the phase space and in a 
fixed energy surface, in more detail.

First, we consider the ``energetics'' of the reaction. In order for $q_1$ to change sign we must 
have $H_r >0$. Also, it is clear from the form of $H_{b_1}$ and $H_{b_2}$ that $H_{b_1} \geqslant 0$ and $H_{b_2} \geqslant 0$. Therefore, reaction requires $H = H_r + H_{b_1} + H_{b_2} > 0$. 

The  energy surface in which reaction occurs  is given by

\begin{equation}
\frac{\lambda}{2} \left(p_1^2 - q_1^2 \right) + \frac{\omega_2}{2} \left(p_2^2 + q_2^2, \right)  + \frac{\omega_3}{2} \left(p_3^2 + q_3^2, \right)= H_r + H_{b_1} + H_{b_2} = H > 0, \quad H_r > 0, \, H_{b_1}, \, H_{b_2} \geqslant 0.
\label{3DoFES}
\end{equation}

\noindent
The intersection of  the DS, $q_1=0$,  with this energy surface is given by:

\begin{equation}
\frac{\lambda}{2} \, p_1^2  + \frac{\omega_2}{2} \left(p_2^2 + q_2^2, \right) + \frac{\omega_3}{2} \left(p_3^2 + q_3^2, \right)= H_r + H_{b_1} + H_{b_2} = H > 0, \quad H_r > 0, \, H_{b_1}, \, H_{b_2} \geqslant 0.
\label{3DoDS}
\end{equation}

This is the isoenergetic DS. It is a 4-sphere in the six dimensional $(q_1, p_1, q_2, p_2, q_3, p_3)$ phase space. It has two hemispheres that are DSs for the forward and backward reactions, respectively:

\begin{equation}
\frac{\lambda}{2} \, p_1^2  + \frac{\omega_2}{2} \left(p_2^2 + q_2^2, \right)  + \frac{\omega_3}{2} \left(p_3^2 + q_3^2, \right)= H_r + H_{b_1} + H_{b_2} = H > 0, \quad p_1 >0, \quad \mbox{forward DS},
\end{equation}

\begin{equation}
\frac{\lambda}{2} \, p_1^2  + \frac{\omega_2}{2} \left(p_2^2 + q_2^2, \right)  + \frac{\omega_3}{2} \left(p_3^2 + q_3^2, \right)= H_r + H_{b_1} + H_{b_2} = H > 0, \quad p_1 <0, \quad \mbox{backward DS}.
\end{equation}

Thes two hemispheres ``meet'' at  $p_1 =0$:

\begin{equation}
\frac{\omega_2}{2} \left(p_2^2 + q_2^2, \right) + \frac{\omega_3}{2} \left(p_3^2 + q_3^2, \right)= H_{b_1} + H_{b_2} \geqslant  0,  \quad  \mbox{NHIM},
\label{eqn:sep_quad_ham3dof_nhim}
\end{equation}

which is a {\em normally hyperbolic invariant 3 sphere}. {\em invariance} follows from 
\eqref{eqn:hameq3} since for $q_1 = p_1 =0$ then $\dot{q}_1 = \dot{p}_1 =0$. Thus $q_1$ and $p_1$ 
always remain zero, and   trajectories with these initial conditions remain on \eqref{eqn:sep_quad_ham3dof_nhim}, 
i.e., $q_1 = p_1 =0$ is invariant. It is {\em normally hyperbolic} since the directions normal to 
\eqref{eqn:sep_quad_ham3dof_nhim}, i.e. $q_1-p_1$, are linearized saddle like dynamics.

In order to understand the relationship between the NHIM and the index-1 saddle point we note  that for $H_r =0, \, H_{b_1} =0$ and $H_{b_2} = 0$ the NHIM reduces to the point $(q_1, p_1, q_2, p_2, q_3, p_3) = (0, 0, 0, 0, 0 ,0)$, which is the index-1 saddle point  on the energy surface $H_r + H_{b_1} + H_{b_2} =0$. Therefore  as the total energy is increased from $0$, with $H_{b_1}$   increasing from zero and/or $H_{b_2}$   increasing from zero, we see that the NHIM ``grows'' from the index-1 saddle point on the zero energy surface into an invariant 3 sphere. This shows how the ``influence'' of the index-1 saddle point is carried to higher energy sufaces on which the saddle point does not exist.

The stable and unstable manifolds of the NHIM are given the following equations:

\begin{eqnarray}
\mathcal{W}^u\left( \mathcal{M}(h) \right) & = & \left \{ (q_1, p_1, q_2, p_2, q_3, p_3) \, \vert \, q_1=p_1, \, \frac{\omega_2}{2} \left(p_2^2 + q_2^2 \right) + \frac{\omega_3}{2} \left(p_3^2 + q_3^2 \right) 
= H_{b_1} + H_{b_2} > 0 \right \}, \label{eqn:quad_ham3dof_umani} \\
\mathcal{W}^s\left( \mathcal{M}(h) \right) & = & \left \{ (q_1, p_1, q_2, p_2, q_3, p_3) \, \vert q_1= -p_1, \, \frac{\omega_2}{2} \left(p_2^2 + q_2^2 \right) + \frac{\omega_3}{2} \left(p_3^2 + q_3^2 \right) 
=H_{b_1} + H_{b_2} > 0  \right \}, \label{eqn:quad_ham3dof_smani} 
\end{eqnarray}

\noindent
and they are four dimensional on a fixed five dimensional energy surface. The topology of the manifolds is the product of a line
($q_1 = p_1$ or $q_1 = -p_1$) with a 3 sphere  
($\frac{\omega_2}{2} \left(p_2^2 + q_2^2 \right) + \frac{\omega_3}{2} \left(p_3^2 + q_3^2 \right) = H_{b_1} + H_{b_2} >0$), that is $\mathbb{R} \times \mathbb{S}^3$ and sometimes such geometry is referred to as ``spherical cylinders''. Since the  lines $q_1 = p_1$ and $q_1 = -p_1$ correspond to the contour $H_r = 0$,  in the six dimensional phase space, these manifolds have energy $H = H_r + H_{b_1}  + H_{b_2} = 0 + H_{b_1}  + H_{b_2} > 0$.

To understand the dynamics governed by the invariant manifolds of the NHIM, let us choose an initial condition $(q_1, p_1, q_2, p_2, q_3, p_3)$ on $\mathcal{W}^u\left(\mathcal{M}(h)\right)$. Then as $t \rightarrow +\infty$ the $(q_2, p_2, q_3, p_3)$ components of the trajectory with this initial condition evolve quasiperiodically and the $(q_1, p_1)$ components grow  at an exponential rate. Similarly, if we  choose an initial condition on $\mathcal{W}^s\left(\mathcal{M}(h)\right)$, then as $t \rightarrow +\infty$ the $(q_2, p_2, q_3, p_3)$ components of the trajectory  evolve quasiperiodically and the $(q_1, p_1)$ components decay to zero  at an exponential rate as $t \rightarrow \infty$.  In other words, trajectories starting on $\mathcal{W}^u \left(\mathcal{M}(h)\right)$ decay at an exponential  rate to the NHIM as $t \rightarrow - \infty$ and trajectories starting on $\mathcal{W}^s\left(\mathcal{M}(h)\right)$ decay at an exponential  rate to the NHIM as $t \rightarrow  + \infty$.


\subsubsection{Detecting NHIM and its manifolds}

We consider isoenergetic two-dimensional surfaces parametrized by two coordinates and compute the Lagrangian descriptor in a square domain of size 2 units around the origin. We discretize the coordinates of the two dimensional surface and pick constant values for three of the four remaining coordinates, and use the total energy equation to solve for the sixth coordinate. Due to the form of the Hamiltonian~\eqref{ham3}, obtaining the coordinate from the constant energy condition is simply solving a quadratic equation. 


\newpoint~Isoenergetic two-dimensional surface parametrized by $(q_1, p_1)$~\textemdash~On the constant energy surface, $H(q_1, p_1, q_2, p_2, q_3, p_3) = h$, we compute Lagrangian descriptor on a two-dimensional surface parametrized by $(q_1, p_1)$ coordinates by defining 
\begin{align}
U_{q_1p_1}^+ = \left\{ (q_1, p_1, q_2, p_2, q_3, p_3) \; | \; \right. & \left. \kern - \nulldelimiterspace q_2 = 0,  p_2 = 0, q_3 = 0, \dot{q}_3 > 0 : \right. \nonumber \\
& \left. \kern - \nulldelimiterspace p_3(q_1, p_1, q_2, p_2, q_3; h) > 0 \right\}
\end{align}
where
\begin{align}
p_3(q_1, p_1, q_2 = 0, p_2 = 0, q_3 = 0; h) = \sqrt{\frac{2}{\omega_3}\left( h - \frac{\lambda}{2}\left( p_1^2 - q_1^2 \right) \right)} 
\end{align}
The intersection of the two-dimensional surface $U_{q_1p_1}^+$ with the NHIM~\eqref{eqn:sep_quad_ham3dof_nhim} becomes

\begin{align}
\mathcal{M}(h) \cap U_{q_1p_1}^+ = \left\{ (q_1, p_1, q_2, p_2, q_3, p_3) \; \vert \; \right. & \left. \kern - \nulldelimiterspace q_1 = 0, p_1 = 0, q_2 = 0,  p_2 = 0, q_3 = 0, \dot{q}_3 > 0 : \right. \nonumber \\ 
& \left.\kern - \nulldelimiterspace  p_3(q_1, p_1, q_2, p_2, q_3; h) > 0 \right\}.
\end{align}
%
Thus, the NHIM is located at the origin $(0,0)$ and marked by a red cross in the LD plot (Fig.~\ref{fig:LD_SQH_3dof_q1p1}).  

Next, the intersection of the two-dimensional surface with the unstable~\eqref{eqn:quad_ham3dof_umani} and stable manifolds~\eqref{eqn:quad_ham3dof_smani} is given by

\begin{align}
\mathcal{W}^u(\mathcal{M}(h)) \cap U_{q_1p_1}^+ = \left\{ (q_1, p_1, q_2, p_2, q_3, p_3) \; \vert \; \right. & \left. \kern - \nulldelimiterspace q_1 = p_1, q_2 = 0, p_2 = 0, q_3 = 0, \dot{q_3} > 0 : \right. \nonumber \\  
& \left.\kern - \nulldelimiterspace p_3(q_1, p_1, q_2, p_2, q_3; h) > 0 \right\}, \\
\mathcal{W}^s(\mathcal{M}(h)) \cap U_{q_1p_1}^+ = \left\{ (q_1, p_1, q_2, p_2, q_3, p_3) \; \vert \; \right. & \left. \kern - \nulldelimiterspace q_1 = -p_1, q_2 = 0, p_2 = 0, q_3 = 0,\dot{q_3} > 0 : \right. \nonumber \\  
& \left.\kern - \nulldelimiterspace p_3(q_1, p_1, q_2, p_2, q_3;h) > 0 \right\},
\end{align}

which are one-dimensional for a fixed energy, and represent lines passing through the origin shown as dashed red (unstable) and white (stable) lines, respectively, in Fig.~\ref{fig:LD_SQH_3dof_q1p1}. The only points of local minima in the LD plot (Fig.~\ref{fig:LD_SQH_3dof_q1p1}) also lie along the lines passing through the origin and correspond to the manifolds of the NHIM.

\newpoint~Isoenergetic two-dimensional surface parametrized by $(q_2, p_2)$~\textemdash~On the constant energy surface, $H(q_1, p_1, q_2, p_2, q_3, p_3) = h$, we compute Lagrangian descriptor on a two-dimensional surface parametrized by $(q_2, p_2)$ coordinates by defining 
\begin{align}
U_{q_2p_2}^+ = \left\{ (q_1, p_1, q_2, p_2, q_3, p_3) \; | \; \right. & \left. \kern - \nulldelimiterspace q_1 = 0, q_3 = 0, p_3 = 0, \dot{q}_1 > 0 : \right. \nonumber \\
& \left. \kern - \nulldelimiterspace p_1(q_1, q_2, p_2, q_3, p_3; h) > 0 \right\}
\end{align}
where
\begin{align}
p_1(q_1 = 0, q_2, p_2, q_3 = 0, p_3 = 0; h) = \sqrt{\frac{2}{\lambda}\left( h - \frac{\omega_2}{2}\left( p_2^2 + q_2^2 \right) \right)} 
\end{align}

The intersection of the two-dimensional surface $U_{q_2p_2}^+$ with the NHIM~\eqref{eqn:sep_quad_ham3dof_nhim} becomes

\begin{align}
\mathcal{M}(h) \cap U_{q_2p_2}^+ = \left\{ (q_1, p_1, q_2, p_2, q_3, p_3) \; \vert \; \right. & \left. \kern - \nulldelimiterspace q_1 = 0, p_1 = 0, p_3 = 0, q_3 = 0, \dot{q}_3 > 0 : \right. \nonumber \\ 
& \left.\kern - \nulldelimiterspace  \frac{\omega_2}{2}\left( p_2^2 + q_2^2  \right) = h \right\}.
\end{align}

Thus, the NHIM is the circle of radius $\sqrt{2h/\omega_2}$ and marked by a dashed line in the LD plot (Fig.~\ref{fig:LD_SQH_3dof_q2p2}).

Next, the intersection of the two-dimensional surface with the unstable~\eqref{eqn:quad_ham3dof_umani} and stable manifolds~\eqref{eqn:quad_ham3dof_smani} is given by

\begin{align}
\mathcal{W}^u(\mathcal{M}(h)) \cap U_{q_2p_2}^+ = \left\{ (q_1, p_1, q_2, p_2, q_3, p_3) \; \vert \; \right. & \left. \kern - \nulldelimiterspace q_1 = p_1, q_1 = 0, q_3 = 0, p_3 = 0, \dot{q_3} > 0 : \right. \nonumber \\  
& \left.\kern - \nulldelimiterspace \frac{\omega_2}{2}\left( p_2^2 + q_2^2  \right) = h \right\}, \\
\mathcal{W}^s(\mathcal{M}(h)) \cap U_{q_2p_2}^+ = \left\{ (q_1, p_1, q_2, p_2, q_3, p_3) \; \vert \; \right. & \left. \kern - \nulldelimiterspace q_1 = -p_1, q_1 = 0, q_3 = 0, p_3 = 0,\dot{q_3} > 0 : \right. \nonumber \\  
& \left.\kern - \nulldelimiterspace \frac{\omega_2}{2}\left( p_2^2 + q_2^2  \right) = h \right\},
\end{align}

which are one-dimensional for a fixed energy, and marked by dashed red (unstable) and white (stable) lines, respectively, in Fig.~\ref{fig:LD_SQH_3dof_q2p2}. The only points of local minima in the LD plot (Fig.~\ref{fig:LD_SQH_3dof_q2p2}) also lie along the lines passing through the origin and correspond to the manifolds of the NHIM.

\begin{figure}[!ht] 
	\centering
	\subfigure[]{\includegraphics[width=0.32\textwidth]{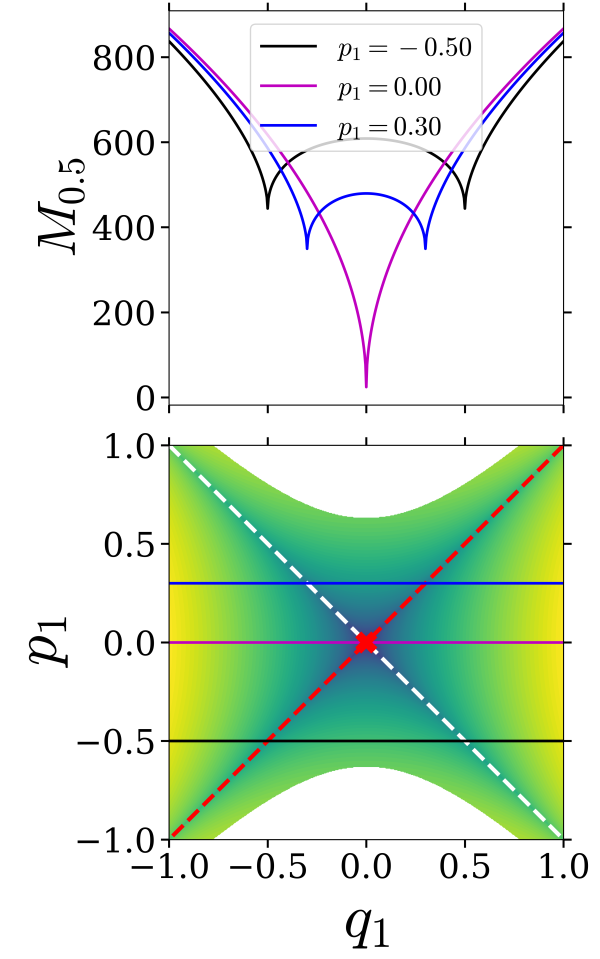}\label{fig:LD_SQH_3dof_q1p1}}
	\subfigure[]{\includegraphics[width=0.32\textwidth]{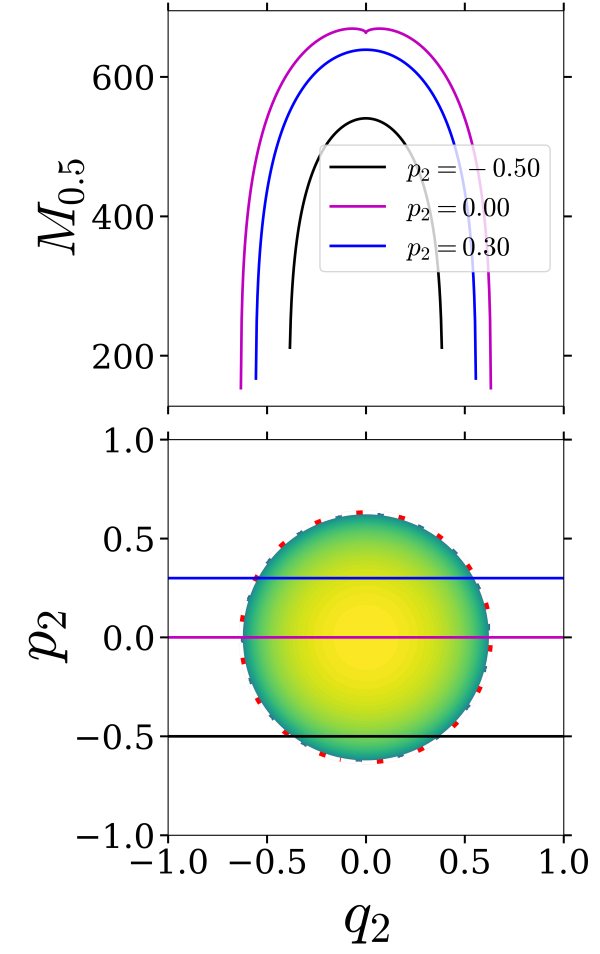}\label{fig:LD_SQH_3dof_q2p2}}
	\subfigure[]{\includegraphics[width=0.32\textwidth]{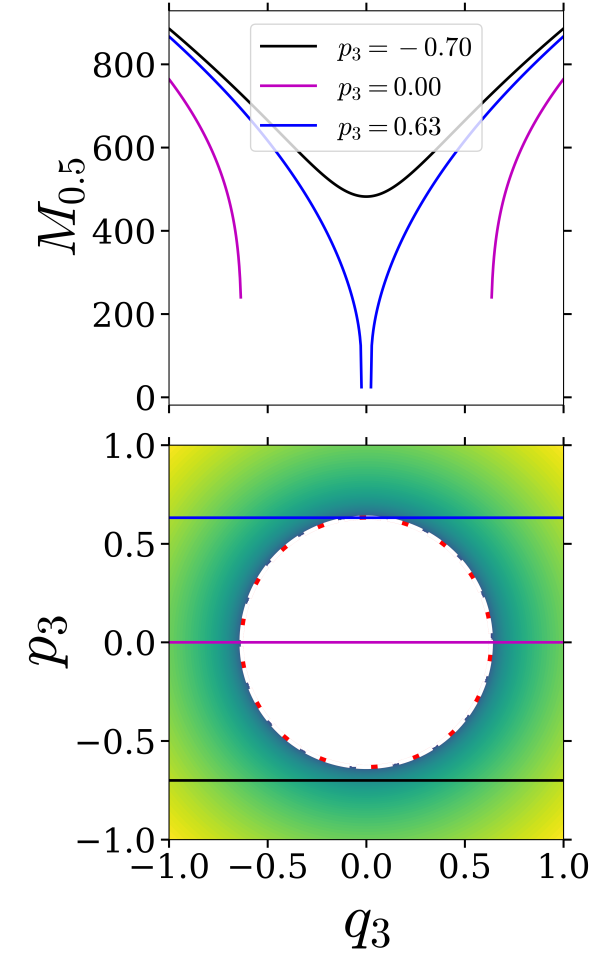}\label{fig:LD_SQH_3dof_q3p3}}
	\caption{Lagrangian descriptor plot of the separable quadratic Hamiltonian vector field~\eqref{eqn:hameq3} on the isoenergetic two-dimensional surface \protect\subref{fig:LD_SQH_3dof_q1p1} $U_{q_1p_1}^{+}$, \protect\subref{fig:LD_SQH_3dof_q2p2} $U_{q_2p_2}^{+}$, \protect\subref{fig:LD_SQH_3dof_q3p3} $U_{q_3p_3}^{+}$. The parameters used are $\lambda = \omega_2 = \omega_3 = 1.0$, $h = 0.2$, and $\tau = 10$.}
\end{figure}


\newpoint~Isoenergetic two-dimensional surface parametrized by $(q_3, p_3)$~\textemdash~On the constant energy surface, $H(q_1, p_1, q_2, p_2, q_3, p_3) = h$, we compute Lagrangian descriptor on a two-dimensional surface parametrized by $(q_3, p_3)$ coordinates by defining 

\begin{align}
U_{q_3p_3}^+ = \left\{ (q_1, p_1, q_2, p_2, q_3, p_3) \; | \; \right. & \left. \kern - \nulldelimiterspace p_1 = 0, q_2 = 0, p_2 = 0, \dot{p_1} > 0 : \right. \nonumber \\
& \left. \kern - \nulldelimiterspace q_1(p_1, q_2, p_2, q_3, p_3; h) > 0 \right\}
\end{align}

where

\begin{align}
q_1(p_1 = 0, q_2 = 0, p_2 = 0, q_3, p_3; h) = \sqrt{\frac{2}{\lambda}\left( \frac{\omega_3}{2}\left( p_3^2 + q_3^2 \right) - h \right)} 
\end{align}

The intersection of the two-dimensional surface $U_{q_3p_3}^+$ with the NHIM~\eqref{eqn:sep_quad_ham3dof_nhim} becomes

\begin{align}
\mathcal{M}(h) \cap U_{q_3p_3}^+ = \left\{ (q_1, p_1, q_2, p_2, q_3, p_3) \; \vert \; \right. & \left. \kern - \nulldelimiterspace q_1 = 0, p_1 = 0, p_2 = 0, q_2 = 0, \dot{p_1} > 0 : \right. \nonumber \\ 
& \left.\kern - \nulldelimiterspace  \frac{\omega_3}{2}\left( p_3^2 + q_3^2  \right) = h \right\}.
\end{align}

Thus, the NHIM is the circle of radius $\sqrt{2h/\omega_3}$ and marked by a dashed line in the LD plot (Fig.~\ref{fig:LD_SQH_3dof_q3p3}).  

Next, the intersection of the two-dimensional surface with the unstable~\eqref{eqn:quad_ham3dof_umani} and stable manifolds~\eqref{eqn:quad_ham3dof_smani} is given by

\begin{align}
\mathcal{W}^u(\mathcal{M}(h)) \cap U_{q_3p_3}^+ = \left\{ (q_1, p_1, q_2, p_2, q_3, p_3) \; \vert \; \right. & \left. \kern - \nulldelimiterspace q_1 = p_1, p_1 = 0, q_2 = 0, p_2 = 0, \dot{p_1} > 0 : \right. \nonumber \\  
& \left.\kern - \nulldelimiterspace \frac{\omega_3}{2}\left( p_3^2 + q_3^2  \right) = h \right\}, \\
\mathcal{W}^s(\mathcal{M}(h)) \cap U_{q_3p_3}^+ = \left\{ (q_1, p_1, q_2, p_2, q_3, p_3) \; \vert \; \right. & \left. \kern - \nulldelimiterspace q_1 = -p_1, p_1 = 0, q_2 = 0, p_2 = 0,\dot{p_1} > 0 : \right. \nonumber \\  
& \left.\kern - \nulldelimiterspace \frac{\omega_3}{2}\left( p_3^2 + q_3^2  \right) = h \right\},
\end{align}

which are one-dimensional for a fixed energy and represent circles of radius $\sqrt{2h/\omega_3}$, and marked by dashed red (unstable) and white (stable) lines, respectively, in Fig.~\ref{fig:LD_SQH_3dof_q3p3}. The only points of minima and singularity in the LD plot (Fig.~\ref{fig:LD_SQH_3dof_q3p3}) is along the circle and thus identify the manifolds of the NHIM.

\newpoint~Isoenergetic two-dimensional surface parametrized by $(q_1, q_2)$~\textemdash~On the constant energy surface, $H(q_1, p_1, q_2, p_2, q_3, p_3) = h$, we compute Lagrangian descriptor on a two-dimensional surface parametrized by $(q_1, q_2)$ coordinates by defining 
\begin{align}
U_{q_1q_2}^+ = \left\{ (q_1, p_1, q_2, p_2, q_3, p_3) \; | \; \right. & \left. \kern - \nulldelimiterspace p_1 = 0, p_2 = 0, q_3 = 0, \dot{q}_3 > 0 : \right. \nonumber \\
& \left. \kern - \nulldelimiterspace p_3(q_1, p_1, q_2, p_2, q_3; h) > 0 \right\}
\end{align}
where
\begin{align}
p_3(q_1, p_1 = 0, q_2, p_2 = 0, q_3 = 0; h) = \sqrt{\frac{2}{\omega_3}\left( h - \left( \frac{\omega_2}{2} q_2^2 - \frac{\lambda}{2} q_1^2 \right)  \right)} 
\end{align}

The intersection of the two-dimensional surface $U_{q_1q_2}^+$ with the NHIM~\eqref{eqn:sep_quad_ham3dof_nhim} is given by

\begin{align}
\mathcal{M}(h) \cap U_{q_1q_2}^+ = \left\{ (q_1, p_1, q_2, p_2, q_3, p_3) \; \vert \; \right. & \left. \kern - \nulldelimiterspace q_1 = 0, p_1 = 0, p_2 = 0, q_3 = 0, \right. \nonumber \\ 
& \left.\kern - \nulldelimiterspace \frac{\omega_2}{2} q_2^2 + \frac{\omega_3}{2} p_3^2  = h : p_3(q_1, p_1, q_2, p_2, q_3; h) > 0 \right\}.
\end{align}

Thus, on the isoenergetic two-dimensional surface $U_{q_1q_2}^+$, the NHIM is the line $q_1 = 0$ and marked by a dashdot line in the LD plot (Fig.~\ref{fig:LD_SQH_3dof_q1q2}).  

Next, the intersection of the unstable~\eqref{eqn:quad_ham3dof_umani} and stable manifolds~\eqref{eqn:quad_ham3dof_smani} with two-dimensional surface is given by

\begin{align}
\mathcal{W}^u(\mathcal{M}(h)) \cap U_{q_1q_2}^+ = \left\{ (q_1, p_1, q_2, p_2, q_3, p_3) \; \vert \; \right. & \left. \kern - \nulldelimiterspace q_1 = p_1, p_1 = 0, p_2 = 0, q_3 = 0, \right. \nonumber \\  
& \left.\kern - \nulldelimiterspace \frac{\omega_2}{2} q_2^2 + \frac{\omega_3}{2} p_3^2  = h :  p_3(q_1, p_1, q_2, p_2, q_3; h) > 0 \right\}, \\
\mathcal{W}^s(\mathcal{M}(h)) \cap U_{q_1q_2}^+ = \left\{ (q_1, p_1, q_2, p_2, q_3, p_3) \; \vert \; \right. & \left. \kern - \nulldelimiterspace q_1 = -p_1, p_1 = 0, p_2 = 0, q_3 = 0, \right. \nonumber \\  
& \left.\kern - \nulldelimiterspace \frac{\omega_2}{2} q_2^2 + \frac{\omega_3}{2} p_3^2  = h :  p_3(q_1, p_1, q_2, p_2, q_3; h) > 0 \right\},
\end{align}

which are one-dimensional for a fixed energy, and marked by dashed red (unstable) and white (stable) lines, respectively, in Fig.~\ref{fig:LD_SQH_3dof_q1q2}. The only points of minima in the LD plot (Fig.~\ref{fig:LD_SQH_3dof_q1q2}) also lie along this line at $q_1 = 0$ and identify the manifolds of the NHIM.

\newpoint~Isoenergetic two-dimensional surface parametrized by $(q_2, q_3)$~\textemdash~On the constant energy surface, $H(q_1, p_1, q_2, p_2, q_3, p_3) = h$, we compute Lagrangian descriptor on a two-dimensional surface parametrized by $(q_2, q_3)$ coordinates by defining 
\begin{align}
U_{q_2q_3}^+ = \left\{ (q_1, p_1, q_2, p_2, q_3, p_3) \; | \; \right. & \left. \kern - \nulldelimiterspace q_1 = 0, p_2 = 0, p_3 = 0, \dot{q}_1 \geqslant 0 : \right. \nonumber \\
& \left. \kern - \nulldelimiterspace p_1(q_1, q_2, p_2, q_3, p_3; h) \geqslant 0 \right\}
\end{align}
where
\begin{align}
p_1(q_1 = 0, q_2, p_2 = 0, q_3, p_3 = 0; h) = \sqrt{\frac{2}{\lambda}\left( h - \left( \frac{\omega_2}{2} q_2^2 + \frac{\omega_3}{2} q_3^2 \right)  \right)} 
\end{align}

The intersection of the two-dimensional surface $U_{q_2q_3}^+$ with the NHIM~\eqref{eqn:sep_quad_ham3dof_nhim} is given by

\begin{align}
\mathcal{M}(h) \cap U_{q_2q_3}^+ = \left\{ (q_1, p_1, q_2, p_2, q_3, p_3) \; \vert \; \right. & \left. \kern - \nulldelimiterspace q_1 = 0, p_1 = 0, p_2 = 0, p_3 = 0, \right. \nonumber \\ 
& \left.\kern - \nulldelimiterspace \frac{\omega_2}{2} q_2^2 + \frac{\omega_3}{2} q_3^2  = h, p_1(q_1, q_2, p_2, q_3, p_3; h) \geqslant 0  \right\}.
\end{align}

Thus, on the isoenergetic two-dimensional surface $U_{q_2q_3}^+$, the NHIM is an ellipse of semi-major axis $\sqrt{2h/\omega_2}$ and semi-minor axis $\sqrt{2h/\omega_3}$ if $\omega_2 < \omega_3$, and vice-versa, otherwise. In our case, $\omega_2 = \omega_3$, the intersection of the NHIM with the isoenergetic two-dimensional surface becomes a circle and marked by a dashdot line in the LD plot (Fig.~\ref{fig:LD_SQH_3dof_q2q3}).  

Next, the intersection of the unstable~\eqref{eqn:quad_ham3dof_umani} and stable manifolds~\eqref{eqn:quad_ham3dof_smani} with the isoenergetic two-dimensional surface is given by

\begin{align}
\mathcal{W}^u(\mathcal{M}(h)) \cap U_{q_2q_3}^+ = \left\{ (q_1, p_1, q_2, p_2, q_3, p_3) \; \vert \; \right. & \left. \kern - \nulldelimiterspace q_1 = p_1, q_1 = 0, p_2 = 0, p_3 = 0, \right. \nonumber \\  
& \left.\kern - \nulldelimiterspace \frac{\omega_2}{2} q_2^2 + \frac{\omega_3}{2} q_3^2  = h, p_1(q_1, q_2, p_2, q_3, p_3; h) \geqslant 0  \right\}, \\
\mathcal{W}^s(\mathcal{M}(h)) \cap U_{q_2q_3}^+ = \left\{ (q_1, p_1, q_2, p_2, q_3, p_3) \; \vert \; \right. & \left. \kern - \nulldelimiterspace q_1 = -p_1, q_1 = 0, p_2 = 0, p_3 = 0, \right. \nonumber \\  
& \left.\kern - \nulldelimiterspace \frac{\omega_2}{2} q_2^2 + \frac{\omega_3}{2} q_3^2  = h, p_1(q_1, q_2, p_2, q_3, p_3; h) \geqslant 0  \right\},
\end{align}

which have the same geometry as the NHIM on the isoenergetic two-dimensional surface, $U_{q_2q_3}$, and marked by dashed red (unstable) and white (stable) lines, respectively, in Fig.~\ref{fig:LD_SQH_3dof_q2q3}. The only points of minima in the LD plot (Fig.~\ref{fig:LD_SQH_3dof_q2q3}) also lie along this circle and identify the manifolds of the NHIM. 

\begin{figure}[!ht]
	\centering
	\subfigure[]{\includegraphics[width=0.32\textwidth]{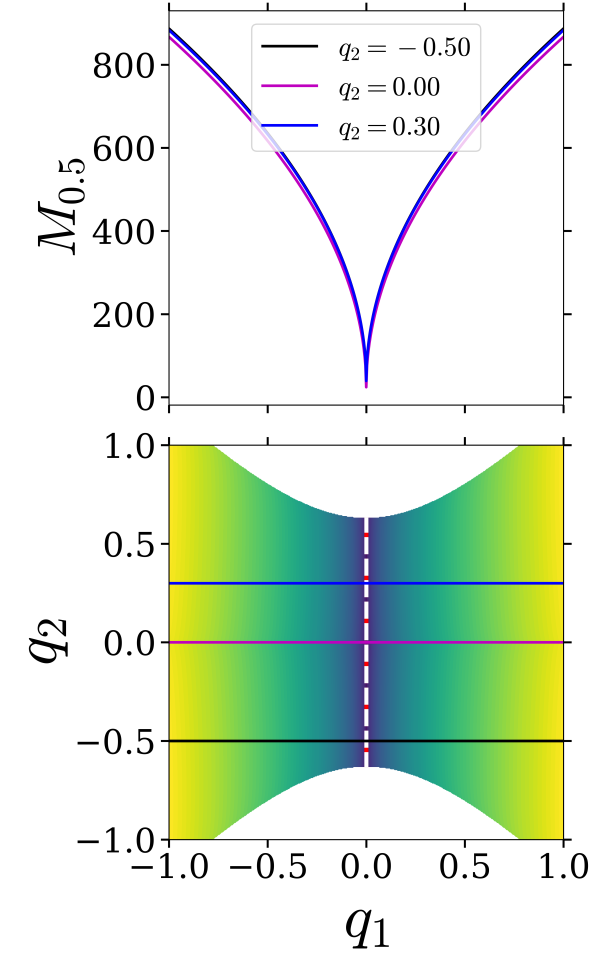}\label{fig:LD_SQH_3dof_q1q2}}
	\subfigure[]{\includegraphics[width=0.32\textwidth]{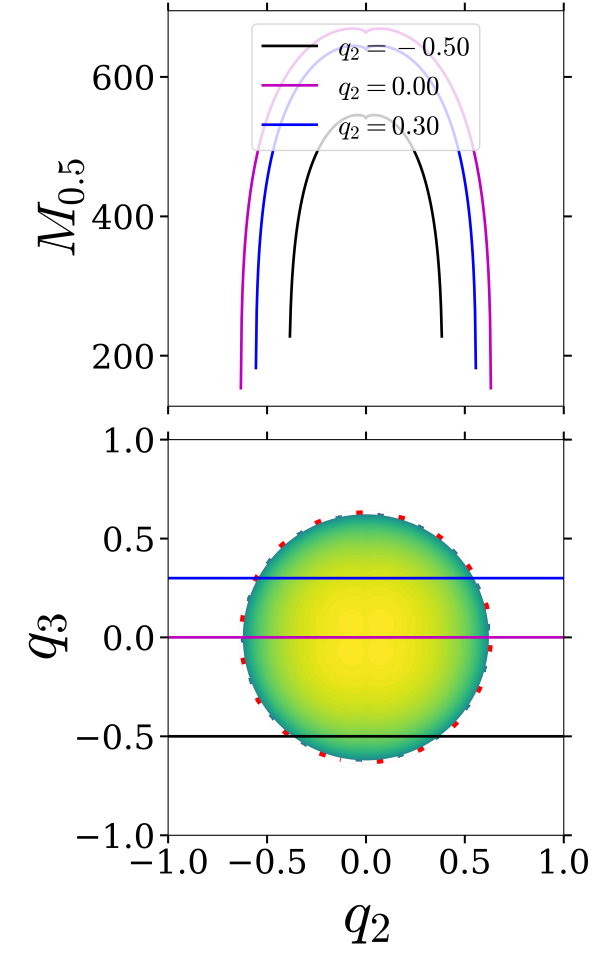}\label{fig:LD_SQH_3dof_q2q3}}
	\subfigure[]{\includegraphics[width=0.32\textwidth]{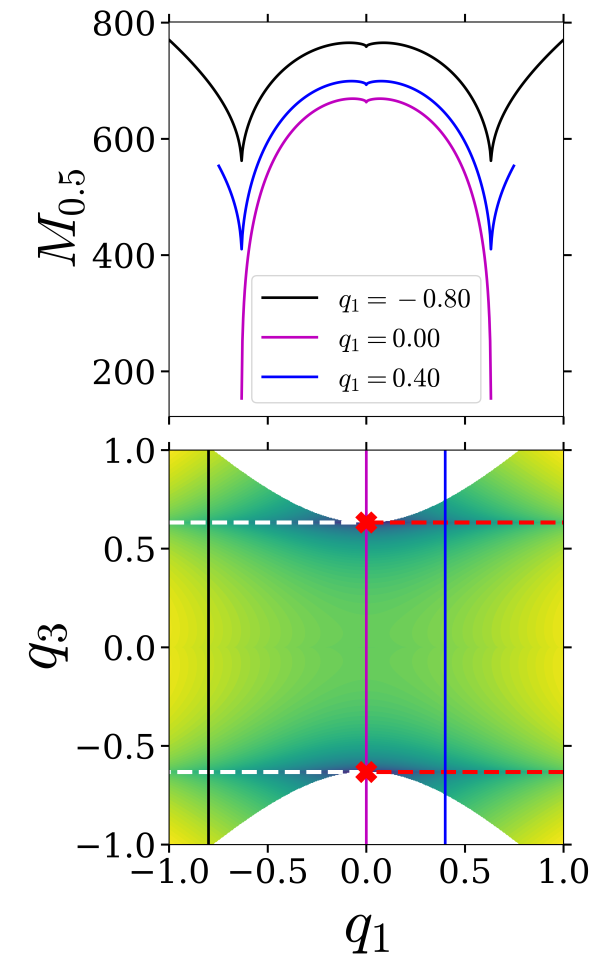}\label{fig:LD_SQH_3dof_q3q1}}
	\caption{Lagrangian descriptor plot of the separable quadratic Hamiltonian vector field~\eqref{eqn:hameq3} on the isoenergetic two-dimensional surface \protect\subref{fig:LD_SQH_3dof_q1q2} $U_{q_1q_2}^{+}$, \protect\subref{fig:LD_SQH_3dof_q2q3} $U_{q_2q_3}^{+}$, \protect\subref{fig:LD_SQH_3dof_q3q1} $U_{q_3q_1}^{+}$. The parameters used are $\lambda = \omega_2 = \omega_3 = 1.0$, $h = 0.2$, and $\tau = 10$.}
\end{figure}

\newpoint~Isoenergetic two-dimensional surface parametrized by $(q_3, q_1)$~\textemdash~On the constant energy surface, $H(q_1, p_1, q_2, p_2, q_3, p_3) = h$, we compute Lagrangian descriptor on a two-dimensional surface parametrized by $(q_3, q_1)$ coordinates by defining

%
%
%
%

\begin{align}
U_{q_3q_1}^+ = \left\{ (q_1, p_1, q_2, p_2, q_3, p_3) \; | \; \right. & \left. \kern - \nulldelimiterspace q_2 = 0, p_2 = 0, p_3 = 0, \dot{q_1} > 0 : \right. \nonumber \\
& \left. \kern - \nulldelimiterspace p_1(q_1, q_2, p_2, q_3, p_3; h) > 0 \right\}
\end{align}

where

\begin{align}
p_1(q_1, q_2 = 0, p_2 = 0, q_3, p_3 = 0; h) = \sqrt{\frac{2}{\lambda}\left( h - \left( \frac{\omega_3}{2} q_3^2 - \frac{\lambda}{2} q_1^2 \right)  \right)} 
\end{align}

The intersection of the two-dimensional surface $U_{q_3q_1}^+$ with the NHIM~\eqref{eqn:sep_quad_ham3dof_nhim} is given by

\begin{align}
\mathcal{M}(h) \cap U_{q_3q_1}^+ = \left\{ (q_1, p_1, q_2, p_2, q_3, p_3) \; \vert \; \right. & \left. \kern - \nulldelimiterspace q_1 = 0, p_1 = 0, q_2 = 0, p_2 = 0, p_3 = 0, \right. \nonumber \\ 
& \left.\kern - \nulldelimiterspace p_1(q_1, q_2, p_2, q_3, p_3; h) > 0, \frac{\omega_3}{2} q_3^2  = h \right\}.
\end{align}

which represents two points $q_3 = \pm \sqrt{2h/\omega_3}$ on the line $q_1 = 0$, and marked by a red cross in Fig.~\ref{fig:LD_SQH_3dof_q3q1}. These points are also identified by the minima in the Lagrangian descriptor values as shown by one-dimensional slices at constant $q_1$.

Next, the intersection of the unstable~\eqref{eqn:quad_ham3dof_umani} and the stable manifold~\eqref{eqn:quad_ham3dof_smani} with the isoenergetic two-dimensional surface is given by

\begin{align}
\mathcal{W}^u(\mathcal{M}(h)) \cap U_{q_3q_1}^+ = \left\{ (q_1, p_1, q_2, p_2, q_3, p_3) \; \vert \; \right. & \left. \kern - \nulldelimiterspace q_1 = p_1, q_2 = 0, p_2 = 0, p_3 = 0, \right. \nonumber \\  
& \left.\kern - \nulldelimiterspace p_1(q_1, q_2, p_2, q_3, p_3; h) > 0, \frac{\omega_3}{2} q_3^2  = h \right\}, \\
\mathcal{W}^s(\mathcal{M}(h)) \cap U_{q_3q_1}^+ = \left\{ (q_1, p_1, q_2, p_2, q_3, p_3) \; \vert \; \right. & \left. \kern - \nulldelimiterspace q_1 = -p_1, q_2 = 0, p_2 = 0, p_3 = 0, \right. \nonumber \\  
& \left.\kern - \nulldelimiterspace p_1(q_1, q_2, p_2, q_3, p_3; h) > 0, \frac{\omega_3}{2} q_3^2  = h \right\},
\end{align}

where each manifold represent lines parallel to $q_1$ axis. The unstable manifold lies on the $q_1 > 0$ plane and the stable manifold lies on the $q_1 < 0$ plane, and are marked by dashed red (unstable) and dashed white (stable) lines in the Fig.~\ref{fig:LD_SQH_3dof_q3q1}. These manifolds are again identified by points of minima in the Lagrangian descriptor values as shown by one-dimensional slices at constant $q_1$.

\subsection{Coupled quadratic Hamiltonian} 

To couple the coordinates in the separable quadratic Hamiltonian~\eqref{ham3}, we introduce a linear transformation such that it satisfies the symplectic condition~\eqref{eqn:symp_cond}. 

Let us consider the symplectic transformation
\begin{equation}
C =
\begin{pmatrix}
0 & 0 & 0 & 1 & 0 & 0 \\
0 & 0 & 0 & 0 & 1 & 0 \\
0 & 0 & 0 & 0 & 0 & 1 \\
-1 & 0 & 0 & 1 & 1 & 1 \\
0 & -1 & 0 & 1 & 1 & 1 \\
0 & 0 & -1 & 1 & 1 & 1 \\
\end{pmatrix}
\label{eqn:three_dof_C}
\end{equation}

The change of coordinates is given by
\begin{equation}
\begin{aligned}
q_1 &= p_x \\
p_1 &= -x + p_x + p_y + p_z \\
q_2 &= p_y \\
p_2 &= -y + p_x + p_y + p_z \\
q_3 &= p_z \\
p_3 &= -z + p_x + p_y + p_z
\end{aligned}
\end{equation}
This transformation applied to~\eqref{ham3} gives the Hamiltonian in the transformed 
coordinates
\begin{equation}
\begin{aligned}
\mathcal{H}(x, p_x, y, p_y, z, p_z) = \frac{\lambda}{2} \left[ (-x + p_x + p_y + p_z)^2 -  p_x^2 \right] + 
& 
\frac{\omega_2}{2} \left[ (-y + p_x + p_y + p_z)^2 + p_y^2 \right] \\ + &\frac{\omega_3}{2} \left[ (-z + p_x + p_y + p_z)^2 + p_z^2 \right] 
\label{eqn:ham_symp_3dof}
\end{aligned}
\end{equation}
which gives the vector field
\begin{equation}
\begin{aligned}
\dot{x} & = \frac{\partial \mathcal{H}}{\partial p_x} = \frac{\lambda}{2} \left[ 2(-x + p_x + 
p_y + p_z) - 2p_x \right] +  \frac{\omega_2}{2} \left[ 2(-y + p_x + p_y + p_z) \right] + 
\frac{\omega_3}{2} \left[ 2(-z + p_x + p_y + p_z) \right] \\
& = \lambda\left( -x + p_y + p_z \right) +  \omega_2 \left(-y + p_x + p_y + p_z \right) + 
\omega_3 \left(-z + p_x + p_y + p_z \right) \\
\dot{p_x} & = - \frac{\partial \mathcal{H}}{\partial x} \\ & = \lambda \left ( -x + p_x + p_y + p_z 
\right) \\
\dot{y} & = \frac{\partial \mathcal{H}}{\partial p_y} = \frac{\lambda}{2} \left[ 2(-x + p_x + 
p_y + p_z) \right] +  \frac{\omega_2}{2} \left[ 2p_y + 2(-y + p_x + p_y + p_z) \right] + 
\frac{\omega_3}{2} \left[ 2(-z + p_x + p_y + p_z) \right] \\
& = \lambda \left(-x + p_x + p_y + p_z \right) +  \omega_2 \left (-y + p_x + 2p_y + p_z \right) + 
\omega_3 \left(-z + p_x + p_y + p_z \right) \\
\dot{p_y} & = - \frac{\partial \mathcal{H}}{\partial y} \\ & = \omega_2 \left( -y + p_x + p_y + p_z 
\right) \\
\dot{z} & = \frac{\partial \mathcal{H}}{\partial p_z} = \frac{\lambda}{2} \left[ 2(-x + p_x + 
p_y + p_z) \right] +  \frac{\omega_2}{2} \left[ 2(-y + p_x + p_y + p_z) \right] + 
\frac{\omega_3}{2} \left[ 2p_z + 2(-z + p_x + p_y + p_z) \right] \\
& = \lambda \left(-x + p_x + p_y + p_z \right) +  \omega_2 \left(-y + p_x + p_y + p_z \right) + 
\omega_3 \left( -z + p_x + p_y + 2p_z \right) \\
\dot{p_z} & = - \frac{\partial \mathcal{H}}{\partial z} \\ & = \omega_3 \left( -z + p_x + p_y + p_z 
\right)
\end{aligned}
\label{eqn:eom_symp_3dof}
\end{equation}
where the equilibrium point is at $(0,0,0,0,0,0)$ and its total energy is $0$. The Jacobian at this equilibrium point has eigenvalues $\lambda, -\lambda, i\omega_2, -i \omega_2, i\omega_3, -i \omega_3$ (as shown in the Appendix~\ref{sect:linear_symp_transf_jac}) and is saddle $\times$ center $\times$ center, thus index-1. The energy surface is a five dimensional surface given by~\eqref{eqn:ham_symp_3dof} in the six dimensional phase 
space.

In the decoupled (separable) system~\eqref{ham3}, the dividing surface is defined by $q_1 = 0$ which becomes $p_x = 0$ in the coupled (non-separable) coordinates. The dividing surface for a fixed energy $\mathcal{H}(x,p_x,y,p_y,z,p_z) = h$ is given by
\begin{equation}
{\rm DS} =  \frac{\lambda}{2} \left[ \left( -x + p_y + p_z \right)^2 \right] + \frac{\omega_2}{2} \left[ p_y^2 
+ \left( -y + p_y + p_z \right)^2  \right] + \frac{\omega_3}{2} \left[ p_z^2 + \left( -z + p_y + 
p_z \right)^2  \right] = h, 
\end{equation}

On this dividing surface, the NHIM is defined by $p_1 = 0$ in the separable coordinates, which 
gives $-x + p_x + p_y  + p_z= 0$ that is $p_y + p_z = x$ in the non-separable coordinates, and can be 
expressed as
\begin{equation}
\mathcal{M}(h) = \frac{\omega_2}{2} \left[ (x - p_z)^2 + \left( x - y \right)^2 \right] + \frac{\omega_3}{2} \left[ p_z^2 
+ \left( x - z \right)^2 \right] = h.
\label{eqn:nhim_symp_3dof}
\end{equation}

Next, in the coupled coordinates, the stable and unstable manifolds of the NHIM are given by
\begin{align}
\mathcal{W}^{\rm u}(\mathcal{M}(h)) = \left\{  (x, y, z, p_x, p_y, p_z) \; | \; \right. & \left. \kern - \nulldelimiterspace x = p_y + p_z,  \right. \nonumber \\ & \left.\kern - \nulldelimiterspace \;  \frac{\omega_2}{2} \left( (-y + p_x + p_y + p_z )^2 + p_y^2 \right) + \right. \nonumber \\
&\mathrel{\phantom{=}} \left.\kern - \nulldelimiterspace \frac{\omega_3}{2}\left( (-z + p_x + p_y + p_z )^2 + p_z^2 \right)  = h \right\}, \label{eqn:umani_nsqh_3dof}\\
\mathcal{W}^{\rm s}(\mathcal{M}(h)) = \left\{  (x, y, z, p_x, p_y, p_z) \; | \; \right. & \left. \kern - \nulldelimiterspace x = 2p_x + p_y + p_z, \right. \nonumber \\ &  \left.\kern - \nulldelimiterspace  \;  \frac{\omega_2}{2} \left( (-y + p_x + p_y + p_z )^2 + p_y^2 \right) + \right. \nonumber \\
&\mathrel{\phantom{=}} \left.\kern - \nulldelimiterspace \frac{\omega_3}{2}\left( (-z + p_x + p_y + p_z )^2 + p_z^2 \right)  = h  \right\} \label{eqn:smani_nsqh_3dof}
\end{align}

As noted earlier, invertible linear symplectic transformations preserve the normal hyperbolicity and the transversality of the Hamiltonian vector field to the dividing surface. Hence the NHIM~\eqref{eqn:nhim_symp_3dof} still has a geometry $\mathbb{S}^3$ and the invariant manifolds~\eqref{eqn:umani_nsqh_3dof} and \eqref{eqn:smani_nsqh_3dof} are also $\mathbb{R} \times \mathbb{S}^3$. 

\subsubsection{Detecting NHIM and its manifolds}

Now we illustrate the procedure for detecting NHIM and its stable and unstable manifolds using features in Lagrangian descriptor on isoenergetic two-dimensional surfaces for the non-separable quadratic Hamiltonian~\eqref{eqn:ham_symp_2dof}.


\newpoint~Isoenergetic two-dimensional surface parametrized by $(x, p_x)$~\textemdash~On the constant energy surface, $\mathcal{H}(x, p_x, y, p_y, z, p_z) = h$~\eqref{eqn:ham_symp_3dof}, we compute the Lagrangian descriptor on a two-dimensional surface parametrized by $(x, p_x)$ coordinates by defining 

\begin{align}
U_{xp_x}^+ = \left\{ (x, p_x, y, p_y, z, p_z)  \; \vert \; \right. & \left. \kern - \nulldelimiterspace y = 0, z = 0, p_y = 0, \right. \nonumber \\ 
& \left.\kern - \nulldelimiterspace p_z(x, p_x, y, p_y, z; h) > 0 \, : \, \dot{z}(x, y, z, p_x, p_y, p_z) > 0 \right\}, 
\label{eqn:sos_symp_3dof_xpx}
\end{align}
%

Thus, on the five-dimensional energy surface, the intersection of the three-dimensional NHIM~\eqref{eqn:nhim_symp_3dof} with the two-dimensional surface is zero-dimensional and given by

\begin{align}
\mathcal{M}(h) \cap U_{xp_x}^+ = \left\{ (x, p_x, y, p_y, z, p_z) \; \vert \; \right. & \left. \kern - \nulldelimiterspace y = 0, z = 0, p_y = 0, p_x = 0, p_y + p_z = x, \right. \nonumber \\ 
& \left.\kern - \nulldelimiterspace \frac{\omega_2}{2} \left( \left( x - p_z \right)^2 + x^2 \right) + \frac{\omega_3}{2}\left( p_z^2 + x^2 \right) = h, \right. \nonumber \\
& \left.\kern - \nulldelimiterspace  \dot{z}(x, p_x, y, p_y, z, p_z) > 0 : x > 0\right\}. 
\label{eqn:nhim_x_coord_3dof}
\end{align}
which is a point on the line $p_x = 0$, and marked by a red cross in Fig.~\ref{fig:NSQH_3dof_xpx}. This is also identified by the location of minima in the Lagrangian descriptor values as shown by the one-dimensional slice in Fig.~\ref{fig:NSQH_3dof_xpx} and agrees with $x = \sqrt{2h/\left(\omega_2 + 2\omega_3 \right)}$ to within the grid resolution.

Next, the intersection of the unstable~\eqref{eqn:umani_nsqh_3dof} and stable~\eqref{eqn:smani_nsqh_3dof} manifolds with the isoenergetic two-dimensional surface~\eqref{eqn:sos_symp_3dof_xpx} becomes

\begin{align}
\mathcal{W}^{\rm u}(\mathcal{M}(h)) \cap U_{xp_x}^+ = \left\{  (x, p_x, y, p_y, z, p_z) \; | \; \right. & \left. \kern - \nulldelimiterspace y = 0, z = 0, p_y = 0, x = p_y + p_z, \right. \nonumber \\ 
& \left.\kern - \nulldelimiterspace \frac{\omega_2}{2} \left( x + p_x \right)^2 + \frac{\omega_3}{2}\left( ( x + p_x )^2 + x^2 \right) = h, \right. \nonumber \\ 
& \left.\kern - \nulldelimiterspace \dot{z}(x, p_x, y, p_y, z, p_z) > 0 : \right. \nonumber \\ 
& \left.\kern - \nulldelimiterspace (\omega_2 + 2\omega_3)x + (\lambda + \omega_2 + \omega_3)p_x > 0 \right\}, \\
\mathcal{W}^{\rm s}(\mathcal{M}(h)) \cap U_{xp_x}^+ = \left\{  (x, p_x, y, p_y, z, p_z) \; | \; \right. & \left. \kern - \nulldelimiterspace y = 0, z = 0, p_y = 0, x = 2p_x + p_y + p_z, \right. \nonumber \\ 
& \left.\kern - \nulldelimiterspace \frac{\omega_2}{2} \left( x - p_x \right)^2 + \frac{\omega_3}{2}\left( ( x - p_x )^2 + (x - 2p_x)^2 \right) = h, \right. \nonumber \\ 
& \left.\kern - \nulldelimiterspace \dot{z}(x, p_x, y, p_y, z, p_z) > 0 : \right. \nonumber \\ 
& \left.\kern - \nulldelimiterspace (\omega_2 + 2\omega_3)x + ( -\lambda - \omega_2 -3\omega_3)p_x > 0 \right\}
\end{align}

which are one-dimensional and are shown as a dashed red line (unstable) and a dashed white line (stable) in Fig.~\ref{fig:NSQH_3dof_xpx}. These lines are also identified by the minima in the Lagrangian descriptor values as shown by one-dimensional slice at constant $p_x$.

\begin{figure}[!ht]
	\centering
	\subfigure[]{\includegraphics[width=0.32\textwidth]{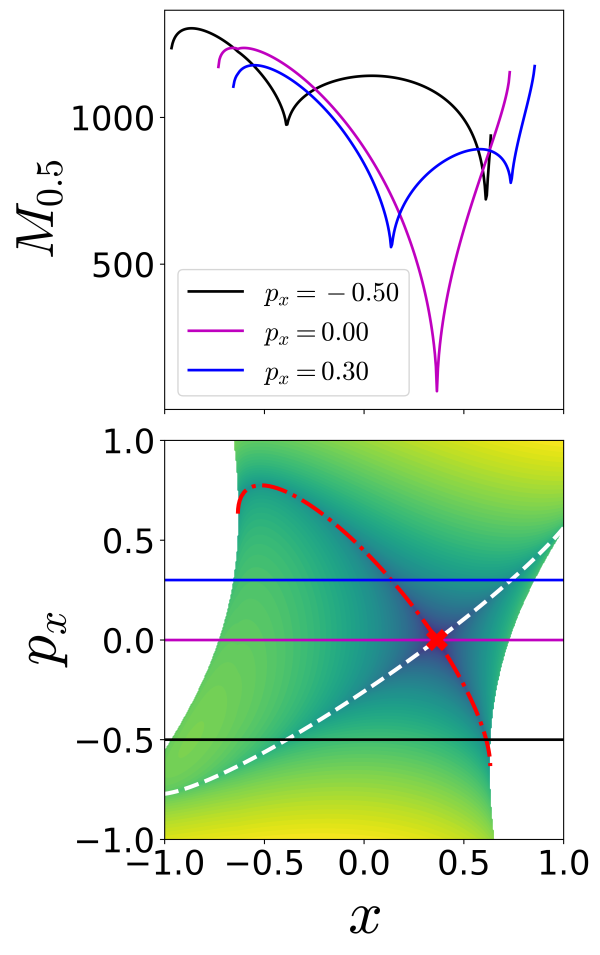}\label{fig:NSQH_3dof_xpx}}
	\subfigure[]{\includegraphics[width=0.32\textwidth]{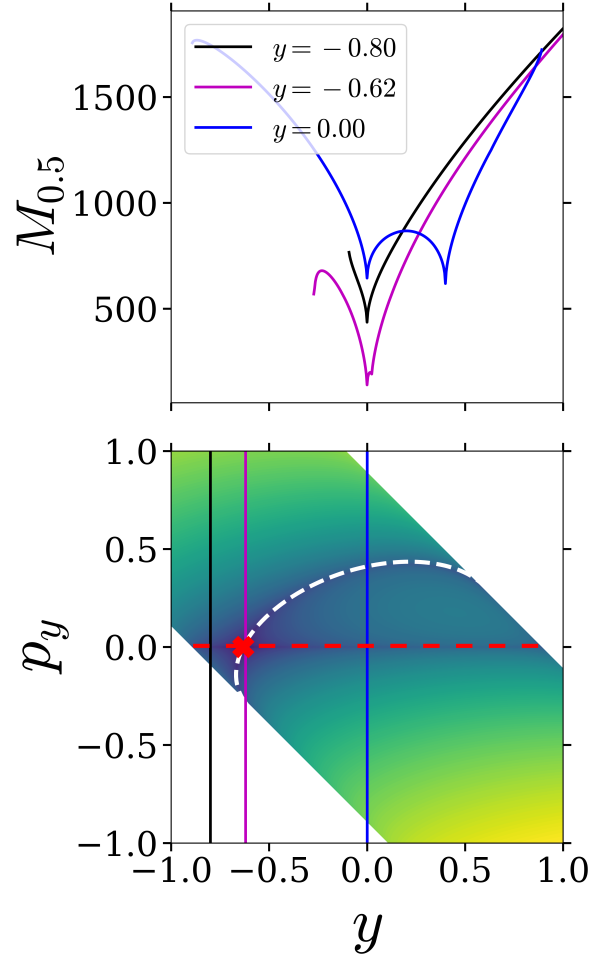}\label{fig:NSQH_3dof_ypy}}
	\subfigure[]{\includegraphics[width=0.32\textwidth]{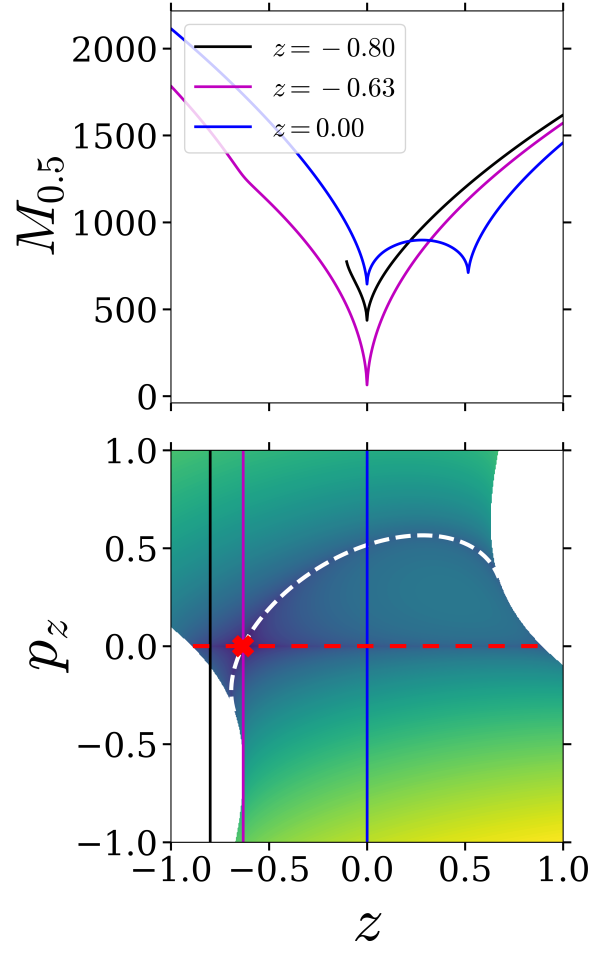}\label{fig:NSQH_3dof_zpz}}
	\caption{Lagrangian descriptor plot of the non-separable quadratic Hamiltonian vector field~\eqref{eqn:eom_symp_3dof} on the \protect\subref{fig:NSQH_3dof_xpx} $U_{xp_x}^+$, \protect\subref{fig:NSQH_3dof_ypy} $U_{yp_y}^+$, \protect\subref{fig:NSQH_3dof_zpz} $U_{zp_z}^+$. The parameters used are $\lambda = \omega_2 = \omega_3 = 1.0$, integration time of $\tau = 10$, and total energy $h = 0.2$.}
	\label{fig:linear_symp_trans_3dof_M400x400_1}
\end{figure}

\newpoint~Isoenergetic two-dimensional surface parametrized by $(y,p_y)$~\textemdash~On a constant energy surface, $\mathcal{H}(x, p_x, y, p_y, z, p_z) = h$, we compute Lagrangian descriptor on a two-dimensional surface by defining

\begin{align}
U_{yp_y}^+ = \left\{ (x, p_x, y, p_y, z, p_z) \; \vert \; \right. & \left. \kern - \nulldelimiterspace x = 0, \; z = 0, \; p_z = 0, \right. \nonumber \\ 
& \left. \kern - \nulldelimiterspace p_x(x,y,p_y,z,p_z;h) > 0 : \dot{x}(x, p_x, y, p_y, z, p_z) > 0 \right\},
\label{eqn:sos_symp_3dof_ypy}
\end{align} 
%


The intersection of the NHIM~\eqref{eqn:nhim_symp_3dof} with the two-dimensional surface is given by  

\begin{align}
\mathcal{M}(h) \cap U_{yp_y}^+ = \left\{ (x, p_x, y, p_y, z, p_z) \; \vert \; \right. & \left. \kern - \nulldelimiterspace x = 0, z = 0, p_z = 0, p_x = 0, p_y + p_z = x, \right. \nonumber \\ 
& \left.\kern - \nulldelimiterspace \frac{\omega_2}{2} y^2 = h, \dot{x}(x, p_x, y, p_y, z, p_z) > 0 \, : \, y < 0 \right\}. 
\label{eqn:nhim_ypy_3dof}
\end{align}

which is a point on the isoenergetic two-dimensional surface~\eqref{eqn:sos_symp_3dof_ypy}, and shown as red cross in Fig.~\ref{fig:NSQH_3dof_ypy}. This point is also identified by the minima in the Lagrangian descriptor values as shown by the one-dimensional slices for constant $y$.

Next, the intersection of the unstable~\eqref{eqn:umani_nsqh_3dof} and stable~\eqref{eqn:smani_nsqh_3dof} manifolds with the isoenergetic two-dimensional surface~\eqref{eqn:sos_symp_3dof_ypy} becomes

\begin{align}
\mathcal{W}^{\rm u}(\mathcal{M}(h)) \cap U_{ypy}^+ = \left\{  (x, p_x, y, p_y, z, p_z) \; | \; \right. & \left. \kern - \nulldelimiterspace x = 0, z = 0, p_z = 0 , x = p_y + p_z, \right. \nonumber \\ 
& \left. \kern - \nulldelimiterspace \frac{\omega_2}{2}\left( -y + p_x \right)^2 + \frac{\omega_3}{2} p_x^2 = h, \right. \nonumber \\
& \left.\kern - \nulldelimiterspace   \dot{x}(x, p_x, y, p_y, z, p_z) > 0 \, : \, -\omega_2y + \lambda p_y > 0 \right\} \\ 
\mathcal{W}^{\rm s}(\mathcal{M}(h)) \cap U_{ypy}^+ = \left\{  (x, p_x, y, p_y, z, p_z) \; | \; \right. & \left. \kern - \nulldelimiterspace x = 0, z = 0, p_z = 0, x = 2p_x + p_y + p_z, \right. \nonumber \\  
& \left. \kern - \nulldelimiterspace \frac{\omega_2}{2}\left( \left( -y + \frac{p_y}{2} \right)^2 + p_y^2 \right) + \frac{\omega_3}{2} \frac{p_y^2}{4} = h : \right. \nonumber \\
& \left.\kern - \nulldelimiterspace \dot{x}(x, p_x, y, p_y, z, p_z) > 0 \, : \, -2\omega_2 y + (2\lambda + \omega_2 + \omega_3)p_y > 0 \right\}
\end{align}

which are lines on the two-dimensional isoenergetic surface~\eqref{eqn:sos_symp_3dof_ypy} and are also identified by the minima in the Lagrangian descriptor values as shown by the one-dimensional slices at constant $y$. 


\newpoint~Isoenergetic two-dimensional surface parametrized by $(z,p_z)$~\textemdash~On a constant energy surface, $\mathcal{H}(x, p_x, y, p_y, z, p_z) = h$~\eqref{eqn:ham_symp_3dof}, we compute Lagrangian descriptor on a two-dimensional surface by defining

\begin{align}
U_{zp_z}^+ = \left\{ (x, p_x, y, p_y, z, p_z) \; \vert \; \right. & \left. \kern - \nulldelimiterspace x = 0, \; y = 0, \; p_y = 0, \right. \nonumber \\ & \left. \kern - \nulldelimiterspace p_x(x,y,p_y,z,p_z;h) > 0 : \dot{x}(x, p_x, y, p_y, z, p_z) > 0 \right\},
\label{eqn:sos_symp_3dof_zpz}
\end{align} 
%

The intersection of the NHIM~\eqref{eqn:nhim_symp_3dof} with the two-dimensional surface is given by  

\begin{align}
\mathcal{M}(h) \cap U_{zp_z}^+ = \left\{ (x, p_x, y, p_y, z, p_z) \; \vert \; \right. & \left. \kern - \nulldelimiterspace x = 0, y = 0, p_y = 0, p_x = 0, p_y + p_z = x, \right. \nonumber \\ 
& \left.\kern - \nulldelimiterspace \frac{\omega_2}{2} p_z^2 + \frac{\omega_3}{2}\left( p_z^2 + z^2 \right) = h, \right. \nonumber \\
& \left.\kern - \nulldelimiterspace  \dot{x}(x, p_x, y, p_y, z, p_z) > 0 \, : \, z < 0 \right\}. 
\label{eqn:nhim_zpz_3dof}
\end{align}

which is a point and shown as red cross in Fig.~\ref{fig:NSQH_3dof_zpz}. This point is also identified by the minima in the Lagrangian descriptor values as shown by the one-dimensional slices for constant $z$.

Next, the intersection of the four-dimensional unstable~\eqref{eqn:umani_nsqh_3dof} and stable~\eqref{eqn:smani_nsqh_3dof} manifolds with the isoenergetic two-dimensional surface~\eqref{eqn:sos_symp_3dof_zpz} becomes

\begin{align}
\mathcal{W}^{\rm u}(\mathcal{M}(h)) \cap U_{zpz}^+ = \left\{  (x, p_x, y, p_y, z, p_z) \; | \; \right. & \left. \kern - \nulldelimiterspace x = 0, y = 0, p_y = 0, p_z = 0, \right. \nonumber \\ 
& \left. \kern - \nulldelimiterspace \frac{\omega_2}{2}p_x^2 + \frac{\omega_3}{2}\left( (-z + p_x)^2 \right) = h :  \right. \nonumber \\ 
& \left. \kern - \nulldelimiterspace  \dot{x}(x, p_x, y, p_y, z, p_z) > 0 \, : \, \omega_2 p_x + \omega_3( -z + p_x ) > 0 \right\} \\ 
\mathcal{W}^{\rm s}(\mathcal{M}(h)) \cap U_{zpz}^+ = \left\{  (x, p_x, y, p_y, z, p_z) \; | \; \right. & \left. \kern - \nulldelimiterspace x = 0,  y = 0, p_y = 0, p_z + 2p_x = 0, \right. \nonumber \\  
& \left. \kern - \nulldelimiterspace \frac{\omega_2}{2}\left( \frac{p_z}{2} \right)^2 + \frac{\omega_3}{2}\left( \left(-z + \frac{p_z}{2}\right) ^2 + p_z^2\right) = h :  \right. \nonumber \\ 
& \left. \kern - \nulldelimiterspace  \dot{x}(x, p_x, y, p_y, z, p_z) > 0 \, : \, -\omega_3 z + (\lambda + \omega_2/2 + \omega_3/2)p_z > 0 \right\} 
\end{align}

which are one-dimensional and shown as dashed red (unstable) and dashed white (stable) lines in Fig.~\ref{fig:NSQH_3dof_zpz}. These invariant manifolds are also identified by minima in the Lagrangian descriptor values as shown by the one-dimensional slice at constant $z$.

\newpoint~Isoenergetic two-dimensional surface parametrized by $(x,p_z)$~\textemdash~On a constant energy surface, $\mathcal{H}(x, p_x, y, p_y, z, p_z) = h$, we compute the Lagrangian descriptor on a two-dimensional surface by defining

\begin{align}
U_{xp_z}^+ = \left\{ (x, p_x, y, p_y, z, p_z) \; \vert \; \right. & \left. \kern - \nulldelimiterspace p_x = 0, y = 0, z = 0, \right. \nonumber \\ & \left. \kern - \nulldelimiterspace p_y = p_y(x,p_x,y,p_y,z;h) : \dot{y}(x, p_x, y, p_y, z, p_z) > 0 \right\},
\label{eqn:sos_symp_3dof_xpz}
\end{align} 

The intersection of the NHIM~\eqref{eqn:nhim_symp_3dof} with the two-dimensional surface is given by

\begin{align}
\mathcal{M}(h) \cap U_{xp_z}^+ = \left\{ (x, p_x, y, p_y, z, p_z) \; \vert \; \right. & \left. \kern - \nulldelimiterspace p_x = 0, y = 0, z = 0, p_y + p_z = x, \right. \nonumber \\ 
& \left.\kern - \nulldelimiterspace \frac{\omega_2}{2}\left( (x-p_z)^2 + x^2 \right) + \frac{\omega_3}{2}\left( p_z^2 + x^2 \right) = h, \right. \nonumber \\  
& \left.\kern - \nulldelimiterspace  \dot{y}(x, p_x, y, p_y, z, p_z) > 0 : \omega_2 (2x - p_z) + \omega_3x > 0 \right\}. 
\label{eqn:nhim_xpz_3dof}
\end{align}

which is also identified by the minima in the Lagrangian descriptor values as shown in the one-dimensional slices at constant $p_z$.

Next, the intersection of the unstable~\eqref{eqn:umani_nsqh_3dof} and stable~\eqref{eqn:smani_nsqh_3dof} manifolds with the isoenergetic two-dimensional surface~\eqref{eqn:sos_symp_3dof_xpz} becomes

\begin{align}
\mathcal{W}^{\rm u}(\mathcal{M}(h)) \cap U_{xp_z}^+ = \left\{  (x, p_x, y, p_y, z, p_z) \; | \; \right. & \left. \kern - \nulldelimiterspace p_x = 0, y = 0, z = 0, x = p_y + p_z, \right. \nonumber \\ 
& \left. \kern - \nulldelimiterspace \frac{\omega_2}{2} \left( x^2 + (x - p_z)^2 \right)  + \frac{\omega_3}{2} \left( x^2 + p_z^2 \right)= h, \right. \nonumber \\
& \left.\kern - \nulldelimiterspace  \dot{y}(x, p_x, y, p_y, z, p_z) > 0 : \omega_2 (2x - p_z) + \omega_3x > 0 \right\}, \\
\mathcal{W}^{\rm s}(\mathcal{M}(h)) \cap U_{xp_z}^+ = \left\{  (x, p_x, y, p_y, z, p_z) \; | \; \right. & \left. \kern - \nulldelimiterspace p_x = 0, y = 0, z = 0, x = 2p_x + p_y + p_z, \right. \nonumber \\  
& \left. \kern - \nulldelimiterspace \frac{\omega_2}{2} \left( x^2 + (x - p_z)^2 \right)  + \frac{\omega_3}{2} \left( x^2 + p_z^2 \right)= h, \right. \nonumber \\
& \left.\kern - \nulldelimiterspace  \dot{y}(x, p_x, y, p_y, z, p_z) > 0 : \omega_2 (2x - p_z) + \omega_3x > 0 \right\}
\end{align}

which project as the same one-dimensional curve as the NHIM and identified by the minima in the Lagrangian descriptor values as shown by the one-dimensional slices at constant $p_z$ in Fig.~\ref{fig:NSQH_3dof_xpz}.

\begin{figure}[!ht]
	\centering
	\subfigure[]{\includegraphics[width=0.32\textwidth]{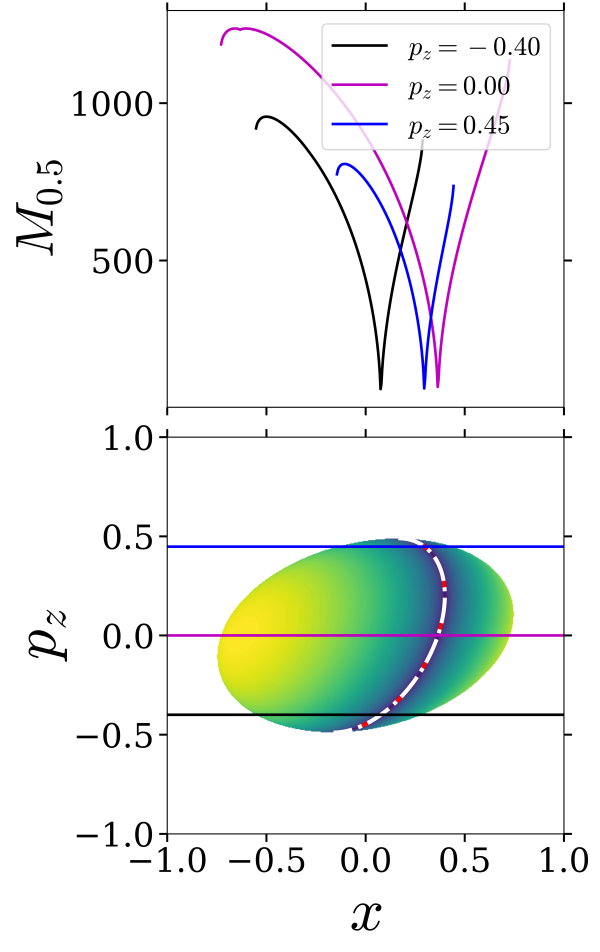}\label{fig:NSQH_3dof_xpz}}
	\subfigure[]{\includegraphics[width=0.32\textwidth]{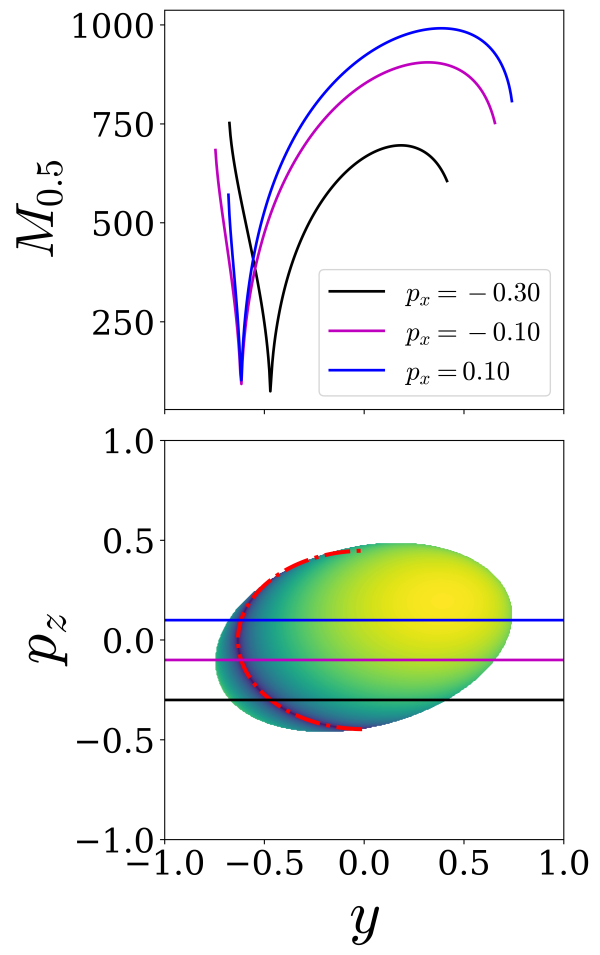}\label{fig:NSQH_3dof_ypz}}
	\subfigure[]{\includegraphics[width=0.32\textwidth]{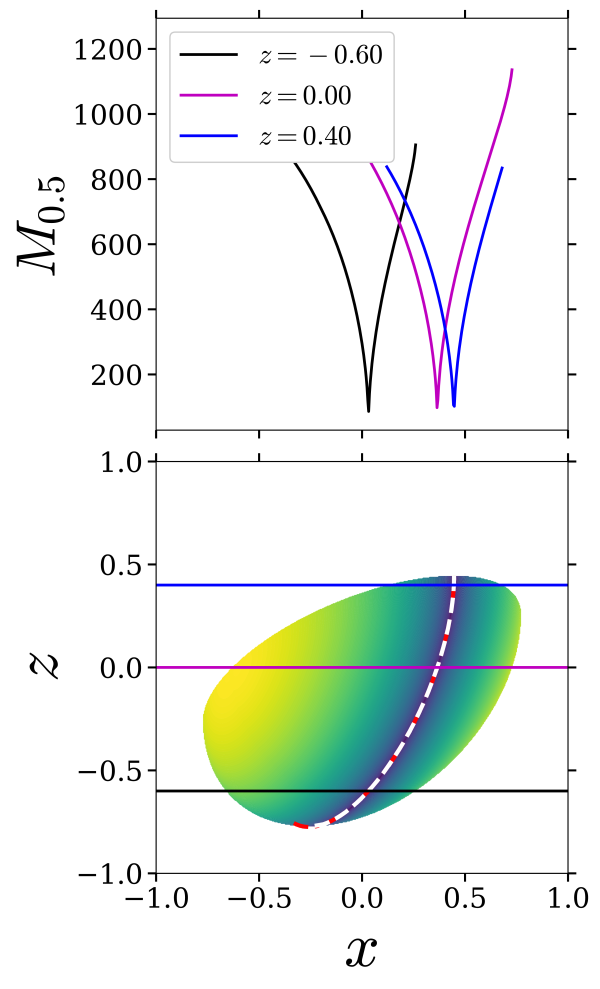}\label{fig:NSQH_3dof_xz}}
	\caption{LD plot of the transformed Hamiltonian vector field~\eqref{eqn:eom_symp_3dof} on the  \protect\subref{fig:NSQH_3dof_ypz} $U_{xp_z}^+$,  \protect\subref{fig:NSQH_3dof_ypz} $U_{yp_z}^+$, \protect\subref{fig:NSQH_3dof_xz} $U_{xz}^+$. The parameters used are $\lambda = \omega_2 = \omega_3 = 1.0$, integration time of $\tau = 10$, and total energy $h = 0.2$.}
	\label{fig:linear_symp_trans_3dof_M400x400_2}
\end{figure}

\newpoint~Isoenergetic two-dimensional surface parametrized by $(y,p_z)$~\textemdash~Next, on a constant energy surface, $\mathcal{H}(x, p_x, y, p_y, z, p_z) = h$, we compute Lagrangian descriptor on a two-dimensional surface by defining 

\begin{align}
U_{yp_z}^+ = \left\{ (x, p_x, y, p_y, z, p_z) \; \vert \; \right. & \left. \kern - \nulldelimiterspace x = 0, p_x = 0, z = 0, \right. \nonumber \\
& \left.\kern - \nulldelimiterspace  p_y(x,p_x,y,z,p_z;h) > 0 \, : \, \dot{x}(x,p_x,y,p_y,z,p_z) > 0 \right\}, 
\label{eqn:sos_symp_3dof_ypz}
\end{align} 
%

Thus, the intersection of the NHIM~\eqref{eqn:nhim_symp_3dof} with the two-dimensional surface is given by

\begin{align}
\mathcal{M}(h) \cap U_{yp_z}^+ = \left\{ (x, p_x, y, p_y, z, p_z) \; \vert \; \right. & \left. \kern - \nulldelimiterspace x = 0, p_x = 0, z = 0, p_y + p_z = x, \right. \nonumber \\ 
& \left.\kern - \nulldelimiterspace  \frac{(\omega_2 + \omega_3)}{2}p_z^2 + \frac{\omega_2}{2}y^2 = h, \right. \nonumber \\  
& \left.\kern - \nulldelimiterspace  \dot{x}(x,p_x,y,p_y,z,p_z) > 0 \, : \, y < 0 \right\}. 
\label{eqn:nhim_ypz_3dof}
\end{align}


which represents a portion of the ellipse and shown as the dashdot red curve in Fig.~\ref{fig:NSQH_3dof_ypz}. As shown by the one-dimensional slices for constant $p_z$, Eqn.~\eqref{eqn:nhim_ypz_3dof} is identified by the minima in the Lagrangian descriptor values. 

Next, the intersection of the unstable~\eqref{eqn:umani_nsqh_3dof} and stable~\eqref{eqn:smani_nsqh_3dof} manifolds with the isoenergetic two-dimensional surface~\eqref{eqn:sos_symp_3dof_ypz} becomes

\begin{align}
\mathcal{W}^{\rm u}(\mathcal{M}(h)) \cap U_{ypz}^+ = \left\{  (x, p_x, y, p_y, z, p_z) \; | \; \right. & \left. \kern - \nulldelimiterspace x = 0, p_x = 0, z = 0, x = p_y + p_z, \right. \nonumber \\ 
& \left. \kern - \nulldelimiterspace \frac{\omega_2}{2} y^2  + \frac{(\omega_2 + \omega_3)}{2} p_z^2 = h, \right. \nonumber \\
& \left.\kern - \nulldelimiterspace  \dot{x}(x,p_x,y,p_y,z,p_z) > 0 : y < 0 \right\}, \\
\mathcal{W}^{\rm s}(\mathcal{M}(h)) \cap U_{ypz}^+ = \left\{  (x, p_x, y, p_y, z, p_z) \; | \; \right. & \left. \kern - \nulldelimiterspace x = 0, p_x = 0, z = 0, x = 2p_x + p_y + p_z, \right. \nonumber \\  
& \left. \kern - \nulldelimiterspace \frac{\omega_2}{2} y^2  + \frac{(\omega_2 + \omega_3)}{2} p_z^2 = h, \right. \nonumber \\
& \left.\kern - \nulldelimiterspace  \dot{x}(x,p_x,y,p_y,z,p_z) > 0 : y < 0 \right\}
\end{align}

which are shown as red dashed line (unstable) and white dashed line (stable) in Fig.~\ref{fig:NSQH_3dof_ypz}. We note here that the intersection of the manifolds with the isoenergetic two-dimensional surface~\eqref{eqn:sos_symp_3dof_ypz} overlays on the NHIM's intersection~\eqref{eqn:nhim_ypz_3dof}. These curves are also identified by the minima in the Lagrangian descriptor values as shown by the one-dimensional slice for constant $p_z$.

\newpoint~Isoenergetic two-dimensional surface parametrized by $(x,z)$~\textemdash~Next, on a constant energy surface, $\mathcal{H}(x, p_x, y, p_y, z, p_z) = h$~\eqref{eqn:ham_symp_3dof}, we compute the Lagrangian descriptor on a two-dimensional surface by defining

\begin{align}
U_{xz}^+ = \left\{ (x, p_x, y, p_y, z, p_z) \; \vert \; \right. & \left. \kern - \nulldelimiterspace y = 0, p_x = 0, p_z = 0, \right. \nonumber \\ & \left. \kern - \nulldelimiterspace p_y = p_y(x,p_x,y,z,p_z;h) : \dot{p_z}(x, p_x, y, p_y, z, p_z) > 0 \right\},
\label{eqn:sos_symp_3dof_xz}
\end{align} 

The intersection of the NHIM~\eqref{eqn:nhim_symp_3dof} with the two-dimensional surface is given by  

\begin{align}
{\rm NHIM} \cap U_{xz}^+ = \left\{ (x, p_x, y, p_y, z, p_z) \; \vert \; \right. & \left. \kern - \nulldelimiterspace y = 0, p_x = 0, p_z = 0, p_y + p_z = x, \right. \nonumber \\ 
& \left.\kern - \nulldelimiterspace \omega_2 x^2 + \frac{\omega_3}{2}\left( x - z \right)^2 = h, \right. \nonumber \\
& \left.\kern - \nulldelimiterspace  \dot{p_z}(x, p_x, y, p_y, z, p_z) > 0 : \omega_3(-z + x) > 0 \right\}. 
\label{eqn:nhim_xz_3dof}
\end{align}

which represents a portion of an ellipse and is identified by the minima in the Lagrangian descriptor values as shown by the one-dimensional slices for constant $z$ in Fig.~\ref{fig:NSQH_3dof_xz}. 

Next, the intersection of the unstable~\eqref{eqn:umani_nsqh_3dof} and stable~\eqref{eqn:smani_nsqh_3dof} manifolds with the isoenergetic two-dimensional surface~\eqref{eqn:sos_symp_3dof_ypz} becomes

\begin{align}
W^{\rm u}({\rm NHIM}) \cap U_{xz}^+ = \left\{  (x, p_x, y, p_y, z, p_z) \; | \; \right. & \left. \kern - \nulldelimiterspace y = 0, p_x = 0, p_z = 0, x = p_y + p_z, \right. \nonumber \\ 
& \left. \kern - \nulldelimiterspace \omega_2x^2  + \frac{\omega_3}{2} \left( -z + x \right)^2 = h, \right. \nonumber \\
& \left.\kern - \nulldelimiterspace  \dot{p_z}(x, p_x, y, p_y, z, p_z) > 0 : \omega_3(-z + x) > 0 \right\}, \\
W^{\rm s}({\rm NHIM}) \cap U_{xz}^+ = \left\{  (x, p_x, y, p_y, z, p_z) \; | \; \right. & \left. \kern - \nulldelimiterspace y = 0, p_x = 0, p_z = 0, x = 2p_x + p_y + p_z, \right. \nonumber \\  
& \left. \kern - \nulldelimiterspace \omega_2x^2  + \frac{\omega_3}{2} \left( -z + x \right)^2 = h, \right. \nonumber \\
& \left.\kern - \nulldelimiterspace  \dot{p_z}(x, p_x, y, p_y, z, p_z) > 0 : \omega_3(-z + x) > 0 \right\}
\end{align}

which are also be identified by the minima in the Lagrangian descriptor values as shown by the one-dimensional slice for constant $z$ in Fig.~\ref{fig:NSQH_3dof_xz}.

\section{Summary and Outlook}
\label{sec:summ}
In this article, we assessed Lagrangian descriptor (LD) based detection of high dimensional phase space structures, that is 
normally hyperbolic invariant manifolds, their unstable and stable manifolds, in the two and three degrees-of-freedom 
quadratic Hamiltonian systems with index-1 saddle. This is done using a sytematic comparison of numerical results with 
analytical expressions for a comprehensive set of coordinates; see Appendix~\ref{sect:LD_add_slices} for more pairs of 
coordinates that support the results herein.. It is to be noted that in $N$ DoF system, there are $2N(2N-1)(N-1)$ 
two-dimensional sections, so the detection strategy should include inspecting a few pairs of coordinates. Based on our 
investigation, a judicious choice of coordinate for the two-dimensional surface is to use one of the reaction coordinates in 
phase space, that is $q_1$ or $p_1$, and the other one from the remaining $2N-1$ coordinates.
 
The closed-form analytical expressions of the NHIM and their manifolds are precisely identified by the minima and singular features in the Lagrangian descriptor values. This provides further numerical evidence for the use of LDs in detecting high dimensional phase space structures with a simple computation that can be implemented along with trajectory integration. So at least for the form of Hamiltonian considered here, one can rely with certainity on using Lagrangian descriptor to detect the NHIM and its stable and unstable manifolds. Furthermore, this detection approach has the potential to be combined with machine learning type methods~\cite{feldmaier_invariant_2019}. Future work on this approach will include specific nonlinear systems that are inspired by applications in celestial mechanics, ship dynamics, structural mechanics, and chemical reaction dynamics.

\section*{Acknowledgements} 
We acknowledge the support of EPSRC Grant No. ~EP/P021123/1 and  ONR Grant No.~N00014-01-1-0769. We would like to thank Vladim{\'i}r Kraj{\v n}{\'a}k for useful discussions and feedback.


\pagebreak
\renewcommand{\thesubsection}{\Alph{subsection}}
\setcounter{section}{0}
\section*{Appendix}
\addcontentsline{toc}{section}{Appendix}


\subsection{Eigenvalues of linear symplectic transformed system}
\label{sect:linear_symp_transf_jac}

The Jacobian of the vector field~\eqref{eqn:eom_symp_2dof} near the equilibrium point is given by
\begin{equation}
\mathbb{J} =
\begin{pmatrix}
-\lambda & -\omega_2 & \omega_2 & (\lambda + \omega_2) \\
-\lambda & -\omega_2 & (\lambda + \omega_2) & (\lambda + 2\omega_2) \\
-\lambda & 0 & \lambda & \lambda \\
0 & -\omega_2 & \omega_2 & \omega_2
\end{pmatrix}	
\end{equation}
which gives the characteristic polynomial
\begin{equation}
\det(\mathbb{J} - \beta \mathbb{I}) := \beta^4 - \lambda^2 \beta^2 + \omega_2^2 \beta^2 - \lambda^2 \omega_2^2    
\end{equation}
which has solutions $\lambda, -\lambda, i \omega_2, -i \omega_2$.

The Jacobian of the vector field~\eqref{eqn:eom_symp_3dof} near the equilibrium point is given by
\begin{equation}
\mathbb{J} =
\begin{pmatrix}
-\lambda & -\omega_2 & -\omega_3 & (\omega_2 + \omega_3) & (\lambda + \omega_2 + \omega_3) & 
(\lambda + \omega_2 + \omega_3) \\
-\lambda & -\omega_2 & -\omega_3 & (\lambda + \omega_2 + \omega_3) & (\lambda + 2\omega_2 + 
\omega_3) & (\lambda + \omega_2 + \omega_3)\\
-\lambda & -\omega_2 & -\omega_3 & (\lambda + \omega_2 + \omega_3) & (\lambda + \omega_2 + 
\omega_3) & (\lambda + \omega_2 + 2\omega_3) \\
-\lambda & 0 & 0 & \lambda & \lambda & \lambda \\
0 & -\omega_2 & 0 & \omega_2 & \omega_2 & \omega_2 \\
0 & 0 & -\omega_3 & \omega_3 & \omega_3 & \omega_3\\
\end{pmatrix}	
\end{equation}
which gives the characteristic polynomial
\begin{equation}
\det(\mathbb{J} - \beta \mathbb{I}) := \beta^6 - \beta^4 \lambda^2 + \beta^4 \omega_2^2 + \beta^4 
\omega_3^2 - \beta^2 \lambda^2 \omega_2^2 - \beta^2 \lambda^2 \omega_3^2 + \beta^2 \omega_2^2 
\omega_3^2 - \lambda^2 \omega_2^2 \omega_3^2 
\end{equation}
which has solutions $\lambda, -\lambda, i \omega_2, -i \omega_2, i \omega_3, - i \omega_3$.

\subsection{Examples of symplectic transformations}
\label{sect:examples_C}

Along with the transformations $C: (x, y, p_x, p_y) \rightarrow (q_1, q_2, p_1, p_2)$ used for demonstrating the procedure 
of detecting NHIM, its unstable, and stable manifolds, we have verified the following transformations for symplectic 
condition. 

\textbf{Two degrees of freedom}

\begin{equation}
C = \begin{bmatrix}
0 & 0 & 1 & 0 \\
0 & 0 & 0 & 1 \\
-1 & 0 & 1 & 1 \\
0 & -1 & 1 & 1 \\
\end{bmatrix}, \;
C = \begin{bmatrix}
0 & 0 & 1 & 0 \\
0 & 0 & 0 & 1 \\
-1 & 0 & 1 & 0 \\
0 & -1 & 0 & 1 \\
\end{bmatrix}, \;
C = \begin{bmatrix}
0 & 0 & 1 & 0 \\
0 & 0 & 0 & 1 \\
-1 & 0 & 0 & 1 \\
0 & -1 & 1 & 0 \\
\end{bmatrix}
\end{equation}

\textbf{Three degrees of freedom}

The following  $C: (x, y, z, p_x, p_y, p_z) \rightarrow (q_1, q_2, q_3, p_1, p_2, p_3)$ have been verfied to be symplectic.

\begin{equation}
C = \begin{bmatrix}
0 & 0 & 0 & 1 & 0 & 0 \\
0 & 0 & 0 & 0 & 1 & 0 \\
0 & 0 & 0 & 0 & 0 & 1 \\
-1 & 0 & 0 & 1 & 1 & 1 \\
0 & -1 & 0 & 1 & 1 & 1 \\
0 & 0 & -1 & 1 & 1 & 1
\end{bmatrix}, \;
C = \begin{bmatrix}
0 & 0 & 0 & 1 & 0 & 0 \\
0 & 0 & 0 & 0 & 1 & 0 \\
0 & 0 & 0 & 0 & 0 & 1 \\
-1 & 0 & 0 & 1 & 0 & 0 \\
0 & -1 & 0 & 0 & 1 & 0 \\
0 & 0 & -1 & 0 & 0 & 1
\end{bmatrix}, \;
C = \begin{bmatrix}
0 & 0 & 0 & 1 & 0 & 0 \\
0 & 0 & 0 & 0 & 1 & 0 \\
0 & 0 & 0 & 0 & 0 & 1 \\
-1 & 0 & 0 & 0 & 1 & 1 \\
0 & -1 & 0 & 1 & 0 & 1 \\
0 & 0 & -1 & 1 & 1 & 0
\end{bmatrix}
\end{equation}

\subsection{Lagrangian Descriptor on isoenergetic two-dimensional surfaces}
\label{sect:LD_add_slices}

\subsubsection{Decoupled quadratic Hamiltonian: 2 DoF}
\begin{figure}[!th]
	\centering
	\subfigure[]{\includegraphics[width=0.32\textwidth]{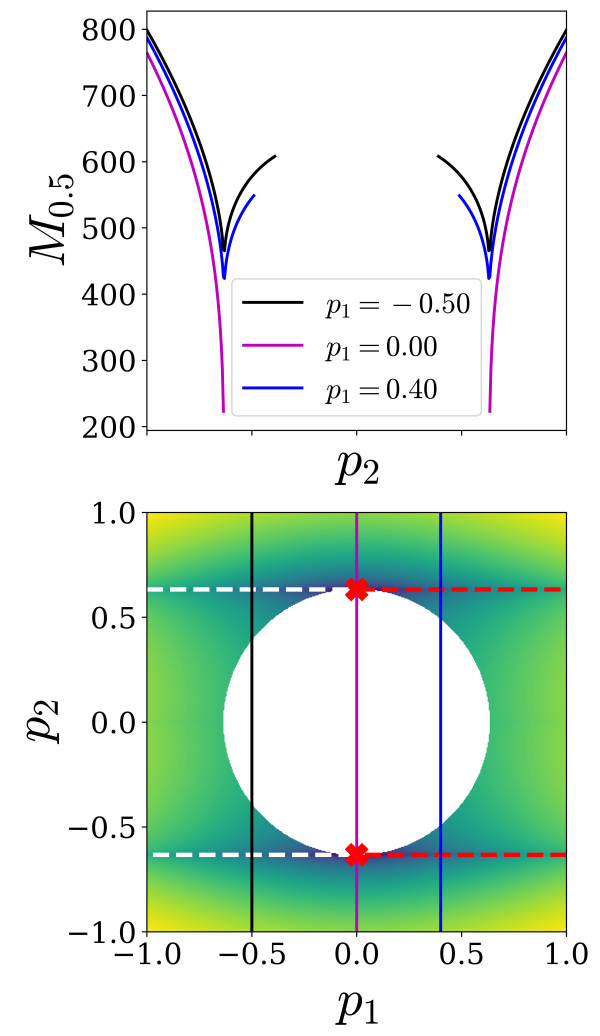}\label{fig:LD_SQH_2dof_p1p2}}
	\subfigure[]{\includegraphics[width=0.32\textwidth]{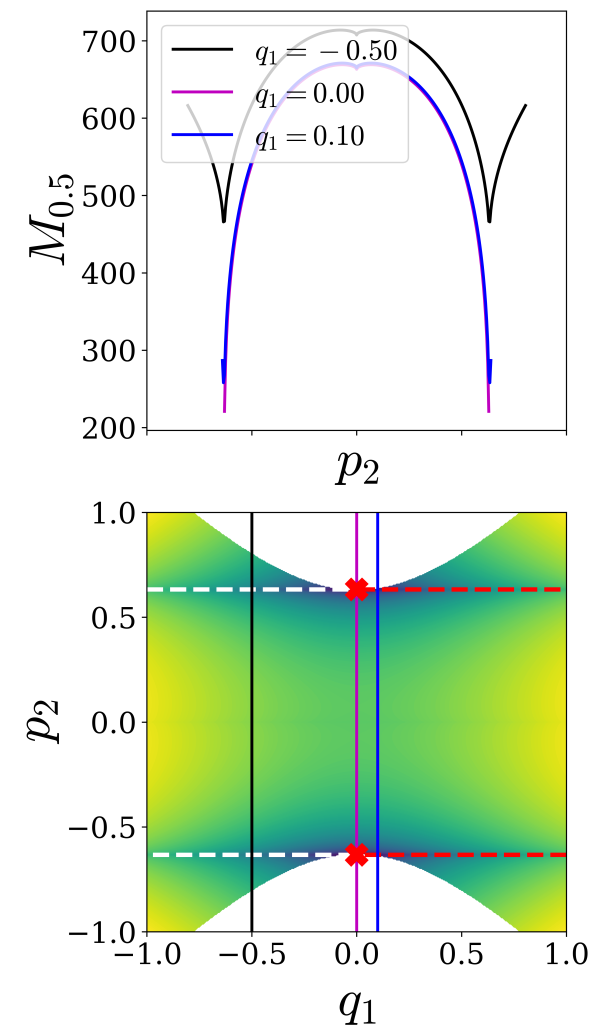}\label{fig:LD_SQH_2dof_q1p2}}
	\subfigure[]{\includegraphics[width=0.32\textwidth]{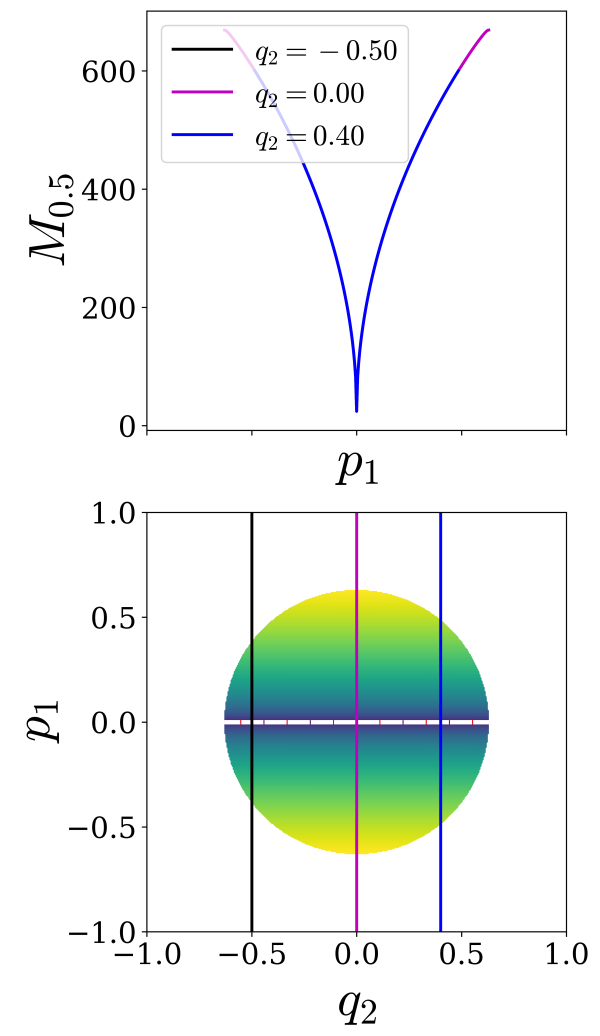}\label{fig:LD_SQH_2dof_q2p1}}
	\caption{Lagrangian descriptor slice for the two degrees of freedom separable quadratic Hamiltonian on the isoenergetic 
	two-dimensional surfaces \protect\subref{fig:LD_SQH_2dof_p1p2} $U_{p_1p_2}^+$, \protect\subref{fig:LD_SQH_2dof_q1p2} 
	$U_{q_1p_2}^+$, \protect\subref{fig:LD_SQH_2dof_q2p1} $U_{q_2p_1}^+$. Parameters used are $\lambda = 1.0, \omega_2 = 
	1.0$, $h = 0.2$, and integration time $\tau = 10$ is fixed.}
	\label{fig:LD_SQH_2dof_2}
\end{figure}



\newpoint~Isoenergetic two-dimensional surface parametrized by $(p_1, p_2)$~\textemdash~Similarly, on the constant energy 
surface, $H(q_1, p_1, q_2, p_2) = h$, we take the Lagrangian descriptor slice in $(q_2, p_2)$ coordinates by defining the 
two dimensional surface
\begin{equation}
U_{p_1p_2}^+ = \left\{ (q_1, p_1, q_2, p_2) \; | \; q_2 = 0, q_1(q_2, p_1, p_2; h) \geqslant 0 \right\}
\end{equation}
where
\begin{align}
q_1(p_1, q_2 = 0, p_2; h) = & \sqrt{\frac{2}{\lambda}\left( \left( \frac{\lambda}{2}p_1^2 + \frac{\omega_2}{2}p_2^2 \right) 
- h\right)} 
\end{align}

The intersection of the NHIM~\eqref{eqn:sep_quad_ham2dof_nhim} with the two-dimensional surface $U_{p_1p_2}^+$ is

\begin{align}
\mathcal{M}(h) \cap U_{p_1p_2}^+ = \left\{ (q_1, p_1, q_2, p_2) \; \vert \; \right. & \left. \kern - \nulldelimiterspace p_1 
= 0, q_1 = 0, q_2 = 0 : \right. \nonumber \\ 
& \left.\kern - \nulldelimiterspace   \frac{\omega_2}{2}\left( p_2^2 + q_2^2 \right) = h \right\}.
\end{align}

which represents two points on the line $p_1 = 0$. Their location is $p_2 = \pm \sqrt{2h/\omega_2}$ which is $\approx 0.632$ 
for $h = 0.2, \omega_2 = 1.0$ and is marked as red crosses in Fig.~\ref{fig:LD_SQH_2dof_p1p2}. These points are also 
locations of minima in LD values as shown by the one-dimensional slices at constant $p_1$.

Next, intersection of the two-dimensional $U_{p_1p_2}^+$ with the stable~\eqref{eqn:quad_ham2dof_smani} and unstable 
manifolds~\eqref{eqn:quad_ham2dof_umani} is given by

\begin{align}
\mathcal{W}^u(\mathcal{M}(h)) \cap U_{p_1p_2}^+ = \left\{ (q_1, p_1, q_2, p_2) \; \vert \; \right. & \left. \kern - 
\nulldelimiterspace p_1 = q_1, q_1 \geqslant 0, q_2 = 0 : \right. \nonumber \\  
& \left.\kern - \nulldelimiterspace \frac{\omega_2}{2}\left( p_2^2 + q_2^2 \right) = h \right\}, \\
\mathcal{W}^s(\mathcal{M}(h)) \cap U_{p_1p_2}^+ = \left\{ (q_1, p_1, q_2, p_2) \; \vert \; \right. & \left. \kern - 
\nulldelimiterspace p_1 = -q_1, q_1 \geqslant 0, q_2 = 0 : \right. \nonumber \\  
& \left.\kern - \nulldelimiterspace \frac{\omega_2}{2}\left( p_2^2 + q_2^2 \right) = h \right\},
\end{align}

where each manifold represent lines parallel to $p_1$ axis for $p_2 = \pm \sqrt{\frac{2h}{\omega_2}}$. The unstable manifold 
lies on the $p_1 > 0$ plane and stable manifold lies on the $p_1 < 0$ plane, and are marked by dashed red (unstable) and 
dashed white (stable) lines in the Fig.~\ref{fig:LD_SQH_2dof_p1p2}. 

%
%


\newpoint~Isoenergetic two-dimensional surface parametrized by $(q_1, p_2)$~\textemdash~On the constant energy surface, 
$H(q_1, p_1, q_2, p_2) = h$, we take the LD slice by defining the two dimensional surface
\begin{equation}
U_{q_1p_2}^+ = \left\{ (q_1, p_1, q_2, p_2) \; | \; q_2 = 0, p_1(q_1, q_2, p_2; h) \geqslant 0 \right\}
\end{equation}
where
\begin{align}
p_1(q_1, q_2 = 0, p_2; h) = & \sqrt{\frac{2}{\lambda}\left( h - \left( \frac{\omega_2}{2}p_2^2 
-\frac{\lambda}{2}q_1^2\right) \right)}
\end{align}

The intersection of the two-dimensional surface $U_{q_1p_2}^+$ with the NHIM~\eqref{eqn:sep_quad_ham2dof_nhim} is given by

\begin{align}
\mathcal{M}(h) \cap U_{q_1p_2}^+ = \left\{ (q_1, p_1, q_2, p_2) \; \vert \; \right. & \left. \kern - \nulldelimiterspace q_1 
= 0, p_1 = 0, q_2 = 0, \right. \nonumber \\ 
& \left.\kern - \nulldelimiterspace   \frac{\omega_2}{2}\left( p_2^2 + q_2^2 \right) = h \right\}.
\end{align}

which gives two points on the $q_1 = 0$ line with $p_2$ satisfying $p_2 = \pm \sqrt{2h/\omega_2}$. This is shown by the 
dashdot line in Fig.~\ref{fig:LD_SQH_2dof_q1p2}.

Next, the intersection of the two-dimensional $U_{q_1p_2}^+$ with the stable~\eqref{eqn:quad_ham2dof_smani} and unstable 
manifolds~\eqref{eqn:quad_ham2dof_umani} is given by

\begin{align}
\mathcal{W}^u(\mathcal{M}(h)) \cap U_{q_1p_2}^+ = \left\{ (q_1, p_1, q_2, p_2) \; \vert \; \right. & \left. \kern - 
\nulldelimiterspace p_1 = q_1, p_1(q_1, q_2, p_2; h) \geqslant 0, q_2 = 0 : \right. \nonumber \\  
& \left.\kern - \nulldelimiterspace \frac{\omega_2}{2}\left( p_2^2 + q_2^2 \right) = h \right\}, \\
\mathcal{W}^s(\mathcal{M}(h)) \cap U_{q_1p_2}^+ = \left\{ (q_1, p_1, q_2, p_2) \; \vert \; \right. & \left. \kern - 
\nulldelimiterspace p_1 = -q_1, p_1(q_1, q_2, p_2; h) \geqslant 0, q_2 = 0 : \right. \nonumber \\  
& \left.\kern - \nulldelimiterspace \frac{\omega_2}{2}\left( p_2^2 + q_2^2 \right) = h \right\},
\end{align}

where each manifold represents lines parallel to the $q_1$ axis. The unstable manifold lies on the $q_1 > 0$ plane and the 
stable manifold lies on the $q_1 < 0$ plane, and are marked by dashed red (unstable) and dashed white (stable) lines in the 
Fig.~\ref{fig:LD_SQH_2dof_q1p2}. These manifolds are identified by points of minima located on the $q_1 = 0$ line.

\newpoint~Isoenergetic two-dimensional surface parametrized by $(q_2, p_1)$~\textemdash~On the constant energy surface, 
$H(q_1, p_1, q_2, p_2) = h$, we take the LD slice by defining the two dimensional surface

\begin{equation}
U_{q_2p_1}^+ = \left\{ (q_1, p_1, q_2, p_2) \; | \; q_1 = 0, p_2(q_1, p_1, q_2; h) \geqslant 0 \right\}
\end{equation}

where

\begin{align}
p_2(q_1 = 0, p_1, q_2; h) = & \sqrt{\frac{2}{\omega_2}\left( h - \left( \frac{\lambda}{2}p_1^2 + \frac{\omega_2}{2}q_2^2 
\right) \right)}
\end{align}

The intersection of the two-dimensional surface $U_{q_2p_1}^+$ with the NHIM~\eqref{eqn:sep_quad_ham2dof_nhim} is given by

\begin{align}
\mathcal{M}(h) \cap U_{q_2p_1}^+ = \left\{ (q_1, p_1, q_2, p_2) \; \vert \; \right. & \left. \kern - \nulldelimiterspace q_1 
= 0, p_1 = 0, p_2(q_1, p_1, q_2; h) \geqslant 0 : \right. \nonumber \\ 
& \left.\kern - \nulldelimiterspace   \frac{\omega_2}{2}\left( p_2^2 + q_2^2 \right) = h \right\}.
\end{align}

which gives points on the $p_1 = 0$ line with $q_2$ satisfying $q_2 = \pm \sqrt{2h/\omega_2 - p_2^2}$. This is shown as a 
dashdot red line in Fig.~\ref{fig:LD_SQH_2dof_q2p1} and is also the line with minima as shown by the one-dimensional slice 
at constant $q_2$.

Next, the intersection of the two-dimensional surface $U_{q_2p_1}^+$ with the stable~\eqref{eqn:quad_ham2dof_smani} and 
unstable manifolds~\eqref{eqn:quad_ham2dof_umani} is given by

\begin{align}
\mathcal{W}^u(\mathcal{M}(h)) \cap U_{q_2p_1}^+ = \left\{ (q_1, p_1, q_2, p_2) \; \vert \; \right. & \left. \kern - 
\nulldelimiterspace p_1 = q_1, p_1(q_1, q_2, p_2; h) \geqslant 0, q_2 = 0 : \right. \nonumber \\  
& \left.\kern - \nulldelimiterspace \frac{\omega_2}{2}\left( p_2^2 + q_2^2 \right) = h \right\}, \\
\mathcal{W}^s(\mathcal{M}(h)) \cap U_{q_2p_1}^+ = \left\{ (q_1, p_1, q_2, p_2) \; \vert \; \right. & \left. \kern - 
\nulldelimiterspace p_1 = -q_1, p_1(q_1, q_2, p_2; h) \geqslant 0, q_2 = 0 : \right. \nonumber \\  
& \left.\kern - \nulldelimiterspace \frac{\omega_2}{2}\left( p_2^2 + q_2^2 \right) = h \right\},
\end{align}

where each manifold represents a line parallel to the $q_1$ axis. The unstable manifold lies on the $q_1 > 0$ plane and 
stable manifold lies on the $q_1 < 0$ plane, and are marked by dashed red (unstable) and dashed white (stable) lines in the 
Fig.~\ref{fig:LD_SQH_2dof_q2p1}.

\subsubsection{Coupled quadratic Hamiltonian: 2 DoF}

\begin{figure}[!ht]
	\centering
	\subfigure[]{\includegraphics[width=0.32\textwidth]{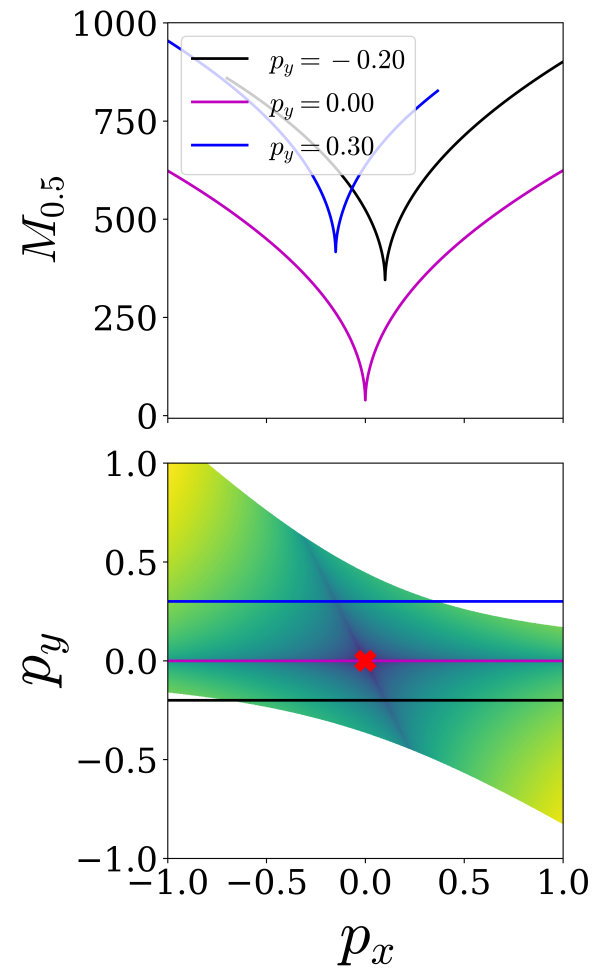}\label{fig:LD_NSQH_2dof_pxpy}}
	\subfigure[]{\includegraphics[width=0.33\textwidth]{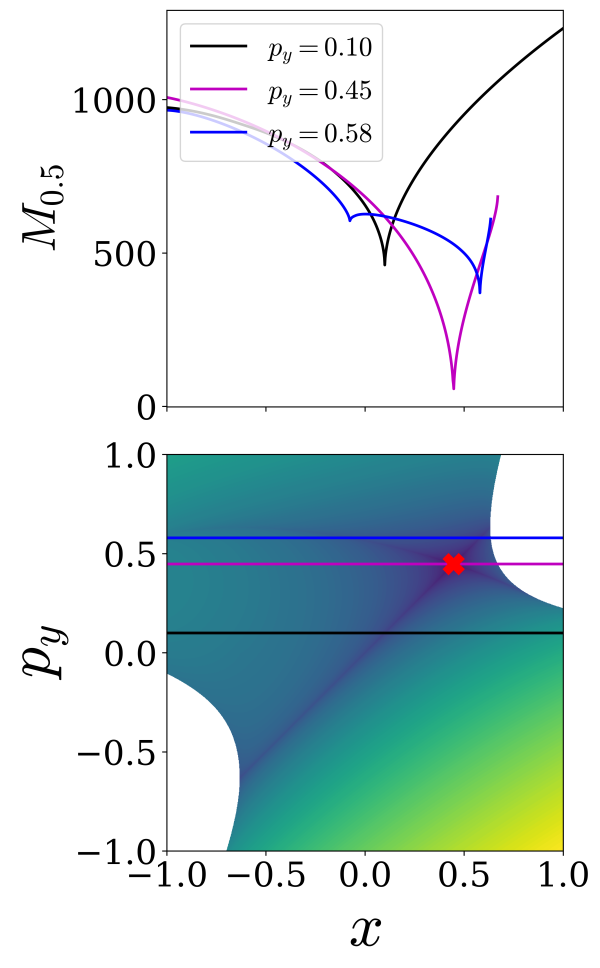}\label{fig:LD_NSQH_2dof_xpy}}
	\subfigure[]{\includegraphics[width=0.33\textwidth]{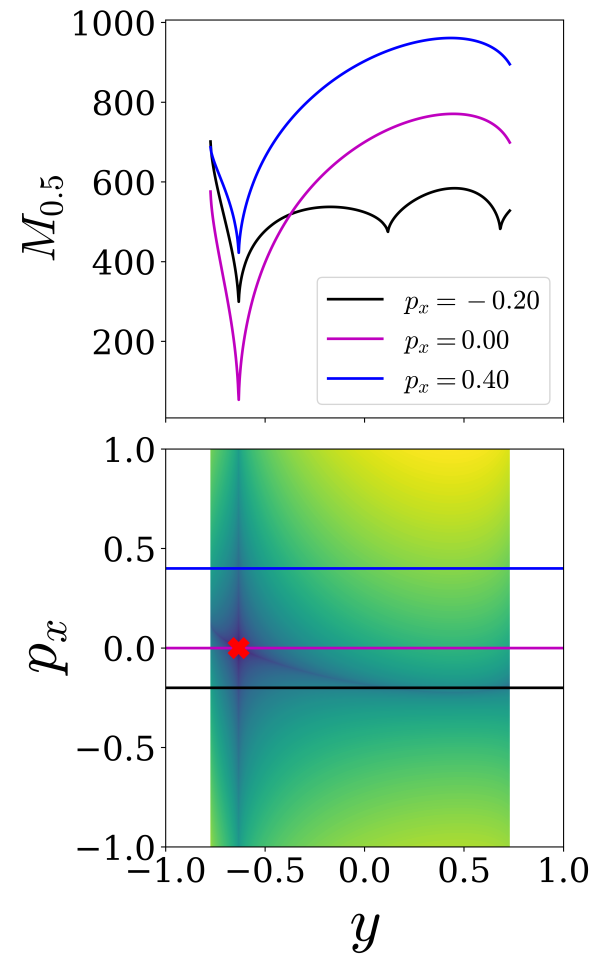}\label{fig:LD_NSQH_2dof_ypx}}
	\caption{Lagrangian descriptor plot of the non-separable quadratic Hamiltonian vector field~\eqref{eqn:eom_symp_2dof} on 
	the isoenergetic two-dimensional surface \protect\subref{fig:LD_NSQH_2dof_pxpy} $U_{p_xp_y}^{+}$, 
	\protect\subref{fig:LD_NSQH_2dof_xpy} $U_{xp_y}^{+}$, \protect\subref{fig:LD_NSQH_2dof_ypx} $U_{yp_x}^{+}$. The 
	intersection of the NHIM and the isoenergetic two-dimensional surfaces is shown as a red cross and the one corresponding 
	to the manifolds is shown as dashed red (unstable) and dashed white (stable) curves. The parameters used are $\lambda = 
	\omega_2 = 1.0$, $h = 0.2$, and $\tau = 10$.}
	\label{fig:linear_symp_trans_2dof_M400x400_2}
\end{figure}

\newpoint~Isoenergetic two-dimensional surface parametrized by $(p_x, p_y)$~\textemdash~On the fixed energy surface 
$\mathcal{H}(x, p_x, y, p_y) = h$, we compute the Lagrangian descriptor by defining a two-dimensional surface

\begin{equation}
U_{p_xp_y}^+ = \left\{ (x,p_x,y,p_y) \; | \; x = 0, \; y(x,p_x,p_y;h) > 0 : \dot{x}(x, p_x, y, p_y) > 0  \right\}
\label{eqn:2dsurf_nsqh_2dof_pxpy}
\end{equation}


Thus, the intersection of the NHIM~\eqref{eqn:nhim_symp_2dof} with this isoenergetic two-dimensional 
surface~\eqref{eqn:2dsurf_nsqh_2dof_ypy} is given by

\begin{align}
\mathcal{M}(h) \cap  U_{p_xp_y}^+ =  \left\{ (x,p_x,y,p_y) \, \vert \, \right. & \left. \kern - \nulldelimiterspace x = 0, 
p_y = x, p_x = 0, \right. \nonumber \\ 
& \left.\kern - \nulldelimiterspace \dot{x}(x, p_x, y, p_y) > 0 : y - p_x < 0, \frac{\omega_2}{2}y^2 =  h  \right\}  
\end{align}

which represents a point at the origin on the two-dimensional surface, $U_{p_xp_y}^+$, and shown as a red cross in 
Fig.~\ref{fig:LD_NSQH_2dof_pxpy}. This point is also identified by the minima in the Lagrangian descriptor values as evident 
by the one-dimensional slices at constant $p_y$.

Next, the intersection of the unstable~\eqref{eqn:umani_nsqh_2dof} and stable manifolds~\eqref{eqn:smani_nsqh_2dof} with the 
isoenergetic surface $U_{yp_y}^+$ manifest as

\begin{align}
\mathcal{W}^{\rm u}(\mathcal{M}(h)) \cap U_{p_xp_y}^+ =  \left\{ (x,p_x,y,p_y) \; \vert \; \right. & \left. \kern - 
\nulldelimiterspace x = 0, x = p_y,  \right. \nonumber \\ 
& \left.\kern - \nulldelimiterspace \dot{x}(x, p_x, y, p_y) > 0 : y - p_x < 0, \frac{\omega_2}{2}\left( - y + p_x \right)^2  
= h \right\} \\
\mathcal{W}^{\rm s}(\mathcal{M}(h)) \cap U_{p_xp_y}^+ =  \left\{ (x,p_x,y,p_y) \; \vert \;  \right. & \left. \kern - 
\nulldelimiterspace x = 0, x = 2p_x + p_y, \right. \nonumber \\ 
& \left.\kern - \nulldelimiterspace \dot{x}(x, p_x, y, p_y) > 0 : \omega_2y + ( \omega_2 + 2\lambda )p_x < 0, \right. 
\nonumber \\ 
& \left.\kern - \nulldelimiterspace \frac{\omega_2}{2}\left( \left( y + p_x \right)^2 + p_y^2 \right)  = h \right\}
\end{align}

which are shown as dashed red (unstable) and dashed white (stable) curves in the Fig.~\ref{fig:LD_NSQH_2dof_pxpy}. These 
curves are also identified by the points of minima in the Lagrangian descriptor values as evident by the one-dimensional 
sections in Fig.~\ref{fig:LD_NSQH_2dof_pxpy}.

%
%
%
%

\newpoint~Isoenergetic two-dimensional surface parametrized by $(x, p_y)$~\textemdash~On the fixed energy surface 
$\mathcal{H}(x,p_x,y,p_y) = h$, we compute the Lagrangian descriptor by defining a two-dimensional surface

\begin{equation}
U_{xp_y}^+ = \left\{ (x,p_x,y,p_y) \; | \; y = 0, \; p_x(x,y,p_y;h) > 0 : \dot{y}(x, p_x, y, p_y) > 0  
\right\}\label{eqn:2dsurf_nsqh_2dof_xpy}
\end{equation}

Thus, the intersection of the NHIM~\eqref{eqn:nhim_symp_2dof} with this isoenergetic two-dimensional 
surface~\eqref{eqn:2dsurf_nsqh_2dof_xpy} is given by

\begin{align}
\mathcal{M}(h) \cap U_{xp_y}^+ =  \left\{ (x,p_x,y,p_y) \, \vert \, \right. & \left. \kern - \nulldelimiterspace y = 0, p_y 
= x, p_x = 0, \right. \nonumber \\ 
& \left.\kern - \nulldelimiterspace \dot{y}(x, p_x, y, p_y) > 0 : (\lambda + \omega_2)p_x + 2\omega_2 x > 0, {\omega_2}x^2 
=  h  \right\}  
\end{align}

which represents a point at the origin on the two-dimensional surface, $U_{xp_y}^+$, and shown as a red cross in 
Fig.~\ref{fig:LD_NSQH_2dof_xpy}. This point is also identified by the minima in the Lagrangian descriptor values as evident 
by the one-dimensional slices at constant $p_y$.

Next, the intersection of the unstable~\eqref{eqn:umani_nsqh_2dof} and stable manifolds~\eqref{eqn:smani_nsqh_2dof} with the 
isoenergetic surface $U_{xp_y}^+$ manifest as

\begin{align}
\mathcal{W}^{\rm u}(\mathcal{M}(h)) \cap U_{xp_y}^+ =  \left\{ (x,p_x,y,p_y) \; \vert \; \right. & \left. \kern - 
\nulldelimiterspace y = 0, x = p_y,  \right. \nonumber \\
& \left.\kern - \nulldelimiterspace \dot{y}(x, p_x, y, p_y) > 0 : 2\omega_2 x + (\lambda + \omega_2)p_x > 0, \right. 
\nonumber \\ 
& \left.\kern - \nulldelimiterspace  \frac{\omega_2}{2}\left(\left( p_x + p_y \right)^2  + p_y^2 \right) = h \right\} \\
\mathcal{W}^{\rm s}(\mathcal{M}(h)) \cap U_{xp_y}^+ =  \left\{ (x,p_x,y,p_y) \; \vert \;  \right. & \left. \kern - 
\nulldelimiterspace y = 0, x = 2p_x + p_y,  \right. \nonumber \\
& \left.\kern - \nulldelimiterspace \dot{y}(x, p_x, y, p_y) > 0 : -\frac{\lambda}{2} x + (\lambda + \frac{3\omega_2}{2})p_y 
> 0, \right. \nonumber \\ 
& \left.\kern - \nulldelimiterspace  \frac{\omega_2}{2}\left(\left( p_x + p_y \right)^2  + p_y^2 \right) = h \right\}
\end{align}

which are shown as dashed red (unstable) and dashed white (stable) curves in the Fig.~\ref{fig:LD_NSQH_2dof_xpy}. These 
curves are also identified by the points of minima in the Lagrangian descriptor values as evident by the one-dimensional 
sections in Fig.~\ref{fig:LD_NSQH_2dof_xpy}.

\newpoint~Isoenergetic two-dimensional surface parametrized by $(y, p_x)$~\textemdash~On the fixed energy surface 
$\mathcal{H} = h$, we compute the Lagrangian descriptor by defining a two-dimensional surface

\begin{equation}
U_{yp_x}^+ = \left\{ (x,p_x,y,p_y) \; | \; x = 0, \; p_y(x,y,p_x;h) > 0 : \dot{x}(x, p_x, y, p_y) > 0  \right\} 
\label{eqn:2dsurf_nsqh_2dof_ypx}
\end{equation}

%
%

Thus, the intersection of NHIM~\eqref{eqn:nhim_symp_2dof} with this isoenergetic two-dimensional 
surface~\eqref{eqn:2dsurf_nsqh_2dof_ypx} is given by

\begin{align}
\mathcal{M}(h) \cap U_{yp_x}^+ =  \left\{ (x,p_x,y,p_y) \, \vert \, \right. & \left. \kern - \nulldelimiterspace x = 0, p_y 
= x, p_x = 0, \right. \nonumber \\ 
& \left.\kern - \nulldelimiterspace \dot{x}(x, p_x, y, p_y) > 0 : y < 0, \frac{\omega_2}{2}y^2 =  h  \right\}  
\end{align}

which represents a point at the origin on the two-dimensional surface, $U_{yp_x}^+$, and shown as a red cross in 
Fig.~\ref{fig:LD_NSQH_2dof_ypx}. This point is also identified by the minima in the Lagrangian descriptor values as evident 
by the one-dimensional slices at constant $p_x$.

Next, the intersection of the unstable~\eqref{eqn:umani_nsqh_2dof} and stable manifolds~\eqref{eqn:smani_nsqh_2dof} with the 
isoenergetic surface $U_{xp_y}^+$ manifest as

\begin{align}
\mathcal{W}^{\rm u}(\mathcal{M}(h)) \cap U_{yp_x}^+ =  \left\{ (x,p_x,y,p_y) \; \vert \; \right. & \left. \kern - 
\nulldelimiterspace x = 0, x = p_y,  \right. \nonumber \\
& \left.\kern - \nulldelimiterspace \dot{x}(x, p_x, y, p_y) > 0 : -y + p_x > 0 , \right. \nonumber \\ 
& \left.\kern - \nulldelimiterspace  \frac{\omega_2}{2} \left( -y + p_x \right)^2   = h \right\} \\
\mathcal{W}^{\rm s}(\mathcal{M}(h)) \cap U_{yp_x}^+ =  \left\{ (x,p_x,y,p_y) \; \vert \;  \right. & \left. \kern - 
\nulldelimiterspace x = 0, x = 2p_x + p_y,  \right. \nonumber \\
& \left.\kern - \nulldelimiterspace \dot{x}(x, p_x, y, p_y) > 0 : -\omega_2 y - (\omega_2 + 2\lambda)p_x > 0, \right. 
\nonumber \\ 
& \left.\kern - \nulldelimiterspace  \frac{\omega_2}{2}\left(\left( -y - p_x \right)^2  + 4p_x^2 \right) = h \right\}
\end{align}

which are shown as dashed red (unstable) and dashed white (stable) curves in the Fig.~\ref{fig:LD_NSQH_2dof_ypx}. These 
curves are also identified by the points of minima and singularities in the Lagrangian descriptor values as evident by the 
one-dimensional sections in Fig.~\ref{fig:LD_NSQH_2dof_ypx}.

\subsubsection{Decoupled quadratic Hamiltonian: 3 DoF}

\newpoint~Isoenergetic two-dimensional surface parametrized by $(p_2, p_3)$~\textemdash~On the constant energy surface, 
$H(q_1, p_1, q_2, p_2, q_3, p_3) = h$, we compute Lagrangian descriptor on a two-dimensional surface parametrized by $(p_2, 
p_3)$ coordinates by defining 

\begin{align}
U_{p_2p_3}^+ = \left\{ (q_1, p_1, q_2, p_2, q_3, p_3) \; | \; \right. & \left. \kern - \nulldelimiterspace p_1 = 0, q_2 = 0, 
q_3 = 0, \dot{p}_1 > 0 : \right. \nonumber \\
& \left. \kern - \nulldelimiterspace q_1(p_1, q_2, p_2, q_3, p_3; h) > 0 \right\}
\end{align}

where

\begin{align}
q_1(p_1 = 0, q_2 = 0, p_2, q_3 = 0, p_3; h) = \sqrt{\frac{2}{\lambda}\left( \left( \frac{\omega_2}{2} p_2^2 + 
\frac{\omega_3}{2} p_3^2 \right) - h \right)} 
\end{align}


The intersection of the two-dimensional surface $U_{p_2p_3}^+$ with the NHIM~\eqref{eqn:sep_quad_ham3dof_nhim} is given by

\begin{align}
\mathcal{M}(h) \cap U_{p_2p_3}^+ = \left\{ (q_1, p_1, q_2, p_2, q_3, p_3) \; \vert \; \right. & \left. \kern - 
\nulldelimiterspace q_1 = 0, p_1 = 0, q_2 = 0, q_3 = 0, \right. \nonumber \\ 
& \left.\kern - \nulldelimiterspace p_1(q_1, q_2, p_2, q_3, p_3; h) > 0, \frac{\omega_2}{2} p_2^2 + \frac{\omega_3}{2} 
p_3^2  = h \right\}.
\end{align}

which represents points on an ellipse, and marked by a dashed line in Fig.~\ref{fig:LD_SQH_3dof_p2p3}. These points are also 
identified by the minima in the Lagrangian descriptor values as shown by one-dimensional slices at constant $p_3$.

Next, the intersection of the unstable~\eqref{eqn:quad_ham3dof_umani} and stable manifolds~\eqref{eqn:quad_ham3dof_smani} 
with the isoenergetic two-dimensional surface $U_{p_2p_3}^+$ is given by

\begin{align}
\mathcal{W}^u(\mathcal{M}(h)) \cap U_{p_2p_3}^+ = \left\{ (q_1, p_1, q_2, p_2, q_3, p_3) \; \vert \; \right. & \left. \kern 
- \nulldelimiterspace q_1 = p_1, p_1 = 0, q_2 = 0, q_3 = 0, \right. \nonumber \\  
& \left.\kern - \nulldelimiterspace p_1(q_1, q_2, p_2, q_3, p_3; h) > 0, \frac{\omega_2}{2} p_2^2 + \frac{\omega_3}{2} 
p_3^2  = h \right\}, \\
\mathcal{W}^s(\mathcal{M}(h)) \cap U_{p_2p_3}^+ = \left\{ (q_1, p_1, q_2, p_2, q_3, p_3) \; \vert \; \right. & \left. \kern 
- \nulldelimiterspace q_1 = -p_1, p_1 = 0, q_2 = 0, q_3 = 0, \right. \nonumber \\  
& \left.\kern - \nulldelimiterspace p_1(q_1, q_2, p_2, q_3, p_3; h) > 0, \frac{\omega_2}{2} p_2^2 + \frac{\omega_3}{2} 
p_3^2  = h \right\},
\end{align}

where both the manifolds represent an ellipse and are marked by dashed red (unstable) and dashed white (stable) lines in the 
Fig.~\ref{fig:LD_SQH_3dof_p2p3}. These manifolds are again identified by points of minima in the Lagrangian descriptor 
values as shown by one-dimensional slices at constant $p_3$.

\newpoint~Isoenergetic two-dimensional surface parametrized by $(q_1, p_3)$~\textemdash~On the constant energy surface, 
$H(q_1, p_1, q_2, p_2, q_3, p_3) = h$, we compute Lagrangian descriptor on a two-dimensional surface parametrized by $(q_1, 
p_3)$ coordinates by defining

\begin{align}
U_{q_1p_3}^+ = \left\{ (q_1, p_1, q_2, p_2, q_3, p_3) \; | \; \right. & \left. \kern - \nulldelimiterspace p_1 = 0, q_2 = 0, 
q_3 = 0, \dot{q}_2 > 0 : \right. \nonumber \\
& \left. \kern - \nulldelimiterspace p_2(q_1, p_1, q_2, q_3, p_3; h) > 0 \right\}
\end{align}

where

\begin{align}
p_2(q_1, p_1 = 0, q_2 = 0, q_3 = 0, p_3; h) = \sqrt{\frac{2}{\omega_2}\left( h - \left( \frac{\omega_3}{2} p_3^2 - 
\frac{\lambda}{2} q_1^2 \right) \right)} 
\end{align}

The intersection of the two-dimensional surface $U_{q_1p_3}^+$ with the NHIM~\eqref{eqn:sep_quad_ham3dof_nhim} is given by

\begin{align}
\mathcal{M}(h) \cap U_{q_1p_3}^+ = \left\{ (q_1, p_1, q_2, p_2, q_3, p_3) \; \vert \; \right. & \left. \kern - 
\nulldelimiterspace q_1 = 0, p_1 = 0, q_2 = 0, q_3 = 0, \right. \nonumber \\ 
& \left.\kern - \nulldelimiterspace p_2(q_1, p_1, q_2, q_3, p_3; h) > 0,  p_3  = \pm \sqrt{ \frac{2}{\omega_3} \left(h - 
\frac{\omega_2}{2} p_2^2 \right) } \right\}.
\end{align}

which represents points on the line $q_1 = 0$, and marked by a dashed line in Fig.~\ref{fig:LD_SQH_3dof_p3q1}. These points 
are also identified by the minima in the Lagrangian descriptor values as shown by one-dimensional slices at constant $p_3$.

Next, the intersection of the unstable~\eqref{eqn:quad_ham3dof_umani} and stable manifolds~\eqref{eqn:quad_ham3dof_smani} 
with the isoenergetic two-dimensional surface $U_{q_1p_3}^+$ is given by

\begin{align}
\mathcal{W}^u(\mathcal{M}(h)) \cap U_{q_1p_3}^+ = \left\{ (q_1, p_1, q_2, p_2, q_3, p_3) \; \vert \; \right. & \left. \kern 
- \nulldelimiterspace q_1 = p_1, p_1 = 0, q_2 = 0, q_3 = 0, \right. \nonumber \\  
& \left.\kern - \nulldelimiterspace p_2(q_1, p_1, q_2, q_3, p_3; h) > 0, p_3  = \pm \sqrt{ \frac{2}{\omega_3} \left(h - 
\frac{\omega_2}{2} p_2^2 \right) } \right\}, \\
\mathcal{W}^s(\mathcal{M}(h)) \cap U_{q_1p_3}^+ = \left\{ (q_1, p_1, q_2, p_2, q_3, p_3) \; \vert \; \right. & \left. \kern 
- \nulldelimiterspace q_1 = -p_1, p_1 = 0, q_2 = 0, q_3 = 0, \right. \nonumber \\  
& \left.\kern - \nulldelimiterspace p_2(q_1, p_1, q_2, q_3, p_3; h) > 0, p_3  = \pm \sqrt{ \frac{2}{\omega_3} \left(h - 
\frac{\omega_2}{2} p_2^2 \right) } \right\},
\end{align}

where both manifolds represent points on the line $q_1 = 0$ which are marked by dashed red (unstable) and dashed white 
(stable) lines in the Fig.~\ref{fig:LD_SQH_3dof_p3q1}. These manifolds are again identified by points of minima in the 
Lagrangian descriptor values as shown by one-dimensional slices at constant $p_3$.

\begin{figure}[!ht]
	\centering
	\subfigure[]{\includegraphics[width=0.32\textwidth]{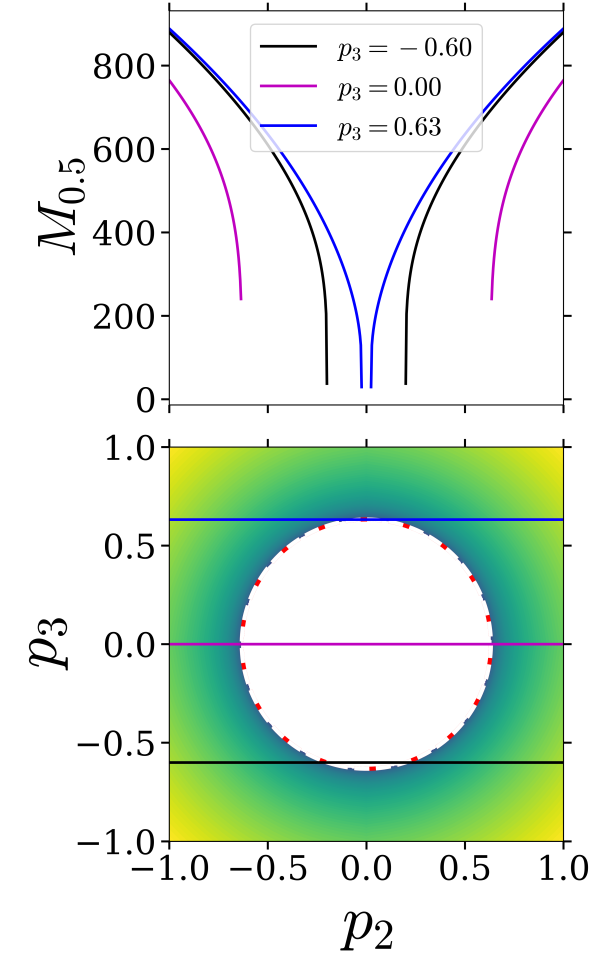}\label{fig:LD_SQH_3dof_p2p3}}
	\subfigure[]{\includegraphics[width=0.32\textwidth]{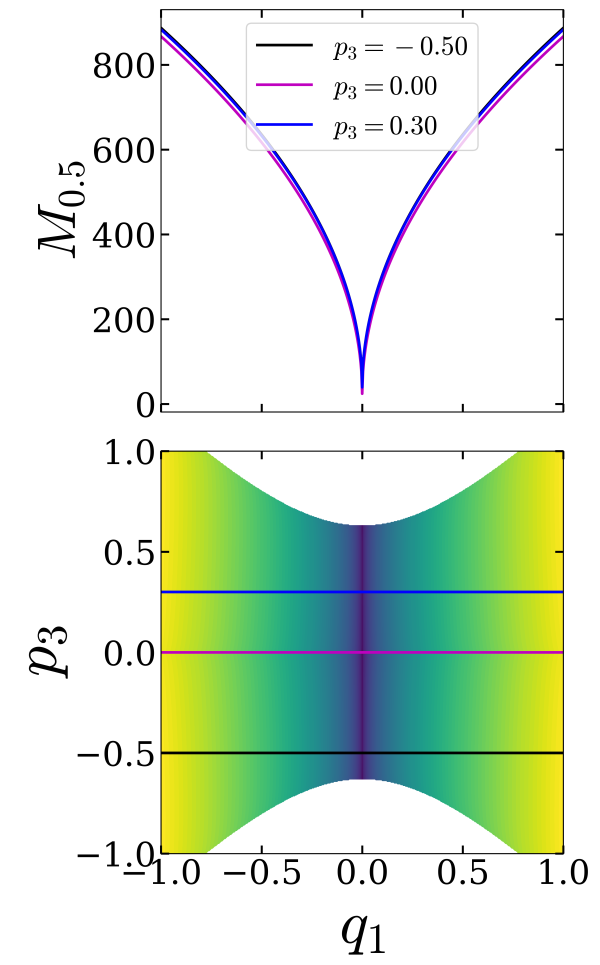}\label{fig:LD_SQH_3dof_p3q1}}
	\subfigure[]{\includegraphics[width=0.32\textwidth]{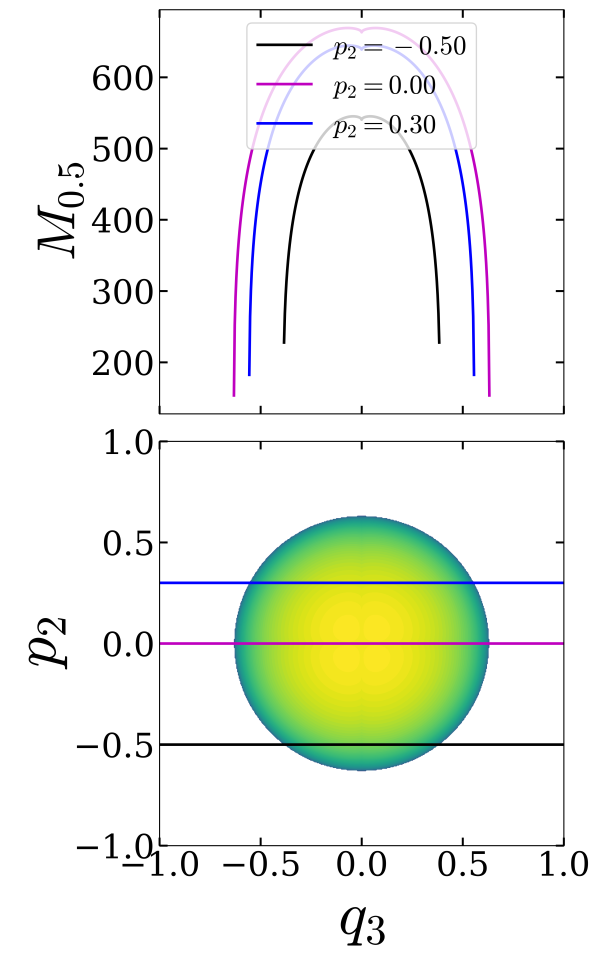}\label{fig:LD_SQH_3dof_q3p2}}
	\caption{Lagrangian descriptor plot of the separable quadratic Hamiltonian vector field~\eqref{eqn:hameq3} on the 
	isoenergetic two-dimensional surface \protect\subref{fig:LD_SQH_3dof_p2p3} $U_{p_2p_3}^{+}$, 
	\protect\subref{fig:LD_SQH_3dof_p3q1} $U_{q_1p_3}^{+}$, \protect\subref{fig:LD_SQH_3dof_q3p2} $U_{q_3p_2}^{+}$. The 
	parameters used are $\lambda = \omega_2 = \omega_3 = 1.0$, $h = 0.2$, and $\tau = 10$.}
\end{figure}

\newpoint~Isoenergetic two-dimensional surface parametrized by $(q_3, p_2)$~\textemdash~On the constant energy surface, 
$H(q_1, p_1, q_2, p_2, q_3, p_3) = h$, we compute Lagrangian descriptor on a two-dimensional surface parametrized by $(q_3, 
p_2)$ coordinates by defining 
\begin{align}
U_{q_3p_2}^+ = \left\{ (q_1, p_1, q_2, p_2, q_3, p_3) \; | \; \right. & \left. \kern - \nulldelimiterspace q_1 = 0, q_2 = 0, 
p_3 = 0, \dot{q}_1 > 0 : \right. \nonumber \\
& \left. \kern - \nulldelimiterspace p_1(q_1, q_2, p_2, q_3, p_3; h) > 0 \right\}
\end{align}
where
\begin{align}
p_1(q_1 = 0, q_2 = 0, p_2, q_3, p_3 = 0; h) = \sqrt{\frac{2}{\lambda}\left( h - \left( \frac{\omega_3}{2} q_3^2 + 
\frac{\omega_2}{2} p_2^2 \right)  \right)} 
\end{align}

The intersection of the two-dimensional surface $U_{q_3p_2}^+$ with the NHIM~\eqref{eqn:sep_quad_ham3dof_nhim} is given by

\begin{align}
\mathcal{M}(h) \cap U_{q_3p_2}^+ = \left\{ (q_1, p_1, q_2, p_2, q_3, p_3) \; \vert \; \right. & \left. \kern - 
\nulldelimiterspace q_1 = 0, p_1 = 0, q_2 = 0, p_3 = 0, \right. \nonumber \\ 
& \left.\kern - \nulldelimiterspace p_1(q_1, q_2, p_2, q_3, p_3; h) > 0,  \frac{\omega_2}{2}p_2 + \frac{\omega_3}{2}q_3^2 = 
h \right\}.
\end{align}

which represents points on an ellipse, and marked by dashed line in Fig.~\ref{fig:LD_SQH_3dof_q3p2}. These points are also 
identified by the minima in the Lagrangian descriptor values as shown by one-dimensional slices at constant $p_2$.

Next, the intersection of the unstable~\eqref{eqn:quad_ham3dof_umani} and stable manifolds~\eqref{eqn:quad_ham3dof_smani} 
with the isoenergetic two-dimensional surface $U_{q_3p_2}^+$ is given by

\begin{align}
\mathcal{W}^u(\mathcal{M}(h)) \cap U_{q_3p_2}^+ = \left\{ (q_1, p_1, q_2, p_2, q_3, p_3) \; \vert \; \right. & \left. \kern 
- \nulldelimiterspace q_1 = p_1, p_1 = 0, q_2 = 0, p_3 = 0, \right. \nonumber \\  
& \left.\kern - \nulldelimiterspace p_1(q_1, q_2, p_2, q_3, p_3; h) > 0, \frac{\omega_2}{2}p_2 + \frac{\omega_3}{2}q_3^2 = 
h  \right\}, \\
\mathcal{W}^s(\mathcal{M}(h)) \cap U_{q_3p_2}^+ = \left\{ (q_1, p_1, q_2, p_2, q_3, p_3) \; \vert \; \right. & \left. \kern 
- \nulldelimiterspace q_1 = -p_1, p_1 = 0, q_2 = 0, q_3 = 0, \right. \nonumber \\  
& \left.\kern - \nulldelimiterspace p_2(q_1, p_1, q_2, q_3, p_3; h) > 0, \frac{\omega_2}{2}p_2 + \frac{\omega_3}{2}q_3^2 = 
h  \right\},
\end{align}

where both manifolds are ellipses and are marked by dashed red (unstable) and dashed white (stable) lines in the 
Fig.~\ref{fig:LD_SQH_3dof_q3p2}. These manifolds are again identified by points of minima in the Lagrangian descriptor 
values as shown by one-dimensional slices at constant $p_2$.

\subsubsection{Coupled quadratic Hamiltonian: 3 DoF}

\begin{figure}[!ht]
	\centering
	\subfigure[]{\includegraphics[width=0.32\textwidth]{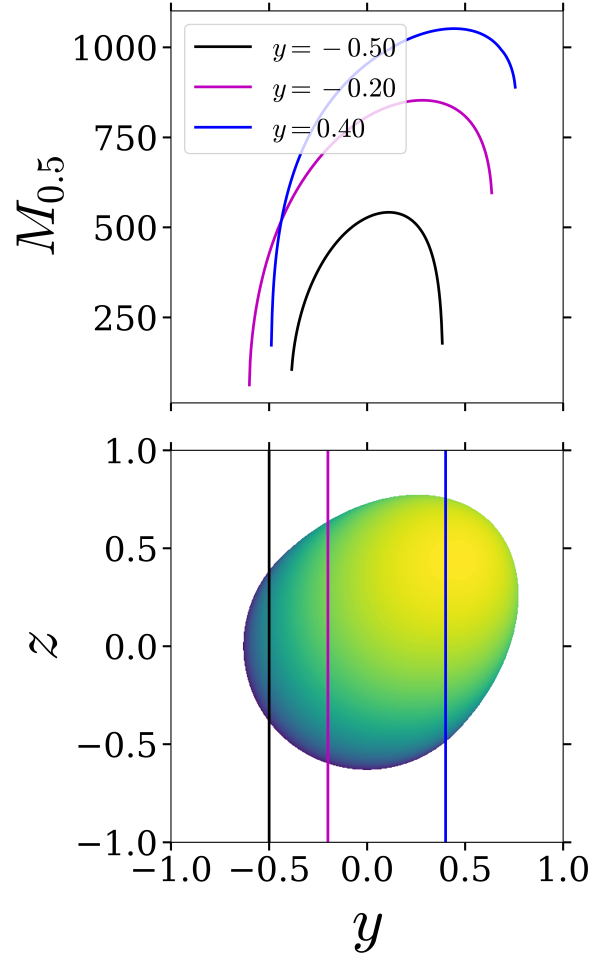}\label{fig:NSQH_3dof_yz}}
	\subfigure[]{\includegraphics[width=0.32\textwidth]{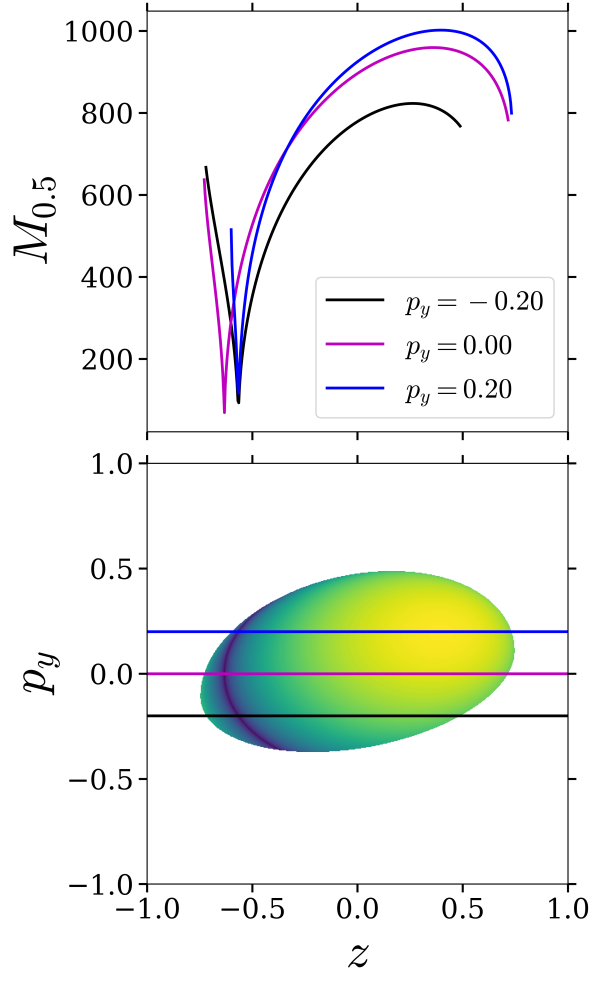}\label{fig:NSQH_3dof_zpy}}
	\subfigure[]{\includegraphics[width=0.32\textwidth]{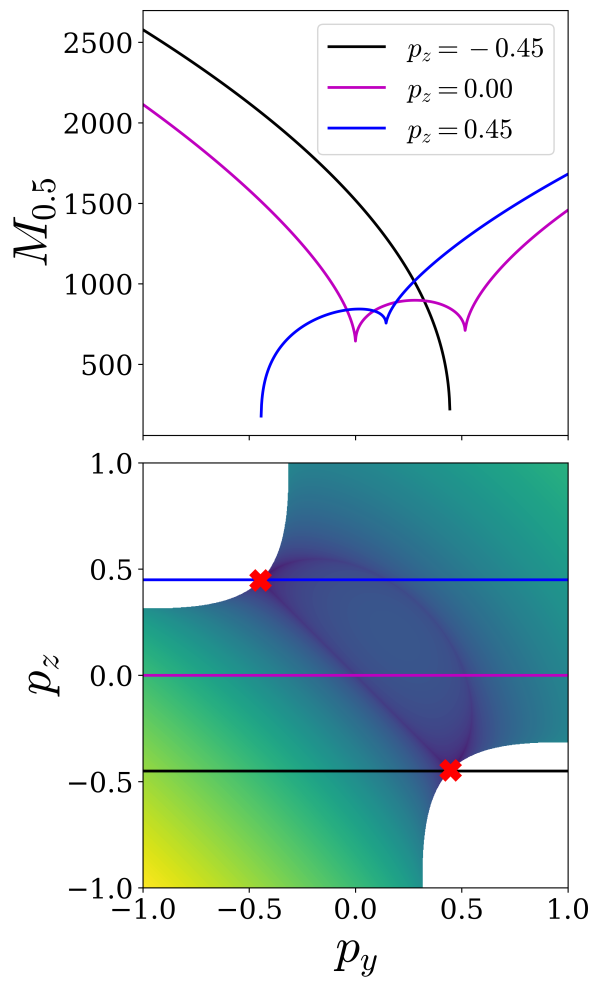}\label{fig:NSQH_3dof_pypz}}
	\caption{LD plot of the transformed Hamiltonian vector field~\eqref{eqn:eom_symp_3dof} on the 
	\protect\subref{fig:NSQH_3dof_yz} $U_{yz}^+$, \protect\subref{fig:NSQH_3dof_zpy} $U_{zp_y}^+$, 
	\protect\subref{fig:NSQH_3dof_pypz} $U_{p_yp_z}^+$. The parameters used are $\lambda = \omega_2 = \omega_3 = 1.0$, 
	integration time of $\tau = 10$, and total energy $h = 0.2$.}
	\label{fig:linear_symp_trans_3dof_M400x400_2}
\end{figure}

\newpoint~Isoenergetic two-dimensional surface parametrized by $(y,z)$~\textemdash~Next, on a constant energy surface, 
$\mathcal{H}(x, p_x, y, p_y, z, p_z) = h$, we compute the Lagrangian descriptor on a two-dimensional surface by defining 

\begin{align}
U_{yz}^+ = \left\{ (x, p_x, y, p_y, z, p_z) \; \vert \; \right. & \left. \kern - \nulldelimiterspace x = 0, p_x = 0, p_y = 
0, \right. \nonumber \\ 
& \left. \kern - \nulldelimiterspace p_z(x,p_x,y,p_y,z;h) > 0 : \dot{p_x}(x, p_x, y, p_y, z, p_z) > 0 \right\},
\label{eqn:sos_symp_3dof_yz}
\end{align} 

The intersection of the NHIM~\eqref{eqn:nhim_symp_3dof} with the two-dimensional surface is given by  

\begin{align}
\mathcal{M}(h) \cap U_{yz}^+ = \left\{ (x, p_x, y, p_y, z, p_z) \; \vert \; \right. & \left. \kern - \nulldelimiterspace x = 
0, p_x = 0, p_y = 0, p_x = 0, p_y + p_z = x, \right. \nonumber \\ 
& \left.\kern - \nulldelimiterspace \frac{\omega_2}{2} \left( p_z^2 + y^2 \right) + \frac{\omega_3}{2}\left( p_z^2 + z^2 
\right) = h, \right. \nonumber \\
& \left.\kern - \nulldelimiterspace  \dot{p_x}(x, p_x, y, p_y, z, p_z) > 0 \, : \, p_z > 0 \right\}. 
\label{eqn:nhim_zpz_3dof}
\end{align}

which is an ellipse as shown in Fig.~\ref{fig:NSQH_3dof_yz} and is also identified by the minima in the Lagrangian 
descriptor values as shown by the one-dimensional slices for constant $y$.

Next, the intersection of the unstable~\eqref{eqn:umani_nsqh_3dof} and stable~\eqref{eqn:smani_nsqh_3dof} manifolds with the 
isoenergetic two-dimensional surface~\eqref{eqn:sos_symp_3dof_ypz} becomes

\begin{align}
\mathcal{W}^{\rm u}(\mathcal{M}(h)) \cap U_{yz}^+ = \left\{  (x, p_x, y, p_y, z, p_z) \; | \; \right. & \left. \kern - 
\nulldelimiterspace x = 0, p_x = 0, p_y = 0, x = p_y + p_z, \right. \nonumber \\ 
& \left. \kern - \nulldelimiterspace \frac{\omega_2}{2} \left( -y + p_z \right)^2  + \frac{\omega_3}{2} \left( \left( -z + 
p_z \right)^2 + p_z^2 \right) = h, \right. \nonumber \\
& \left.\kern - \nulldelimiterspace  \dot{p_x}(x, p_x, y, p_y, z, p_z) > 0 : p_z > 0 \right\}, \\
\mathcal{W}^{\rm s}(\mathcal{M}(h)) \cap U_{yz}^+ = \left\{  (x, p_x, y, p_y, z, p_z) \; | \; \right. & \left. \kern - 
\nulldelimiterspace x = 0, p_x = 0, p_y = 0, x = 2p_x + p_y + p_z, \right. \nonumber \\  
& \left. \kern - \nulldelimiterspace \frac{\omega_2}{2} \left( -y + p_z \right)^2  + \frac{\omega_3}{2} \left( \left( -z + 
p_z \right)^2 + p_z^2 \right) = h, \right. \nonumber \\
& \left.\kern - \nulldelimiterspace  \dot{p_x}(x, p_x, y, p_y, z, p_z) > 0 : p_z > 0 \right\}
\end{align}

which are also be identified by the minima in the Lagrangian descriptor values as shown by the one-dimensional slice for 
constant $y$ in Fig.~\ref{fig:NSQH_3dof_yz}.

\newpoint~Isoenergetic two-dimensional surface parametrized by $(z,p_y)$~\textemdash~Next, on a constant energy surface, 
$\mathcal{H}(x, p_x, y, p_y, z, p_z) = h$~\eqref{eqn:ham_symp_3dof}, we consider a two-dimensional surface by defining

\begin{align}
U_{zp_y}^+ = \left\{ (x, p_x, y, p_y, z, p_z) \; \vert \; \right. & \left. \kern - \nulldelimiterspace x = 0, p_x = 0, y = 0 
\right. \nonumber \\ & \left. \kern - \nulldelimiterspace p_z = p_z(x,p_x,y,p_y,z;h) : \dot{y}(x, p_x, y, p_y, z, p_z) > 0 
\right\},
\label{eqn:sos_symp_3dof_zpy}
\end{align}

Thus, the intersection of the NHIM~\eqref{eqn:nhim_symp_3dof} with the two-dimensional surface is given by

\begin{align}
\mathcal{M}(h) \cap U_{zp_y}^+ = \left\{ (x, p_x, y, p_y, z, p_z) \; \vert \; \right. & \left. \kern - \nulldelimiterspace x 
= 0, p_x = 0, y = 0, p_y + p_z = x, \right. \nonumber \\ 
& \left.\kern - \nulldelimiterspace  \frac{(\omega_2 + \omega_3)}{2}p_y^2 + \frac{\omega_3}{2}z^2 = h, \right. \nonumber \\  
& \left.\kern - \nulldelimiterspace  \dot{y}(x, p_x, y, p_y, z, p_z) > 0 : \omega_2p_y - \omega_3 z > 0 \right\}. 
\label{eqn:nhim_zpy_3dof}
\end{align}

which represents a portion of an ellipse and is identified by the minima in the Lagrangian descriptor values as shown by the 
one-dimensional slices for constant $p_y$ in Fig.~\ref{fig:NSQH_3dof_zpy}.

Next, the intersection of the unstable~\eqref{eqn:umani_nsqh_3dof} and stable~\eqref{eqn:smani_nsqh_3dof} manifolds with the 
isoenergetic two-dimensional surface~\eqref{eqn:sos_symp_3dof_zpy} becomes

\begin{align}
\mathcal{W}^{\rm u}(\mathcal{M}(h)) \cap U_{zp_y}^+ = \left\{  (x, p_x, y, p_y, z, p_z) \; | \; \right. & \left. \kern - 
\nulldelimiterspace x = 0, p_x = 0, y = 0, x = p_y + p_z, \right. \nonumber \\ 
& \left. \kern - \nulldelimiterspace \frac{\omega_2}{2} p_y^2  + \frac{\omega_3}{2} \left( (x-z)^2 + p_z^2 \right) = h, 
\right. \nonumber \\
& \left.\kern - \nulldelimiterspace  \dot{y}(x, p_x, y, p_y, z, p_z) > 0 : \omega_2p_y - \omega_3 z > 0 \right\}, \\
\mathcal{W}^{\rm s}(\mathcal{M}(h)) \cap U_{zp_y}^+ = \left\{  (x, p_x, y, p_y, z, p_z) \; | \; \right. & \left. \kern - 
\nulldelimiterspace x = 0, p_x = 0, y = 0, x = 2p_x + p_y + p_z, \right. \nonumber \\  
& \left. \kern - \nulldelimiterspace \frac{\omega_2}{2} p_y^2  + \frac{\omega_3}{2} \left( (x-z)^2 + p_z^2 \right) = h, 
\right. \nonumber \\
& \left.\kern - \nulldelimiterspace  \dot{y}(x, p_x, y, p_y, z, p_z) > 0 : \omega_2p_y - \omega_3 z > 0 \right\}
\end{align}

which are identified by the minima in the Lagrangian descriptor values in Fig.~\ref{fig:NSQH_3dof_zpy} as evident by the 
one-dimensional slice at constant $p_y$.

\newpoint~Isoenergetic two-dimensional surface parametrized by $(p_y, p_z)$~\textemdash~Next, on a constant energy surface, 
$\mathcal{H}(x, p_x, y, p_y, z, p_z) = h$, we compute the Lagrangin descriptor on a two-dimensional surface by defining

\begin{align}
U_{p_yp_z}^+ = \left\{ (x, p_x, y, p_y, z, p_z) \; \vert \; \right. & \left. \kern - \nulldelimiterspace x = 0, y = 0, z = 
0, \right. \nonumber \\ & \left. \kern - \nulldelimiterspace p_x = p_x(x,p_x,y,z,p_z;h) : \dot{x}(x, p_x, y, p_y, z, p_z) > 
0 \right\},
\label{eqn:sos_symp_3dof_pypz}
\end{align}

The intersection of the NHIM~\eqref{eqn:nhim_symp_3dof} with the two-dimensional surface is given by  

\begin{align}
\mathcal{M}(h) \cap U_{p_yp_z}^+ = \left\{ (x, p_x, y, p_y, z, p_z) \; \vert \; \right. & \left. \kern - \nulldelimiterspace 
x = 0, y = 0, z = 0, p_x = 0, p_y + p_z = x, \right. \nonumber \\ 
& \left.\kern - \nulldelimiterspace \frac{\omega_2}{2} p_z^2 + \frac{\omega_3}{2}p_z^2 = h, \right. \nonumber \\
& \left.\kern - \nulldelimiterspace  \dot{x}(x, p_x, y, p_y, z, p_z) > 0 : p_x > 0 \right\}. 
\label{eqn:nhim_pypz_3dof}
\end{align}

which are two points and shown as the red crosses in Fig.~\ref{fig:NSQH_3dof_pypz}. These are also identified by the minima 
in the Lagrangian descriptor values as shown by the one-dimensional slices for constant $p_z$.

Next, the intersection of the unstable~\eqref{eqn:umani_nsqh_3dof} and stable~\eqref{eqn:smani_nsqh_3dof} manifolds with the 
isoenergetic two-dimensional surface~\eqref{eqn:sos_symp_3dof_pypz} becomes

\begin{align}
\mathcal{W}^{\rm u}(\mathcal{M}(h)) \cap U_{p_yp_z}^+ = \left\{  (x, p_x, y, p_y, z, p_z) \; | \; \right. & \left. \kern - 
\nulldelimiterspace x = 0, y = 0, z = 0, x = p_y + p_z, \right. \nonumber \\ 
& \left. \kern - \nulldelimiterspace \frac{\omega_2}{2} \left( p_x^2 + p_y^2 \right)  + \frac{\omega_3}{2} \left( p_x^2 + 
p_z^2 \right) = h, \right. \nonumber \\
& \left.\kern - \nulldelimiterspace  \dot{x}(x, p_x, y, p_y, z, p_z) > 0 : \right. \nonumber \\
& \left.\kern - \nulldelimiterspace \lambda(p_y + p_z) + (\omega_2 + \omega_3)( p_x + p_y + p_z) > 0 \right\}, \\
\mathcal{W}^{\rm s}(\mathcal{M}(h)) \cap U_{p_yp_z}^+ = \left\{  (x, p_x, y, p_y, z, p_z) \; | \; \right. & \left. \kern - 
\nulldelimiterspace x = 0, y = 0, z = 0, x = 2p_x + p_y + p_z, \right. \nonumber \\  
& \left. \kern - \nulldelimiterspace \frac{\omega_2}{2} \left( -y + p_z \right)^2  + \frac{\omega_3}{2} \left( \left( -z + 
p_z \right)^2 + p_z^2 \right) = h, \right. \nonumber \\
& \left.\kern - \nulldelimiterspace  \dot{x}(x, p_x, y, p_y, z, p_z) > 0 : \right. \nonumber \\
& \left.\kern - \nulldelimiterspace \lambda(p_y + p_z) + (\omega_2 + \omega_3)( p_x + p_y + p_z) > 0 \right\}
\end{align}

which are also identified by the minima in the Lagrangian descriptor values and shown by the one-dimensional slice for 
constant $p_z$ in Fig.~\ref{fig:NSQH_3dof_pypz}.

\end{document}